\documentclass[12pt,english]{smfart}%\documentclass[a4paper,10pt]{article}
\usepackage{smfthm}
 \usepackage{MnSymbol}
\usepackage{endnotes}
\usepackage[ethiop,english]{babel} \newcommand{\ethi}{\selectlanguage{ethiop}}  
\usepackage{hyperref}

\usepackage{epigraph}
\usepackage{upgreek}

\DeclareMathAlphabet\mathbb{U}{msb}{m}{n}

\usepackage{graphicx}
\usepackage[matrix, arrow,all,cmtip,color]{xy}
\newcommand{\bi}{\begin{itemize}}
\newcommand{\ei}{\end{itemize}}
\newcommand{\bd}{\begin{description}}
\newcommand{\ed}{\end{description}}

\newcommand{\bee}{\begin{enumerate}}
\newcommand{\eee}{\end{enumerate}}

\def\lra{\longrightarrow}
\def\ra{\rightarrow}
\def\llrra{\leftrightarrow}

\def\rtt{\,\rightthreetimes\,}

\makeatletter
\newcommand{\xleftrightarrow}[2][]{\ext@arrow 3359\leftrightarrowfill@{#1}{#2}}
\newcommand{\xdashrightarrow}[2][]{\ext@arrow 0359\rightarrowfill@@{#1}{#2}}
\newcommand{\xdashleftarrow}[2][]{\ext@arrow 3095\leftarrowfill@@{#1}{#2}}
\newcommand{\xdashleftrightarrow}[2][]{\ext@arrow 3359\leftrightarrowfill@@{#1}{#2}}
\def\rightarrowfill@@{\arrowfill@@\relax\relbar\rightarrow}
\def\leftarrowfill@@{\arrowfill@@\leftarrow\relbar\relax}
\def\leftrightarrowfill@@{\arrowfill@@\leftarrow\relbar\rightarrow}
\def\arrowfill@@#1#2#3#4{%
  $\m@th\thickmuskip0mu\medmuskip\thickmuskip\thinmuskip\thickmuskip
   \relax#4#1
   \xleaders\hbox{$#4#2$}\hfill
   #3$%
}

\makeatletter
\newcommand{\xRightarrow}[2][]{\ext@arrow 0359\Rightarrowfill@{#1}{#2}}
\newcommand{\xLeftarrow}[2][]{\ext@arrow 0359\Leftarrowfill@{#1}{#2}}
\makeatother

\makeatletter
\newcommand*{\doublerightarrow}[2]{\mathrel{
  \settowidth{\@tempdima}{$\scriptstyle#1$}
  \settowidth{\@tempdimb}{$\scriptstyle#2$}
  \ifdim\@tempdimb>\@tempdima \@tempdima=\@tempdimb\fi
  \mathop{\vcenter{
    \offinterlineskip\ialign{\hbox to\dimexpr\@tempdima+1em{##}\cr
    \rightarrowfill\cr\noalign{\kern-.3ex}
    \rightarrowfill\cr}}}\limits^{\!#1}_{\!#2}}}
\newcommand*{\triplerightarrow}[1]{\mathrel{
  \settowidth{\@tempdima}{$\scriptstyle#1$}
  \mathop{\vcenter{
    \offinterlineskip\ialign{\hbox to\dimexpr\@tempdima+1em{##}\cr
\rightarrowfill\cr\noalign{\kern-.3ex}
    \rightarrowfill\cr\noalign{\kern-.3ex}
    \rightarrowfill\cr}}}\limits^{\!#1}}}

\newcommand*{\XtoXX}[2]{\mathrel{
  \settowidth{\@tempdima}{$\scriptstyle#1$}
  \mathop{\vcenter{
    \offinterlineskip\ialign{\hbox to\dimexpr\@tempdima+1em{##}\cr
\leftarrowfill\cr\noalign{\kern-.3ex}
    \rightarrowfill\cr\noalign{\kern-.3ex}
    \leftarrowfill\cr}}}\limits^{\!#1}_{\!#2}}}

\newcommand*{\XXtoXXX}[1]{\mathrel{
  \settowidth{\@tempdima}{$\scriptstyle#1$}
  \mathop{\vcenter{
    \offinterlineskip\ialign{\hbox to\dimexpr\@tempdima+1em{##}\cr
\leftarrowfill\cr\noalign{\kern-.1ex}
    \rightarrowfill\cr\noalign{\kern-.3ex}
\leftarrowfill\cr\noalign{\kern-.3ex}
    \rightarrowfill\cr\noalign{\kern-.3ex}
    \leftarrowfill\cr}}}\limits^{\!#1}}}

\newcommand*{\XXXtoXXXX}[1]{\mathrel{
  \settowidth{\@tempdima}{$\scriptstyle#1$}
  \mathop{\vcenter{
    \offinterlineskip\ialign{\hbox to\dimexpr\@tempdima+1em{##}\cr
\leftarrowfill\cr\noalign{\kern-.3ex}
    \rightarrowfill\cr\noalign{\kern-.3ex}
\leftarrowfill\cr\noalign{\kern-.3ex}
    \rightarrowfill\cr\noalign{\kern-.3ex}
\leftarrowfill\cr\noalign{\kern-.3ex}
    \rightarrowfill\cr\noalign{\kern-.3ex}
    \leftarrowfill\cr}}}\limits^{\!#1}}}

\makeatother

\makeatother

\makeatother
\def\xra{\xrightarrow}

\def\xra{\xrightarrow}

 \def\inv{^{-1}}

\def\cof{\text{cof}}
\def\lr{\text{lr}}
\def\rl{\text{rl}}
\def\Top{\text{Top}}
\def\preorders{{\text{preorders}}}
\def\Paths{\text{Paths\,}}

\def\antidiscrete{\text{antidiscrete}}
\def\lrl{\text{l}}
\def\rlr{\text{r}}

\def\rrt#1#2#3#4#5#6{\xymatrix{ {#1} \ar[r]^{} \ar@{->}[d]_{#2} & {#4} \ar[d]^{#5} \\ {#3}  \ar[r] \ar@{-->}[ur]^{}& {#6} }}

\def\NN{\Bbb N}
\def\RR{\Bbb R}

\def\card{\,{\mathrm{card}\,}}

\def\AE{\forall\exists}

\def\Dop{\Delta^{\mathrm{op}}}

\def\Sets{\mathrm{Sets}}
\def\sSets{\mathrm{sSets}}

\def\Topp{\mathrm{Top}}

\def\LL{\mathcal L}
\def\NN{\mathbb N}

\def\id{{\text{id}}}
\def\dist{\text{dist}}
\def\diag{\text{diag}}

\def\Filt{{\ethi\ethmath{wA}}}
\def\sFilt{{{\ethi\ethmath{\raisebox{-2.39pt}{nI}\raisebox{2.39pt}\,\!\!wA}}}}

\def\sFilth{\mathit{sFilt}}

\def\Filth{{\ethi\ethmath{wA}}}
\def\sFilth{{{\ethi\ethmath{\raisebox{-2.39pt}{nI}\raisebox{2.39pt}\,\!\!wA}}}}

\def\Ob{\text{Ob\,}}

\def\sSet{\text{sSet}}
\def\diag{\text{diag}}
\def\cart{\text{cart}}

\def\const{\text{const}}

\def\FFilt{{\ethi\ethmath{wE}}}
\def\sFFilt{{\ethi\ethmath{\raisebox{-2.39pt}{nE}\raisebox{2.39pt}\,\!\!wE}}}

\def\Stone{{\ethi{\ethmath{cI}}}}

\def\ttt{{\ethi{\ethmath{pa}}}}

\def\mU{{\ethi{\ethmath{mi}}}}

\def\mmU{{\ethi{\ethmath{mA}}}}

\def\Arch{\text{Arch}}
\def\hommArrowy#1#2#3{\left\{{#2} \xrightarrow [\text{}{#1}]{} {#3}\right\}}

\def\hom#1#2{\left\{#1 \xrightarrow [\text{}]{} #2\right\}}
\def\homm#1#2#3{\left\{{#2} \xrightarrow [\text{}{#1}]{} {#3}\right\}}

\def\Hom#1#2{\left\{#1 \xRightarrow [\text{}]{} #2\right\}}
\def\Homm#1#2#3{\left\{{#2} \xRightarrow [\text{}{#1}]{} {#3}\right\}}

\def\HommArrowy#1#2#3{\left\{{#2} \xRightarrow [\text{}{#1}]{} {#3}\right\}}

\def\hommSets#1#2#3{\text{Hom}_{#1}\left({#2} , {#3}\right)}

\def\HommSets#1#2#3{\text{\underline{Hom}}_{#1}\left({#2} , {#3}\right)}
\def\hom#1#2{\text{Hom}\left(#1 ,  #2\right)}
\def\homm#1#2#3{\text{Hom}_{#1}\left({#2} , {#3}\right)}
\def\Hom#1#2{\text{\underline{Hom}}\left(#1 ,  #2\right)}
\def\Homm#1#2#3{\text{\underline{Hom}}_{#1}\left({#2} , {#3}\right)}
\title{Simplicial sets with a notion of smallness%\\ Analysis situs
}
\author{\tt 6a6ywke\thanks{\tiny
These preliminary notes are intended as an invitation to the topic, and are released in the hope of generating further activity on the subject.
\tiny\newline\tiny 
Warning: Unfortunately, %These are preliminary notes on unfinished work and
the notes  are likely to contain misprints and perhaps mistakes.  
We hope the elementary nature of the material makes them easy to ignore. 
The notes are likely to remain in current state for a while. I will be grateful for corrections 
of mistakes and inaccuracies and generally help in proofreading
 but may not be in a position to make substantial changes.
\url{mishap.sdf.org/6a6ywke/} Corrections to be sent to either  \href{https://t.me/joinchat/GVRrKxbSO8EWehZYReTKeQ}{here} or 
{\tt mi\!\!\!ishap\!\!\!p@sd\!\!\!df.org}
}} 
\begin{document}\selectlanguage{english}\catcode`\_=8\catcode`\^=7 \catcode`\_=8

\begin{abstract} We consider simplicial sets equipped with a notion of smallness,
and observe that this slight ``topological'' extension of the ``algebraic'' simplicial language
allows a concise reformulation of a number of classical notions in topology,
e.g.~continuity, limit of a map or a sequence along a filter, various notions of equicontinuity and
uniform convergence of a sequence of functions; completeness and compactness; 
in algebraic topology, locally trivial bundles as a direct product after base-change
and geometric realisation as a space of discontinuous paths.

In model theory, we observe that indiscernible sequences in a model form a simplicial set  with a notion of smallness 
which can be seen as an analogue of the Stone space of types.

These reformulations are presented as a series of exercises, to emphasise their elementary nature and
that they indeed can be used as exercises to make a student
familiar with computations in basic  simplicial and topological language. 
(Formally, we consider the category of simplicial objects in the category of filters in the sense of Bourbaki.)

This work is unfinished and is likely to remain such for a while, hence we release it as is,
in the small hope that our reformulations may provide interesting examples
of computations in basic  simplicial and topological language 
on material familiar to a student in 
in a first course of topology or category theory. 
%
%\vskip2.39pt
%These notes are intended as an invitation to the topic, and are released in the hope of generating further activity on the subject.
\end{abstract}
\maketitle

{\small\tiny 
\setcounter{tocdepth}{4}
\tableofcontents
}

\section
{Introduction}

We consider simplicial sets equipped with an additional {\em neighbourhood} or {\em being small enough} structure,
%a topological structure now 
which (paraphrasing \href{http://mishap.sdf.org/tmp/Bourbaki_General_Topology.djvu#page=15}{[Bourbaki, General Topology, \S Introduction]} on topological structure) 
%
%To formulate the idea of neighbourhood we started from the vague 
%concept of an element "sufficiently near" another element. Conversely, 
%a topological structure [on a simplicial set] now
`enables us to give a precise meaning to the
  phrase ``such and such a property holds for all [simplices sufficiently small] (orig. points sufficiently near $a$)'' :
  by definition this means that the set of [simplices] (orig. points) which have this property is
  a neighbourhood for the [neighbourhood] (orig. topological) structure in question.' 
In presence of a metric, a `small' simplex would mean a simplex being wholly contained in a small ball.

%%------The following intuition is vague but sometimes useful:  we think of a simplex $s=s_T$ as a {\em discretised homotopy} 
%%$s[1],s[2],..,s[t],..,s[T]$ between its faces, e.g.~points s[1] and s[T], 
%%and the dimension of a simplex as {\em time}. This makes it natural to require that the consecutive %(cotemporary?) 
%%faces 
%%$s[t_1\!\leq\! ... \!\leq\! t_i]$ become arbitrarily small for bounded $t=t_i-t_1$ as time/dimension $T\lra\infty$, 
%%and then treat a simplex $s_T$ (as $T\lra\infty$) as a discretised homotopy
%%%$s[1],s[2],..,s[t],..,s[T]$ 
%%between its faces. Note that the essential asymmetry (direction of time), which is apparently 
%%%in a comment by M.Kontsevich probably related to sFilt
%%a desirable property in the context of $\infty$-categories.
%%
%%

%%...%We think of a simplicial set equipped with a notion of smallness as a notion of {\em a space} which is more flexible than
%%...%the notion of a topological space. We suggest no good name for these spaces and refer to such a space
%%...%as either a simplicial neighbourhood, 
%%...%a neighbourhood structure, 
%%...%a simplicial filter; the reader preferring a short word may want to call it a situs. 
%%...%

We observe this leads to 
%think of a simplicial set equipped with a notion of smallness as 
a notion (category) of {spaces} 
%which is more flexible than the notion of a topological space: 
allowing a concise reformulation in category-theoretic terms of
a number of classical notions in topology,
e.g.~continuity, limit of a map or a sequence along a filter, various notions of equicontinuity and
uniform convergence of a sequence of functions; completeness and compactness; 
in algebraic topology, locally trivial bundles as a direct product after base-change,
and geometric realisation as a space of discontinuous paths.

The notion of smallness on a simplicial set allows to discuss topology as follows.
In a topological or metric space $X$, a point $x$ is near a point $y$ iff the simplex $(x,y)\in E\times E$ 
is small in the simplicial set of Cartesian powers of $X$ equipped with an appropriate smallness structure.. 

Formally, we consider the category of simplicial objects in the category of filters in the sense of Bourbaki, 
and the intuitive description above is formalised as follows. 
In a topological space $X$, a neighbourhood of a point $x$ consists of all points $y$ such that ``the simplex $(x,y)$ is $\varepsilon$-small'' 
for some $\varepsilon$ a neighbourhood according to the filter on $E\times E$. Yet more formally,
$U\ni x$ is a neighbourhood iff there is a subset $\varepsilon \subset E\times E$ big according to the filter on $1$-simplicies 
$E\times E$ such that $y\in U$ whenever $ (x,y)\in \varepsilon$.

For a metric space $X$, the filters on the Cartesian powers can be described explicitly as follows:
a subset $\varepsilon \subset E\times ... \times E$ is big, or a neighbourhood, iff there is $\epsilon>0$ such that a tuple belongs to $\varepsilon$
whenever it consists of points at distance at most $\epsilon$ apart.

The structure of the paper is as follows. 
In \S\ref{sec:1} we define the category $\sFilth$ and give a number of reformulations. 
Its purpose is give the reader a feeling of expressive power of $\sFilth$. 
Notably, following the approach of [Besser],[Grayson],[Drinfeld], we show\footnote{We stress, again, the preliminary nature of these notes} 
that the geometric realisation of a simplicial set 
can be interpreted as an mapping space of discontinuous paths.
In \S\ref{sec:2} we discuss the intuition behind the new notion of space. 
In \S\ref{sec:2.1} in a verbose manner we argue that the description by \href{http://mishap.sdf.org/tmp/Bourbaki_General_Topology.djvu#page=15}{[Bourbaki, Introduction]}
of the  intuition of (basic general) topology transfers to $\sFilth$ almost verbatim. 
In \S\ref{sec:2.2} we offer several vague speculations about intuition of algebraic topology. 
In \S\ref{sec:2.3} we say that $\sFilth$ has objects originating in Ramsey theory and model theory.

In \S\ref{sec:3} we mostly repeat \S\ref{sec:1} somewhat more formally and with more details.
There we present our reformulations as a series of exercises, to emphasise their elementary nature and
that they indeed can be used as exercises to make a student
familiar with computations in basic  simplicial and topological language.

We end in \S\ref{sec:4} by a discussion of Ramsey theory and indiscernibles in model theory.

\section{A sample of definitions and reformulations}\label{sec:1}

\subsection{Constructions in general topology}

\subsubsection{The definition of the category of simplicial sets with a notion of smallness}\label{filth:def} 
Here we state the key definition of the paper. In this section we introduce notation as we go; see \S\ref{notations} for explanation
if necessary. Here we just note that \href{http://mishap.sdf.org/6a6ywke/6a6ywke.pdf}{a version} 
of this paper uses type-theory notation for Hom-sets, which the author finds more readable
for nested formulae: in a category $C$, the set of all maps from an object $X$ to $Y$ is denoted by 
either $\hommArrowy C X Y $ or $\hommSets C X Y$,
and the space of maps from $X$ to $Y$, usually an object of $C$, if defined, is denoted by 
either $\HommArrowy C X Y $ or $\HommSets C X Y$.

The definition of our main category uses 
%The intuition of
the following definition  
\href{http://mishap.sdf.org/tmp/Bourbaki_General_Topology.djvu#page=63}{[Bourbaki, I\S6.1, Def.I]}
%leads to the following definition, 
which is given more formally later in Definition~\ref{def:filt}. 
\newline\includegraphics[width=\linewidth]{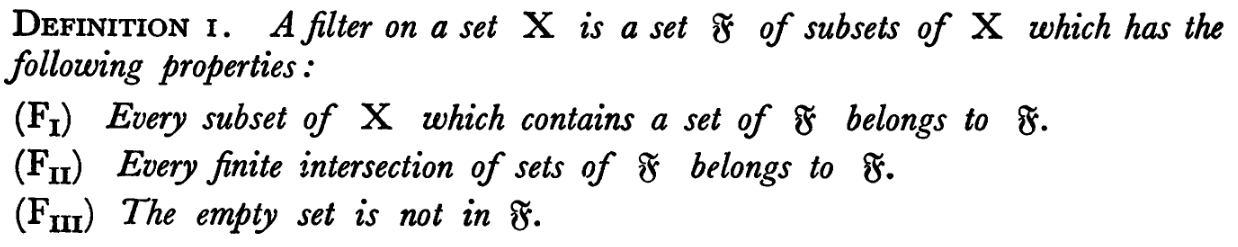}
%%A {\em filter} on a set $X$ is a collection of subsets called {\em neighbourhoods}
% topology 
%such that a subset containing a non-empty open subset is necessary open; 
% non-empty open subsets are called {\em big}. 
%%closed under finite intersection and 
%%such that a subset containing a neighbourhood is necessary a neighbourhood as well.
Subsets in $\mathfrak F$ are called {\em neighbourhoods} or {\em $\mathfrak F$-big}.
Unlike \href{http://mishap.sdf.org/tmp/Bourbaki_General_Topology.djvu}{[Bourbaki]}, we do {\em not} require $(\text{F}_\text{III})$ and allow both $X=\emptyset$ and  $\emptyset\in\frak F$; necessarily $X\in \mathfrak F$.
%to be a neighbourhood and we also allow $X$ to be empty.  

A {\em morphism of filters} is a function of underlying sets 
such that the preimage of a neighbourhood is necessarily a neighbourhood;
we call such maps of filters {\em continuous}.  
%A continuous map is a map of filters iff its image intersects each neighbourhood, hence 
%not each map sending everything to a single point is a morphism of filters. 
%Note that requiring non-emptiness condition is stronger than being continuous: a map sending everying everything to a single point 
%is necessarily continuous and is  a morphism of filters iff its image lies in every neighbourhood.

Let $\Filt$ denote the category of filters. 

\begin{defi}[Simplicial filters $\sFilt$]
Let  $\sFilt=Func(\Dop, \Filt)$ be the category of functors
from 
$\Dop$, the category opposite to the category $\Delta$ of finite linear orders,
to the category $\Filt$ of filters. We refer to its objects 
as either {\em simplicial filters}, {\em simplicial neighbourhoods}, or {\em situses},
for lack of a good name.
\end{defi}

%\subsection{---Examples of generalised spaces}.....
%%...%\subsubsection{---Synopsis: the simplicial neighbourhood associated with a metric space}
%%...%
%%...%Now we give a brief synopsis of our constructions and definitions.  
%%...%The reader is not expected to follow details but rather to get a flavour of our reformulations.
%%...%Precise definitions follow later. 
%%...%

\subsubsection{Neighbourhood structures associated with a metric space}\label{metr-sf-samples} 
See~\ref{top-metr-defs} for a precise definition.

With the set $M$ of points of a metric space associate the simplicial sets represented by $M$:
$$\Dop\lra Sets,\ n\longmapsto \homm{Sets}n M =M^n$$ 
$$\Dop\lra Sets,\ n\longmapsto \homm{Sets}{n+1} M=M\times M^n$$ 
%it represents.%, where as usual $\Dop$ denotes the category opposite to the category $\Delta$ of finite linear orders. 
%it represents.%, where as usual $\Dop$ denotes the category opposite to the category $\Delta$ of finite linear orders. 
Call a subset $\varepsilon\subset M^n=\homm{Sets}n M$  a {\em neighbourhood (of the diagonal)} iff there is $\epsilon>0$
such that $(x_1,..,x_n)\in \varepsilon$ whenever $\dist(x_i,x_j)<\epsilon$ for all $0<i<j\leq n$.
In this way we associate two simplicial neighbourhood with a metric space denoted by $M_\mU$ and $M_\mU[+1]$, resp.
There is an $\sFilt$-morphism $[-1]:M_\mU[+1]\lra M_\mU$ projecting each $M\times M^n$ to $M^n$.  

A $\sFilth$-morphism $f:L_\mU\lra M_\mU$ is a uniformly continuous map $L\lra M$.
%???????????An  $\sFilth$-morphism $f:L_\mU[+1]\lra M_\mU$ is a family $(f_a)_{a\in L}:L\lra M$ of uniformly continuous maps
%for every $\delta>0$ there is $\varepsilon>0$ such that
%$\dist( f_a(x), f_a(y)) < \delta$ whenever $\dist(a,x)< \varepsilon$ and $\dist(a,y)<\varepsilon$.

\subsubsection{Limits and $[+1]==>\id:\Dop\lra \Dop$. Uniform convergence.}\label{limits-samples}
\def\NNcof{\NN_\text{cof}} \def\NNNcof{ \left( n\mapsto \NNcof^n\right)}
Let $\NNcof$ be the filter of cofinite subsets, and let $\NNNcof$ denote the simplicial object $\Dop\lra\Filth,\ n_\leq \mapsto \NNcof^{n}$ 
of Cartesian powers of $\NNcof$.
A {\em Cauchy sequence $(a_n)_{n\in\NN}$} in a metric space $M$ is a morphism 
$\bar a:\NNNcof \lra M_\mU$, $(i_1,..,i_n)\mapsto (a_{i_1},...,a_{i_n})$:
for each $\epsilon>0$ the preimage of $\varepsilon:=\{(x,y):\dist(x,y)<\epsilon\}$ contains  $\delta:=\{(n,m):n,m>N\}\}$
for some $N$ large enough, i.e.~$\dist(a_n,a_m)<\epsilon$ for $n,m>N$.
The sequence $(a_n)_{n\in\NN}$  {\em converges} % to a point $a_\infty$} 
iff the morphism $\bar a: \NNNcof \lra M_\mU$ 
factors as  $\bar a: \NNNcof \xra{\bar a_\infty} M_\mU[+1] \xra{[-1]}  M_\mU$.
Moreover, the
 morphism ${\bar a_\infty}:  \NNNcof \xra{\bar a_\infty} M_\mU[+1]$
is necessarily of form 
$M^n\lra M^{n+1}, (i_0,i_1,...,i_n)\mapsto (a_\infty, a_{i_1},...,a_{i_n})$
where $a_\infty$ is the {\em limit} of sequence $(a_n)_{n\in\NN}$.

To see this, first note that
the underlying sset of  $M_\mU[+1]$ is 
 a disjoint union $M_\mU[+1]=\sqcup_{a\in M} \{a\}\times M_\mU$
of copies of $M_\mU$, and that the underlying sset of $ \NNNcof$ is connected.
Hence, to pick a factorisation of the underlying ssets 
is to pick an $a\in M$. Now, the map $ \NNNcof \lra \{a\}\times M_\mU$
is continuous iff for each $\epsilon>0$ the preimage of $\varepsilon:=\{(a,x):\dist(a,x)<\epsilon\}$ 
contains  $\delta:=\{(n,m):n,m>N\}$, i.e.~for $m>N$ $\dist(a,a_m)<\epsilon$.

A {\em uniformly continuous function} $f:L\lra M$ is a morphism $L_\mU\lra M$, cf.~\S\ref{ascoli}.
Indeed, for every  $\epsilon>0$ the preimage of $\varepsilon:=\{(u,v):\dist(u,v)<\epsilon\}\subset M\times M$ 
contains  ${{\updelta}} :=\{(x,y): \dist(x,y)<\delta\}\subset L\times L$ for some $\delta>0$, 
i.e.~$\dist(f_n(x),f_m(y))<\varepsilon$ whenever  $\dist(x,y)<\delta\}$.
 
A {\em uniformly equicontinuous sequence $(f_i:L\lra M)_{i\in\NN}$ of uniformly continuous functions} 
is a morphism $  \left( n\mapsto \NNcof\right) \times L_\mU \lra M_\mU$, 
$$\NN\times L^n\lra M^n,\ (i,x_1,...,x_n) \mapsto ( f_{i}(x_1),...,f_{i}(x_n)),$$
or, equivalently,  
is a morphism $  \left( n\mapsto (\NNcof^n)_\diag\right) \times L_\mU \lra M_\mU$, 
$$\NN^n\times L^n\lra M^n,\ (i_1,...,i_n,x_1,...,x_n) \mapsto ( f_{i_1}(x_1),...,f_{i_n}(x_n))$$
where  $ \left( n\mapsto (\NNcof^n)_\diag\right)$ denotes  the simplicial object $\Dop\lra\Filth,\ n_\leq \mapsto (\NNcof^{n})_\diag$
of Cartesian powers equipped with  ``the filter of cofinite diagonals'', i.e.~a subset of $\NN^n$ is a neighbourhood of $ (\NNcof^{n})_\diag$
iff it contains the set $\{ (i,i,..,i): i>N\}$ for some $N>0$.

Indeed, for every  $\epsilon>0$ the preimage of $\varepsilon:=\{(u,v):\dist(u,v)<\epsilon\}\subset M\times M$ 
contains  ${\updelta}:=\{(n,x,y):n>N,\dist(x,y)<\delta\}\subset \NN\times \NN\times L\times L $, 
resp.,  ${\updelta}:=\{(n,n,x,y):n>N,\dist(x,y)<\delta\}\subset \NN\times \NN\times L\times L $,
for some $\delta>0$ and $N>0$, 
i.e.~$\dist(f_n(x),f_m(y))<\varepsilon$ whenever $n,m>N$ and $\dist(x,y)<\delta$.

A {\em uniformly Cauchy sequence $(f_i:L\lra M)_{i\in\NN}$ of uniformly continuous functions} 
is a morphism $ \NNNcof \times L_\mU \lra M_\mU$, 
$$\NN^n\times L^n\lra M^n,\ (i_1,...,i_n,x_1,...,x_n) \mapsto ( f_{i_1}(x_1),...,f_{i_n}(x_n)).$$

Indeed, for every  $\epsilon>0$ the preimage of $\varepsilon:=\{(u,v):\dist(u,v)<\epsilon\}\subset M\times M$ 
contains  ${\updelta}:=\{(n,m,x,y):n,m>N,\dist(x,y)<\delta\}\subset \NN\times \NN\times L\times L $
for some $\delta>0$ and $N>0$, 
i.e.~$\dist(f_n(x),f_m(y))<\varepsilon$ whenever $n,m>N$ and $\dist(x,y)<\delta$.

The uniformly equicontinuous sequence $(f_i:L\lra M)_{i\in\NN}$ {\em uniformly 
converges to a uniformly continuous function $f_\infty:L\lra M$} iff this morphism ``lifts by $[-1]$'', 
i.e. fits into a commutative diagram
$$ \xymatrix{   L_\mU[+1]\times    \left( n\mapsto (\NNcof^n)_\diag\right)     \ar@{..>}[r]|-------{%(f_\infty,f_1,f_2,...)
} \ar@{->}[d]|-{[-1]\times \id} &   M_\mU[+1] \ar[d]|-{[-1]\times \id} \\  L_\mU \times   \left( n\mapsto (\NNcof^n)_\diag\right)  \ar[r]|---{(f_1,f_2,...)} &  M_\mU}$$
% $  L_\mU[+1]\times \NNNcof  \lra M_\mU[+1]\xra{[-1]} M_\mU$ 
where the top row morphism is, necessarily, of form 
$$ \NN^n\times L^{n+1}\lra M^{n+1},\ (i_1,...,i_n,x_0,x_1,...,x_n) \mapsto (\, f_\infty(x_0),  f_{i_1}(x_1),...,f_{i_n}(x_n) \,)
$$
To see this, use that $( n\mapsto (\NNcof^n)_\diag)$ is connected and therefore maps into a connected component of $M_\mU[+1]$.

\subsubsection{Complete metric spaces as a lifting property}\label{complete_as_negation}
A metric space $M$ is {\em complete} iff every Cauchy sequence converges, i.e.~the following lifting property\footnote{
%By $f\rtt g$ we denote that a morphism $f$ has the Quillen lifting property with respect to a morphism $g$.
A morphism $i:A\to B$ in a category has {\em the left lifting property} with respect to a
morphism $p:X\to Y$, and $p:X\to Y$ also has {\em the right lifting property} with respect to $i:A\to B$,        
denoted  $i\rtt p$,
 iff for each 
       $f:A\to X$ and                                                                                                                           
       $g:B\to Y$
        such that                                                                                                                     
$p\circ f=g\circ i$
there exists                                                                                                         
       ${ h:B\to X}$ such that                                                                                                                     
       ${h\circ i=f}$ and                                                                                                                         
       ${p\circ h=g}$.
This notion is used to define properties of morphisms starting from an explicitly given class of morphisms,
often a list of (counter)examples,
and a useful intuition is to think that the property of left-lifting against a class $C$ is a kind of negation
of the property of being in $C$, and that right-lifting is also a kind of negation.
} holds: 
$$\emptyset\lra \NNNcof \rtt M_\mU[+1]\lra M_\mU$$
or, in another notation,\footnote{Denote by $P^\lrl$ and $P^\rlr$ the classes (properties) of morphisms having the left, resp.~right, lifting property
with respect to all morphisms with property $P$:
$$P^\lrl:=\{ f\rtt g: g\in P\} \ \ \ \ P^\rlr := \{ f\rtt g: f\in P\}$$
It is convenient to refer to $P^\lrl$ and $P^\rlr$ as the property of {\em left, resp.~right, Quillen negation of property $P$}.
} 
$$\emptyset\lra \NNNcof \in \{ M_\mU[+1]\lra M_\mU \,:\, M\text{ is a complete metric space} \}^\lrl$$
Hence, $\emptyset\lra \NNNcof \in  \{ \RR_\mU[+1]\lra \RR_\mU  \}^\lrl$, and therefore 
\begin{quote} 
 $M_\mU[+1]\lra M_\mU\in  \{ \RR_\mU[+1]\lra \RR_\mU  \}^\lr$ implies $M$ is complete, 
\end{quote} and a little argument shows the converse holds for precompact metric spaces. 

Compactness can also be reformulated as a lifting property, see \S\ref{Compactness-as-lrl} 
for this and other examples. 

\subsection{Elementary constructions in homotopy theory
}

\subsubsection{The unit interval}\label{unit-interval-sample} With the unit interval $[0,1]$ associate the simplicial set  
$$\Dop\lra Sets,\ n_\leq\longmapsto \homm{preorders}{n_\leq}{[0,1]_\leq}$$
Equip it with a neighbourhood structure using the metric: $\varepsilon\subset \homm{preorders}{n_\leq}{[0,1]_\leq}$
is a neighbourhood iff there is $\epsilon>0$ such that 
$(t_1\leq ...\leq t_n)\in\varepsilon$ whenever $t_n < t_1+\epsilon$. This neighbourhood structure can be defined
entirely in terms of the simplicial set itself, cf.~\S\ref{unit:sF} for details:
 $\varepsilon\subset \homm{preorders}{n_\leq}{[0,1]_\leq}$
is a neighbourhood iff for any $\tau>0$ there is $T>\tau>n$ and a simplex $s=(s_1\!\leq\!..\!\leq\!s_T)$ such that
$t[{i_1}\!\leq\!..\!\leq\!{i_n}]=(t_{i_1}\!\leq\!..\!\leq\!t_{i_n})\in\varepsilon$
whenever $T'>0$, the simplex $s$ is a face of a simplex $t=(t_1\!<\!t_2\!<\!..\!<t_{T'})$ and $i_1\leq ..\leq i_n<i_1+\tau$. 
 Denote this simplicial neighbourhood by $[0.1]_\leq$. 

A {\em path} $\gamma:[0,1]\lra M$ in a metric space $M$ is same as a morphism $[0,1]_\leq \lra M_\mU$.  
An automorphism $[0,1]_\leq\lra[0,1]_\leq$ is a non-decreasing (necessarily uniformly) continuous automorphism $[0,1]\lra[0,1]$
of the unit interval.

%TODO: any other construction wrt the unit interval can i give here? 

\subsubsection{Simplicies as $\varepsilon$-discretised homotopies.}\label{disc-homotopies-sample}
A map $f$ is homotopic to a map $g$ iff there is a sequence $f=f_0$, $f_1$, ..., $f_t$,...,$f_T=g$ 
where $f_t$ is as near as we please to $f_{t+1}$, $0\!\leq\!t\!\leq\!t+1\!\leq\!T$.
In $\sFilt$ this is readily formalised by saying that 
the simplex $(f,g)$ is a face of simplex $\overrightarrow f$ with 
 consecutive faces %$f[i,i+1]$ 
as small as we please, 
i.e.~for each neighbourhood $\varepsilon$ in the set of $1$-simplices 
there is a simplex $\overrightarrow f=\overrightarrow{f_\varepsilon}$ %, $T=T_\varepsilon=\dim \overrightarrow f$ 
 in the space  of maps from $X$ to $Y$
such that $(f,g)=\overrightarrow f[0\!\leq\! T]$ for some $t\geq 0$, 
and $\overrightarrow f[t \leq t+1]\in\varepsilon$ for $0\!\leq\!t<T=\dim \overrightarrow f$;
here $\overrightarrow f[0\!\leq\! T]$ denotes the face of simplex $\overrightarrow f$
corresponding to $\Delta$-morphism $0\!\leq\! T:2_\leq\lra T_\leq$, $0\mapsto 0$, $1\mapsto T$,
and similarly for $\overrightarrow f[t \leq t+1]$. 
% lie in the neighbourhood $\varepsilon$. 
This formalisation immediately suggests
%that there is no reason to restrict to $1$-simplices only and 
we should let $\varepsilon$ vary among neighbourhoods of arbitrary dimension $T'$
and rather require that $\overrightarrow f=\overrightarrow {f_{\varepsilon,n}}$ and
$\overrightarrow f[t_0\!\leq\!t_1\!\leq\!...\!\!\leq\!t_{T'}]\in\varepsilon$ whenever  
$t_{T'}\!\leq\!t_0+n$ (where $\varepsilon\subset X_{T'}$ and $n>0$). %See \S\ref{discr_homotpies} for details. 

This leads to the following definition.

For a neighbourhood $\varepsilon\subset M_T$ and $n>0$, a simplex $s:M_{T'}$ is {\em $\varepsilon/n$-fine} 
iff $s[t_0\leq ...\leq t_k]\in \varepsilon$ whenever $0\leq t_0\leq ... \leq t_k\leq t_0+n\leq T'$.
A simplex $s$ is {\em Archimedean} iff it can be split into finitely many arbitrarily small parts, 
i.e.~is a face of some $\varepsilon/n$-fine simplex for every neighbourhood $\varepsilon\subset X_k$ and every $T,n>0$.\footnote
{This definition applies to any object of $\sFilt$ but should likely be modified even for the metric spaces. 
For explanation see the footnote in \S\ref{discr_homotpies}.}
For example, a pair of points $(x,y)\in M\times M$ in a metric space $M$ is an Archimedean simplex in $M_\mU$ 
iff for each $\epsilon>0$ there is 
an $\epsilon$-{\em discretised homotopy} 
$x=x_0,x_1,...,x_l=y$, $\dist(x_t,x_{t+1})<\epsilon$ for $0\leq t<l$, {\em from $x=(x,y)[0]$ to $y=(x,y)[1]$}. 
%In other words, 
%an Archimedean simplex in $M\times M$ is a pair of homotopic points, or rather a pair of points
%which can be connected by a {\em discretised homotopy} of arbitrary finess.

%Note that if $M$ is the space of functions, 

Archimedean simplices of a simplicial filter $X:\sFilth$ 
form a subobject (subfunctor) $X_{\Arch}$, as the definition is invariant. 
%Two functions $f,g:A\lra M$ are homotopic iff $(f,g)$ is an Archimedean simplex 

A well-known lemma says that two functions $f,g:A\lra M$ from an arbitrary topological space $A$ to a metric space $M$
are homotopic iff there is a $\epsilon$-discretised homotopy 
$f=f_2,..., f_n=g$ such that for any $x\in A$ $\dist(f_t(x), f_{t+1}(x))<\epsilon$, under some assumptions on 
the metric space $M$; it is enough to assume that for every $\epsilon>0$ there is $\delta>0$ such that 
every $\epsilon$-ball contains a contractible $\delta$-ball. 
We reformulate this by saying that  
two functions $f,g:A\lra M$ are homotopic iff $(f,g)$ is an Archimedean simplex of the mapping space $Func(A,M)_\mU$ 
with the sup-metric, or, equivalently, %we may say that 
%a homotopy is an Archimedean $1$-simplex of the mapping space 
a $1$-simplex of $(Func(A,M)_\mU)_{\Arch}$.

\subsubsection{Topological spaces as simplicial filters}\label{top-sf-samples} 
See \S\ref{top-sf-intuit} for the intuition and \S\ref{top-metr-defs} for a precise definition.

As with metric spaces, with the set $X$ of points of a topological space associate the simplicial set 
$$\Dop\lra Sets,\ n_\leq \longmapsto \homm{Sets}n X =X^n$$ 
%it represents.%, where as usual $\Dop$ denotes the category opposite to the category $\Delta$ of finite linear orders. 
%it represents.%, where as usual $\Dop$ denotes the category opposite to the category $\Delta$ of finite linear orders. 
Define the filters of  neighbourhoods (of the diagonal) as follows. 
The filter on  $X$ is antidiscrete, as $X$ is the diagonal of itself and thus every neighbourhood has to contain $X$.
A subset  $\varepsilon\subset X\times X$ is a neighbourhood iff $\varepsilon$ contains a set of the form 
$$\bigcup_{x\in X} \{x\}\times U_x
$$
where $U_x\ni x$ is an open neighbourhood of $x$. 
The filter on $X^n$, $n>2$, is the coarsest filter compatible with all the face maps $X^n\lra X\times X$. 
Let $X_\ttt$ denote the simplicial neighbourhood obtained in this way. 

\subsubsection{A forgetful functor to topological spaces}\label{sFtoTop}
The embedding of topological spaces admits an inverse $-_{\ttt\inv}:\sFilth\lra\Topp$ 
defined similarly to the definition of the topology associated with a uniform structure  \href{http://mishap.sdf.org/tmp/Bourbaki_General_Topology.djvu#page=180}{[Bourbaki,II\S1.2,Prop.1,Def.3]},
as follows.

The set of points of $X_{\ttt\inv}$ is the set of points which are $\varepsilon$-small for each neighbourhood $\varepsilon\subset X_0$, i.e.
$X_\text{points}:=\bigcup\limits_{\varepsilon\subset X_0\text{ is a neighbourhood}} \varepsilon$.
The topology is generated by the subsets that together with each point contain all $\varepsilon$-near points
for some $\varepsilon\subset X_1$, i.e.~subsets $U$ with the following property:
there is a neighbourhood $\varepsilon\subset X_1$ such that %neighbourhood base consisting of the following subsets: 
for every point $x\in \varepsilon$ such that both $x[0],x[1]\in X_\text{points}$, it holds that $x[0]\in U$ implies $x[1]\in U$.  

It is easy to check that for a topological space $X$, $(X_\ttt)_{\ttt\inv}=X$, and that 
$([0.1]_\leq)_{\ttt\inv}=[0,1]$ as a topological space.

\subsubsection{Locally trivial bundles}\label{def-bundle-sample} Let $X,B,F$ be %metrisable 
topological spaces. 
A map  $X\xra p B$ is locally trivial iff it becomes a direct product after pull back to the ``local base''  $[-1]:B_\ttt[+1]\lra B_\ttt$, 
i.e.~it fits into the following commutative diagram:
$$\xymatrix{ B_\ttt[+1]\times_{B_\ttt} F_\ttt \ar[rd]  \ar@{..>}[r]|--{(iso)} & B_\ttt[+1]\times_{B_\ttt} X_\ttt  \ar[r] \ar[d] & X_\ttt \ar[d]|{p}\\
{} & B_\ttt[+1] \ar[r] & B 
}$$
That is, a map $X\xra p B$ is {\em locally trivial} iff there is an $\sFilt$-isomorphism $B_\ttt[+1]\times_{B_\ttt} F_\ttt \xra{(iso)} B_\ttt[+1]\times_{B_\ttt} X_\ttt$ over $B_\ttt$. 

Let us verify that this diagram represents the usual definition of local triviality. 
To give a morphism of sSets  $$B_\ttt[+1]\times_{B_\ttt}
X_\ttt=\bigsqcup\limits_{b\in B} \{b\}\times X_\ttt \xra{(iso)} B_\ttt[+1]\times
F_\ttt=\bigsqcup\limits_{b\in B} \{b\}\times B_\ttt\times F_\ttt$$ over $B_\ttt$
is to give for each $b\in B$ a morphism $f_b:X\lra B\times F$; to check this, use that 
ssets $X_\ttt$ and $B_\ttt\times F_\ttt$ are connected. 
The morphism of ssets is an isomorphism of ssets iff each $f_b$ is an isomorphism of sets, i.e.~a  bijection.

Let us now prove that each $f_b$ is a homeomorphism with a neighbourhood of form $U_b\times F$ where $b\in U_b\subset B$ is open.

For each $(b',y')\in B\times F$ pick a neighbourhood $W_{(b',y')}\subset B\times B\ni (b',y')$ 
which is a counterexample to continuity of $f_b$ at the unique preimage of $(b',y')$ 
if it is indeed not continuous at that point. 
The following is a neighbourhood at the set of $2$-simplicies of $B_\ttt[+1]\times F_\ttt$: 
$$\varepsilon:=
\left\{ (b,\, (b',y'),\, (b'',y'')\,) % \in B\times (B\times F)\times (B\times F)  
 \, : \, (b',y')\in B\times F, \, (b'',y'')\in W_{(b',y')}
\right\}
\cup \bigsqcup\limits_{b'\neq b} \{b'\}\times (B\times F)\times (B\times F)$$
By continuity its preimage $\delta_b:=f_b^{\inv}(\varepsilon)$ contains
a set of the form 
$$\{ (b, x', x'')\in B\times X \times X  \,:\, p(x')\in U_b,\ x''\in V_{x'} \} $$
where $U_b\ni b, V_{x'}\ni x'$ are open. %neighbourhoods of these points. 
Hence $f(V_{x'})\subset W_{f_b(x')}$ for all $x'$ such that $p(x')\in U_b\subset B$, 
and, by choice of the neighbourhoods $W_{(b',x')}$, the function $f_b$ is continuous 
over the preimage of $U_b\times F \subset B\times F$.  

A similar argument establishes continuity of $f_b^{\inv}$. 

\subsection{Geometric realisation as path mapping spaces: the approach of Besser, Grayson and Drinfeld}

[Grayson, Remark 2.4.1-2] interprets the geometric realisation of a simplex $\Delta_{N}$
as a space of non-decreasing maps $[0,1]\lra (N+1)_\leq$.
$$|\Delta_{N}|=\{(s_1,..,s_N)\in \RR^N: 0\leq s_1\leq ... \leq s_N\leq 1\} \approx 
\{ s\!:\![0,1]_\leq \lra (N+1)_\leq \} $$ $$
 0\leq s_1\leq ... \leq s_N\leq 1 \ \approx (  [0,s_1)\mapsto 0, ..., [s_{N-1},s_N)\mapsto N-1, [s_N,1]\mapsto N )
$$
with a metric analogous to Levi-Prokhorov or Skorokhod metrics on the spaces of discontinuous functions used in probability theory;
roughly, two functions are close in such a metric iff one can be obtained from the other by a small perturbation of both values and arguments;
in other words, a small neighbourhood of the graph of one function contains the other one.

 We use this observation and 
the construction of geometric realisation by [Drinfeld] to define, for a simplicial set $X$, 
a $\sFilth$-structure on the inner Hom in sSets $$\Homm\sSets {\homm{{\text{preorders}}}-{[0,1]_\leq} } X $$
analogous to the Skorokhod metric. % in probability theory. 
 We then argue
that the metric space associated with this $\sFilth$-object is the geometric realisation of $X$, under some assumptions.
%We assume the reader is familiar with [Grayson, \S2.4] and [Drinfeld]. 

\subsubsection{Drinfeld construction of geometric realisation as a space of paths with %Levi-Prokhorov or 
Skorokhod metric}\label{Skorokhod-Grayson}
As a warm-up, for the reader  familiar with [Grayson, \S2.4] and [Drinfeld],
we sketch a construction of a metric on $$\homm\sSets {\homm{{\text{preorders}}}-{[0,1]_\leq} } X $$
analogous to the Skorokhod metric in probability theory.

A finite subset $F\subset [0,1]$ and an $x\in X(\pi_0([0,1]\setminus F))$ determines a morphism of sSets 
$ \homm{\text{preorders}}{-}{[0,1]_\leq} \lra  X $ as follows:
$$ \homm{\text{preorders}}{n_\leq}{[0,1]_\leq} \lra  X(n_\leq) $$
%$$\overrightarrow t \longmapsto x\left[\xra{\overrightarrow t\,:\,n_\leq\lra [0,1]}\cdot\xra{[0,1]\lra \pi_0([0,1]\setminus F)}\right]\in X(n_\leq)
%$$
$$\overrightarrow t \longmapsto x\left[\xra{n_\leq\xra{\overrightarrow t} [0,1]\lra \pi_0([0,1]\setminus F)}\right]\in X(n_\leq)
$$
where $[0,1]\lra \pi_0([0,1]\setminus F)$ is the obvious map contracting the connected components (we need to make
a convention where to send points of $F$).    

A verification shows that this defines,  
%indeed this defines a morphism of sSets
%$$\homm{\text{preorders}}- {[0,1]_\leq} \lra  X ,$$ and, 
moreover, a map of sets
$$
|X|:=\varinjlim\limits_{{ F\subset [0,1] \text{ finite}}} X(\pi_0([0,1]\setminus F)) 
\lra \homm{sSets}{\homm{\text{preorders}}{-}{[0,1]_\leq}}{ X } 
$$

Conversely, a map $\pi:\homm{\text{preorders}}- {[0,1]_\leq}\lra X$ of ssets determines a system of points  as follows:
%$$x_\theta\in |X|,\theta:n_\leq\lra [0,1],n>0$$ %of points 
$$\pi\left(\xra{\theta\,:\,n_\leq\lra[0,1]}\right)\in X\left([0,1]\setminus\{\theta(0),..,\theta(n-1)\}\right)$$
and thereby a point of $|X|$.

Define the following pseudometric analogous to  {\em Levi-Prokhorov or Skorokhod}\footnote{
This definition is similar to the definition of Skorokhod metric as phrased by [Kolmogorov, \S2, Def.1 of $\varepsilon$-equivalence],
particularly if we consider functions taking values in a discrete metric space $0,1,..,N$, $\rho(n,m):=|n-m|$:
two functions $f,g$ are called $\varepsilon$-equivalent iff there exists $r$ and
\includegraphics[width=\linewidth]{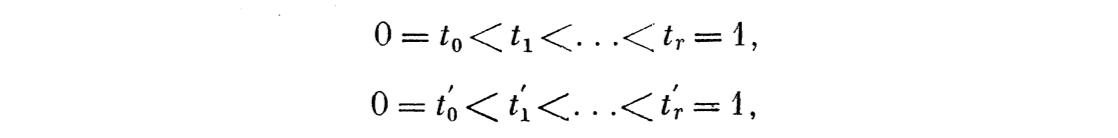}
such that for $k=1,...,r$ the following inequalities hold:
\includegraphics[width=\linewidth]{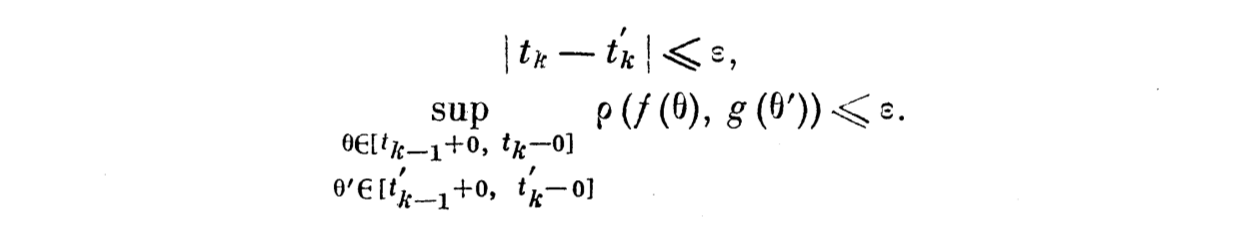}
The original goal of the definition was to define a distance or convergence 
for (distributions of) stochastic processes such that 
a small distortion of either timings of events
or their values results in a small distance.
}
metric
(here we allow distance to be 0):
$$ f, g: \homm{\text{preorders}}{-}{[0,1]_\leq} \lra  X $$
$$\dist(f, g):=\inf\{\epsilon>0 \,:\, \forall n>0 \,
\forall \overrightarrow t=(t_1\leq ...\leq t_n)  \, \exists \overrightarrow t'=(t'_1\leq ...\leq t'_n) \,\,
$$ $$ 
\left(\, f(\overrightarrow t)= g(\overrightarrow t') \,\,\&\,\, |t_1-t'_1|<\epsilon \,\,\,\&\,\,\, ... \,\,\, \& \,\,\,|t_n-t'_n|<\epsilon \,\right)
$$

Let us now compare this construction with  [Grayson, Remark 2.4.1-2] 
for $X=\homm{\text{preorders}}{\cdot}{N_\leq}$ the $(N-1)$-simplex $\Delta_{N-1}$. In this case the Yoneda lemma gives us that a  map 
$[0,1]_\leq = \homm{\text{preorders}}{\cdot}{ [0,1]_\leq}\lra \homm{\text{preorders}}{\cdot}{N_\leq}$ is the same as a map $[0,1]_\leq \lra N_\leq$, 
which, in turn, is essentially the same as a sequence $0\leq s_1\leq ... \leq s_{N-1}\leq  1$, i.e.
a point of the geometric realisation 
$$ |\Delta_{N-1}|= \{ (s_1,..,s_{N-1})\,:\:  0\leq s_1\leq ... \leq s_{N-1}\leq  1\}\subset \RR^{N-1} .$$
We also see that the metric coincides with the metric defined by [Grayson, Remark 2.4.1-2].

\subsubsection{The Skorokhod filter on a Hom-set}\label{skorokhod-filter-def}\label{Skorokhod-mapping-spaces-sample}
For $N>2n$, $\delta\subset X_N$ and $\varepsilon\subset Y_n$, 
a {\em $\varepsilon\delta$-Skorokhod neighbourhood of 
 Hom-set $\homm\sSets  {X_\text{as sSet}} {Y_\text{as sSet}}$} of the underlying simplicial sets
of $X$ and $Y$
 is the subset %of $\Hom X Y$ 
consisting 
of all the function $f:X_\text{as sSet}\lra Y_\text{as sSet}$ with the following property:
\bi\item[] there is a neighbourhood $\delta_0\subset X_n$ such that 
each ``$\delta_0$-small'' $x\in \delta_0$ has a ``$\delta$-small''
``continuation'' $x'\in X_N$, $x=x'[1..N]$ such that its ``tail'' maps into something ``$\varepsilon
$-small'', 
i.e.~$f(x'[N-n+1..N])\in\varepsilon$.
\item [] As a formula, this is 
$$\{ f:X\lra Y\,:\, \exists \delta_0 \subset X_n\, \forall x\in\delta_0 \, \exists x'\in \delta
( x=x'[1...n] \,\,\,\&\,\,\, f(x'[N-n+1,...,N])\in \varepsilon ) \}$$
\ei

The {\em Skorokhod filter on  $\homm\sSets  {X_\text{as sSet}} {Y_\text{as sSet}}$} 
is the filter generated by all the Skorokhod $\varepsilon\delta$-neighbourhoods for 
$N\geq 2n>0$ (sic!), neighbourhoods $\delta\subset X_N$ and $\varepsilon\subset Y_n$.

Let $\Homm\Filt X Y$ denote $\homm\sSets  {X_\text{as sSet}} {Y_\text{as sSet}}$ equipped with the Skorokhod neighbourhood structure.
This allows to  define mapping spaces in $\sFilt$ by equipping the inner Hom of ssets with the Skorokhod filters.
%\begin{defi}[Skorokhod mapping space]ooo The {\em Skorokhod mapping space $\Homm\sFilth X Y$} is the simplicial set $\Homm\sSets XY$ equipped 
\begin{defi}[Mapping space]
 The {\em Skorokhod mapping space $\Homm\sFilth X Y$} is the inner Hom  $\Homm\sSets  {X_\text{as sSet}} {Y_\text{as sSet}}$ 
 of the underlying simplicial sets
of $X$ and $Y$
equipped
with the neighbourhood structure as follows. Equip $\homm{{\text{preorders}}}-{n_\leq}$ with the antidiscrete filter, 
equip  $X\times \homm{{\text{preorders}}}-{n_\leq}$ with the product filter, and, finally, equip the set of $(n-1)$-simplicies  
$\homm\sSets{{X_\text{as sSet}\times \homm{{\text{preorders}}}-{n_\leq}}} {Y_\text{as sSet}}$ with the resulting Skorokhod neighbourhood structure. %filter. 
\end{defi}

\subsubsection{The geometric realisation of a simplex and its Skorokhod space of paths}\label{Skorokhod-paths-space}
We now rephrase [Grayson, Remark 2.4.1-2]  in terms of $\sFilth$.
Let $(\Delta_N)_\diag$ denote the standard simplex $\Delta_N= \homm\sSets-{{N+1}_\leq}$
equipped with the {\em filter of diagonals}, i.e.~the filter on $(\Delta_N)_0$
is antidiscrete and for $n>0$ the filter on $(\Delta_N)_n$ is
the coarsest filter such that the diagonal degeneracy map $(\Delta_N)_0 \lra (\Delta_N)_n$
is continuous.

Let us now follow [Grayson, Remark 2.4.1-2] and see that the Hausdorffisation of the topological space corresponding 
to the Skorokhod space of maps from $[0,1]_\leq$ to $(\Delta_N)_\diag$ is the geometric realisation of $\Delta_N$:
$$
|\Delta_N|=\left(\homm\sFilth{[0,1]_\leq}{ (\Delta_N)_\diag }_{\ttt\inv}\right)_\text{Hausdorff}$$
Let us explicitly describe the underlying simplicial set. 
By Yoneda Lemma $$\homm\sSets{[0,1]_\leq}{ \Delta_{N}}=\homm\leq{[0,1]_\leq}{{(N+1)}_\leq}$$ 
For $M>0$ and $n>0$, 
a map $\overrightarrow f: {[0,1]_\leq\times \Delta_{M}}\lra { \Delta_{N}}$ of sSets is necessarily of form
$$ (t_1\leq ... \leq t_n, n_\leq \xra\theta M_\leq) \longmapsto (f_{\theta(1)}(t_1), f_{\theta(2)}(t_2),..,f_{\theta(n)}(t_n))$$
where $f_i:[0,1]_\leq\lra (N+1)_\leq$, $1\leq i\leq M$.

Let us now verify the following. The Skorokhod filter on $0$-simplicies is antidiscrete.
A subset $U$ of the set of $1$-simplicies is a Skorokhod neighbourhood iff a $1$-simplex $\overrightarrow f=(f,g)\in U$
whenever 
$\dist(f,g)<\updelta$
where Skorokhod distance  $\dist(f,g)$ is defined in \S\ref{Skorokhod-Grayson}.

Let us check this. In dimension $0$, in the definition of the $\varepsilon\delta$-Skorokhod neighbourhood
necessarily each  $0$-simplex of $[0,1]_\leq$, resp. $\Delta_N$, is $\delta$-small, resp.~$\varepsilon$-small,
hence  each function is  $\varepsilon\delta$-Skorokhod-small.
In dimension $1$, we may assume that $\varepsilon$ is as small as possible, i.e.~the diagonal, and that 
$\delta=\{(t_0\leq t_1\leq t_2\leq t_3): t_3\leq t_0+\updelta\}$ and 
$\delta_0=\{(t_0\leq t_1): t_1\leq t_0+\updelta_0\}$
for some $\updelta>0$ and $\updelta_0>0$. 
``Each $x\in\delta_0$'' means we take arbitrary $t_0\leq t_1\leq t_0+\updelta$.
Choosing an $\updelta$-small ``continuation'' $x'$ of $x$ amounts to choosing 
$(t_2\leq t_3)$ such that $t_0\leq t_1\leq t_2\leq t_3\leq t_0+\updelta$, 
and that its ``tail'' maps into something $\varepsilon$-small means that
$f(t_2)=g(t_3)$.

\section{A convenient category for topology}\label{sec:2}

\subsection{The intuition of general topology}\label{top-sf-intuit}\label{sec:2.1}

Here we argue that the category $\sFilth$ of simplicial filters (see \S\ref{filth:def} and \S\ref{def:filt} for a definition) 
is one of the `` structures which give a 
mathematical content to the intuitive notions of {\em limit, continuity and neigh- 
bourhood}'' and that the intuition of general topology as described by  \href{http://mishap.sdf.org/tmp/Bourbaki_General_Topology.djvu#page=15}{[Bourbaki]} 
applies to $\sFilth$ almost verbatim.
% but somewhat more flexible than the standard ones.  
%The category $\sFilth$ as such a structure 
%is ``the subject matter of the present [paper]. 

We do so by paraphrasing the Introduction of \href{http://mishap.sdf.org/tmp/Bourbaki_General_Topology.djvu#page=15}{[Bourbaki]}, which we quote in full for reader's convenience.
\newline\includegraphics[width=\linewidth]{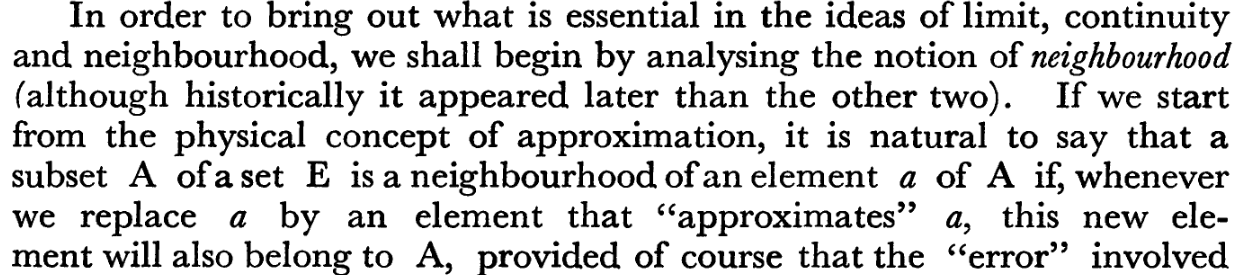}
\newline\includegraphics[width=\linewidth]{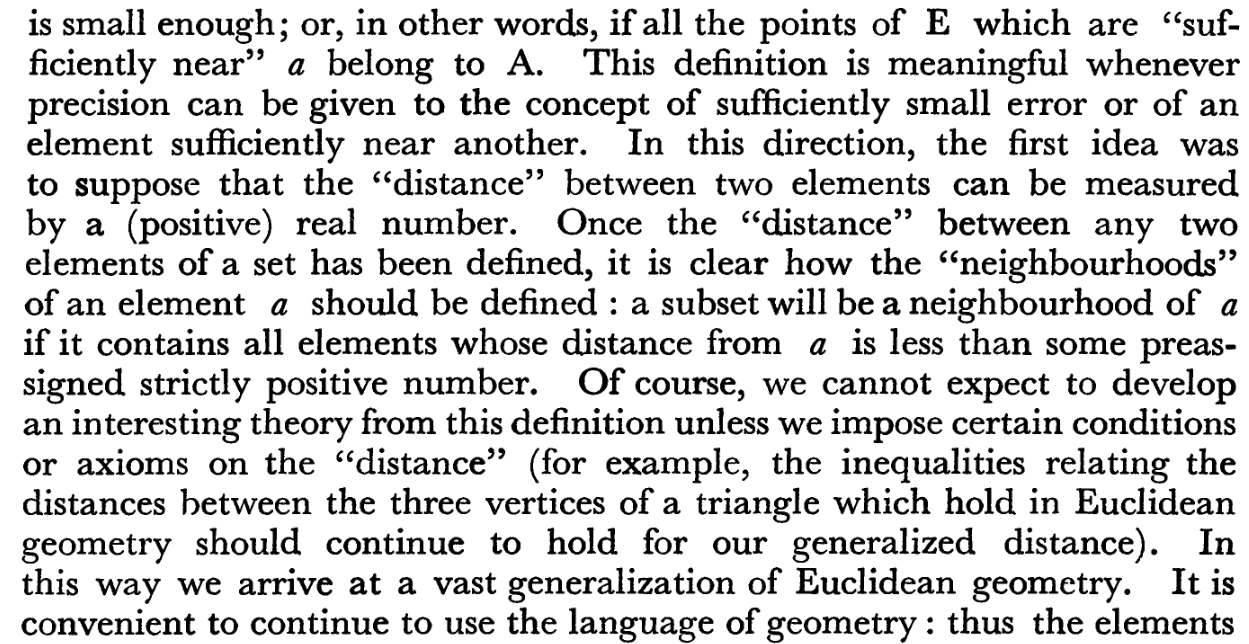}
\newline\includegraphics[width=\linewidth]{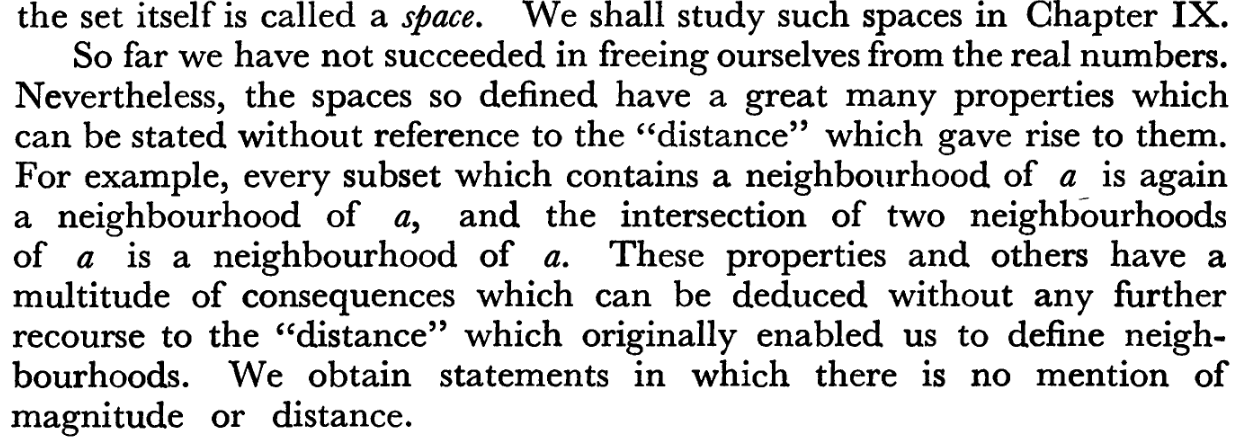}
\newline\includegraphics[width=\linewidth]{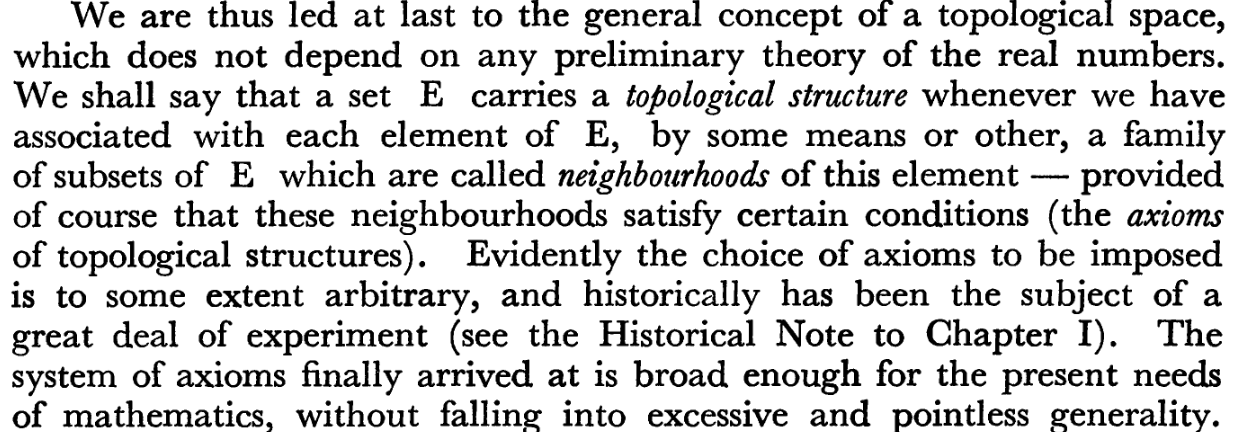}
\newline\includegraphics[width=\linewidth]{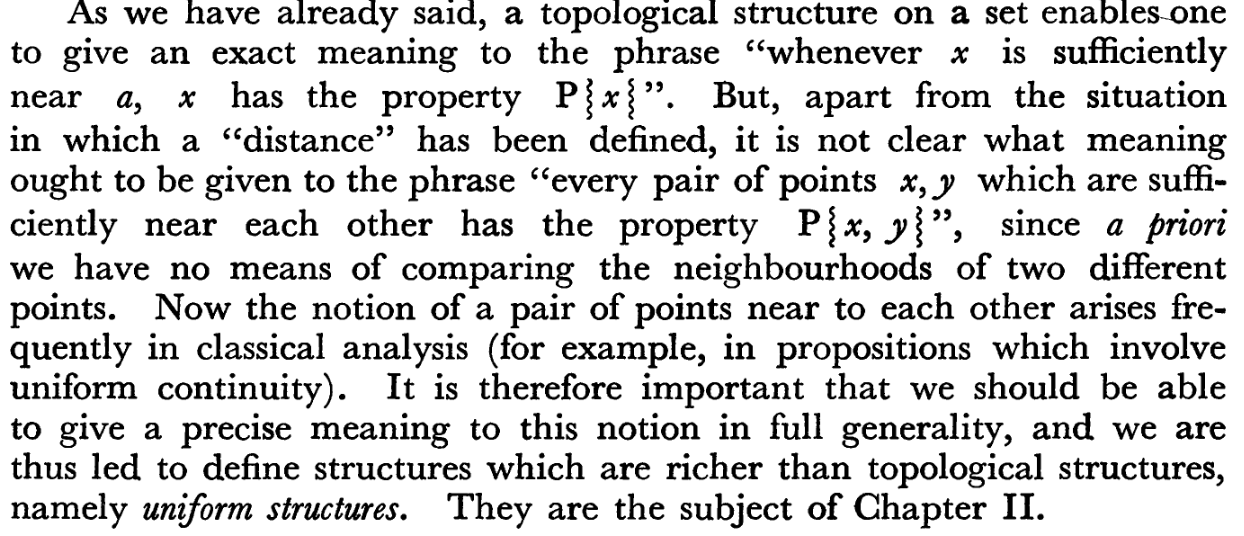}
With help of  $\sFilth$ ``precision can be given to the concept of sufficiently small error''.
In this direction, our idea was to suppose  
an ``error'' is but a pair $(a,a')\in E\times E$ of elements of $E$ 
and to be thought of as {\em a $1$-simplex of a simplicial set}, and that 
to give a precise meaning to the concept of sufficiently small error it is enough
{\em to associate with $E$, by some means or other, 
for each finite Cartesian power of $E$, 
 a family of subsets of
``$n$-simplicies'' $E\times...\times E$ which are called neighbourhoods} ---  provided
of course that these neighbourhoods satisfy certain conditions (namely, form
an object $E_\bullet:\sFilth$, 
i.e.~a contravariant functor from the category of finite linear orders 
to the category of filters.)

Following up ``the first idea'' \href{http://mishap.sdf.org/tmp/Bourbaki_General_Topology.djvu#page=16}{[Bourbaki]}, 
let us first suppose that a set $E$ 
carries a notion of ``distance" between two elements 
which 
can be measured
by a (positive) real number. Once the ``distance" between any two
elements of a set has been defined, it is clear how the ``neighbourhoods" on $E\times E$ in $\sFilth$
%of an element a 
should be defined: a subset of $E\times E$ will be a neighbourhood %of a
if 
for every element $a$ it contains all pairs $(a,a')$ of 
elements whose distance 
$\dist(a,a')$ 
is less than some preas-
signed strictly positive number. 
 Of course, we cannot expect to develop 
an interesting theory from this definition unless we impose certain conditions 
or axioms on the ``distance" (for example, the inequalities relating the 
distances between the three vertices of a triangle which hold in Euclidean 
geometry should continue to hold for our generalized distance). 
%It turns out that the conditions imposed on ``distance" by the definition of a metric
%are enough to obtain an object $E_\bullet$ of $\sFilth$ such that 
%these neighbourhoods form the filter on the set of its $1$-simplicies $E_1=E\times E$; see \S\ref{metric-filth} 
%for details.

 In
this way we arrive at a %vast 
generalization of topology. %Euclidean geometry. 
It is
convenient to continue to use the language of topology: %geometry: 
 thus the $0$-simplicies %elements of a set 
on which a ``distance" has been defined are called {\em points}, and
an object of the category $\sFilth$ %set 
itself is called a {\em space}.

%\hall study such spaces in Chapter IX. 
%%So far we have not succeeded in freeing ourselves from the real numbers. 
Nevertheless, the $\sFilth$-spaces so defined have a great many properties which 
can be stated without reference to the ``distance" which gave rise to them. 
For example, every subset which contains a neighbourhood %of a 
is again a neighbourhood, %of a, 
and the intersection of two neighbourhoods 
%of a 
is a neighbourhood, % of a.
and, more generally, the neighbourhoods form a functor $\Dop \lra \Filth$. 
These properties and others have a 
multitude of consequences which can be deduced without any further 
recourse to the ``distance" which originally enabled us to define neigh- 
bourhoods. We obtain statements in which there is no mention of 
magnitude or distance. 

We are thus led at last to the general concept of a %topological 
$\sFilth$-space, 
which does not depend on any preliminary theory of the real numbers or topology. 
We shall say that a set $E$ carries a %topological 
$\sFilth$-structure whenever we have 
associated with each finite Cartesian power of $E$, by some means or other, a family 
of subsets of $E\times .. \times E$ which are called neighbourhoods % of this element 
--- provided 
of course that these neighbourhoods satisfy certain conditions 
(namely, form a simplicial object  $\Dop\lra \Filth$ in the category of filters; see \S\ref{nei-simplicial} 
for an explanation how to reformulate the axioms of topology in terms of neighbourhoods
as being a simplicial object).
Of course, there are $\sFilth$-spaces which are not associated to a set in this way. 
%the axioms 
%of topological structures)

The goal of this paper is to suggest to the reader that it may be worthwhile to 
view $\sFilth$ as a replacement of the axioms of topology and to consider
the question 
whether the 
system of axioms represented by $\sFilth$ is broad enough for the present needs 
of topology/mathematics, without falling into excessive and pointless generality.

As [Bourbaki] have said, a topological structure on a set enables one 
to give an exact meaning to the phrase 
``whenever $x$ is sufficiently 
near $a$, $x$ has the property $P\{x\}$''.
  But, apart from the situation 
in which a ``distance" has been defined, it is not clear what meaning 
ought to be given to the phrase 
``every pair of points $x$,$y$ which are suffi- 
ciently near each other has the property $P\{ x, y\}$'',
  since a priori 
we have no means of comparing the neighbourhoods of two different 
points. 
%%% 
%%% 
%%% %The $\sFilth$-structure defined above enables one 
%%% to give an exact meaning to the phrase ``whenever $x$ is sufficiently 
%%% near $a$, $x$ has the property $P\{x\}$''. But %, 
%%% %apart from the situation 
%%% %in which a ``distance" has been defined, it is not 
%%% it does not make 
%%% clear what meaning 
%%% ought to be given to the phrase ``every pair of points $x$,$y$ which are suffi- 
%%% ciently near each other has the property $P\{ x, y\}$, since %a priori 
%%% we have 
%%% no means of comparing the neighbourhoods of two different 
%%% points'' using this $\sFilth$-structure. 
%%% 
Now the notion of a pair of points near to each other arises fre- 
quently in classical analysis (for example, in propositions which involve 
uniform continuity). It is therefore important that we should be able 
to give a precise meaning to this notion in full generality, and we are 
thus led to define $\sFilth$-structures which are richer than ones associated with topological structures, %namely 
%$\sFilth$-structures essentially equivalent to
and in fact are  associated with uniform structures which are %. They are 
the subject of Chapter II of \href{http://mishap.sdf.org/tmp/Bourbaki_General_Topology.djvu#page=176}{[Bourbaki]}.

We do this as follows; we speculate that the fact that this is possible is an indication that the notion of an $\sFilth$-space is more flexible
than the usual notion of a topological space.

Whenever a ``distance" has been defined, to give a precise meaning to the notion of a pair of points near to each other,  
we associate with it an $\sFilth$-object such that its filter of $1$-simplicies is defined as follows:  
a subset of $E\times E$ will be a neighbourhood %of a
if
it contains all the pairs $(a,a')$ of
elements whose distance
$\dist(a,a')$ 
is less than some preas-
signed strictly positive number (note that in the previous construction of the $\sFilth$-structure corresponding to a topology, 
this number was allowed to depend on $a$,
and this is why we had no means to compare  the neighbourhoods of two different
points). More generally, a subset of $E\times...\times E$ will be a neighbourhood %of a
if
it contains all the tuples  $(a_1,...,a_n)$ of
elements such that the distance
$\dist(a_i,a_j)$ for all $1\leq i\leq j\leq n$
is less than some preas-
signed strictly positive number.
%\subsection{Continuity} 
\newline\includegraphics[width=\linewidth]{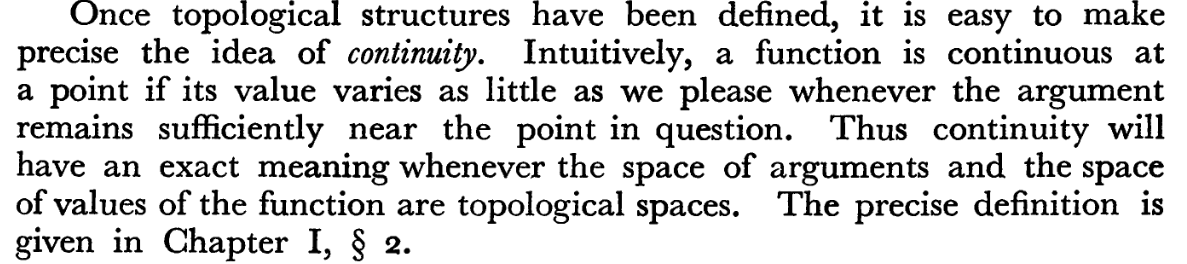}
%\begin{quote}
%Once topological structures have been defined, it is easy to make 
%precise the idea of continuiry. Intuitively, a function is continuous at 
%a point if its value varies as little as we please whenever the argument 
%remains sufficiently near the point in question. 
%\end{quote}
%

Paraphrasing slightly these intuitive words of \href{http://mishap.sdf.org/tmp/Bourbaki_General_Topology.djvu#page=15}{[Bourbaki,p.19]}, we say 
that intuitively, a function is continuous at 
a point if  its value remains the same up to an error as small as we please 
whenever the argument remains the same up to a sufficiently small error.
The precise meaning of this phrase in terms of $\sFilth$ is straightforward:
a map $f:X_\bullet\lra Y_\bullet$ in $\sFilth$ is continuous iff for every $n\geq 0$,
for every neighbourhood $\varepsilon\subset Y_n$ there is a neighbourhood $\delta\subset X_n$ such that 
$f(\delta)\subset\varepsilon$; equivalently, $f^{\inv}(\varepsilon)$ is a neighbourhood in $X_n$.

Note that when we consider the $\sFilth$-structures defined above and take $n=2$,
we recover the standard definitions of continuity and uniform continuity \href{http://mishap.sdf.org/tmp/Bourbaki_General_Topology.djvu#page=179}{[Bourbaki, I\S2.1,II\S2.1]}. 
%
%
%%These intuitive words [Bourbaki,p.19] are formalised as follows, %in $\sFilt$ as follows, and become
%%and become the definition of a morphism  in the category of simplicial sets equipped with a notion of neighbourhood: 
%%%$\sFilt$: 
%%a map $f:X\lra Y$ of sSets is {\em continuous} iff 
%%for every $m>0$, 
%%for every neighbourhood $\delta\subset Y_m$ there is a neighbourhood $\varepsilon\subset X_m$ such that
%%$f(\varepsilon)\subset \delta$. Note that above we only described the notion of a neighbourhood for $m=1$,
%%and in this case the condition on $f:X\times X\lra Y\times Y$ is the definition of continuity, resp.~uniform continuity,
%%when we consider the notion of a neighbourhood associated with a topological, resp.~metric, structure.
%%
The similarity to the exposition of the definition of uniform continuity is particularly startling:
\newline\includegraphics[width=\linewidth]{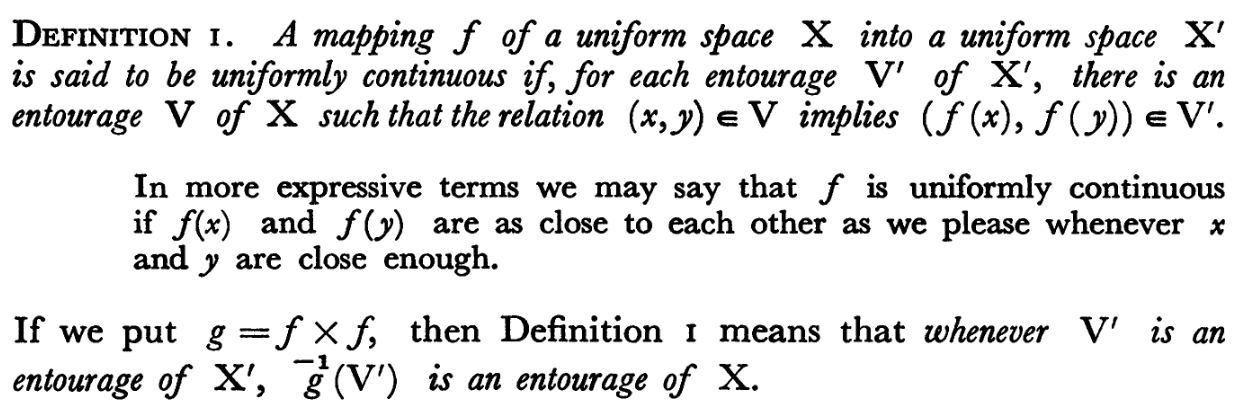}

\subsubsection{The notion of limit via an endofunctor ``shifting'' dimension} 
We now show how to rewrite the definition of a limit of a function 
as a Quillen lifting property involving an endofunctor of $[+1]:\sFilth\lra\sFilth$ ``shifting'' dimension.

% The notion of a filter, which is thus inseparable
Thus here we show that the notion of the ``shift endofunctor $[+1]:\Dop\lra\Dop$
is related to 
%from 
 that of a limit, and later we show that it appears also in other contexts in topology
involving local properties, %; %for example,
%the neighbourhoods of a point in a topological space form a filter.
namely the definition of a locally trivial bundle, cf.~\S\ref{def-bundle-sample} and \S\ref{def-bundle}.
\newline\includegraphics[width=\linewidth]{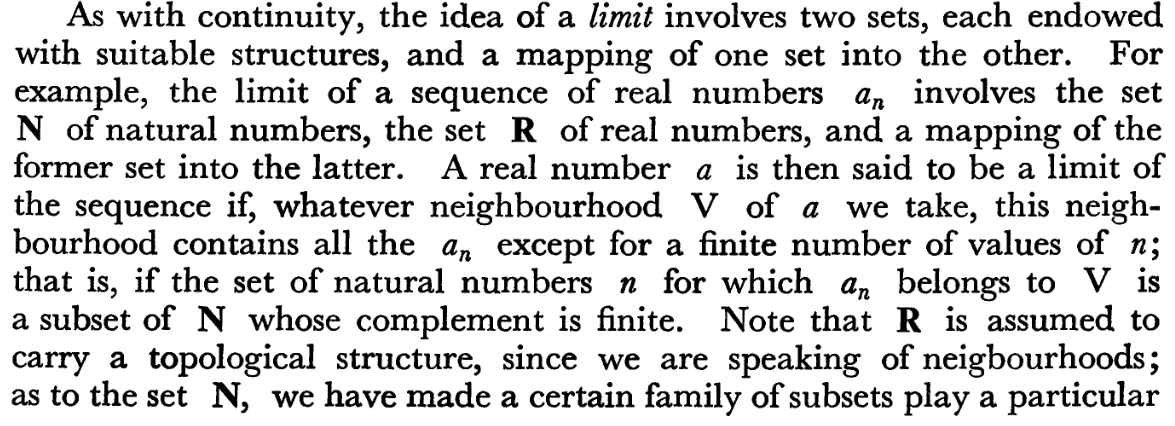}
\newline\includegraphics[width=\linewidth]{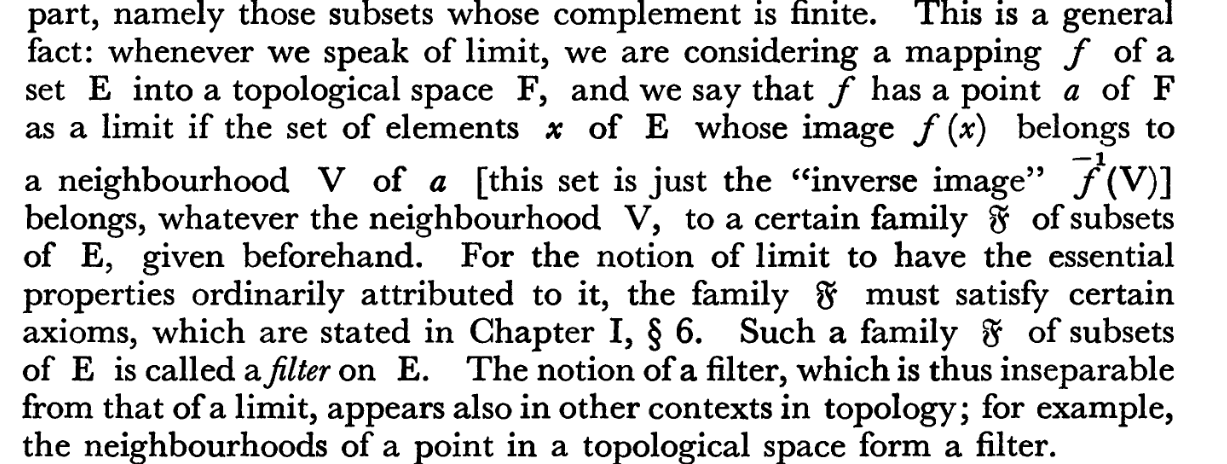}
%This is a general 
%fact: whenever we speak of limit, we are considering a mapping $f$ of a 
%set $E$ into a topological space $F$, and we say that $f$ has a point $a$ of $F$ 
%as a limit if the set of elements $x$ of $E$ whose image $f (x)$ belongs to 
%a neighbourhood $V$ of $a$ [this set is just the "inverse image" $f^{-1}(V)$] 
%belongs, whatever the neighbourhood $V$, to a certain family $\mathfrak F$
% of subsets 
%of $E$, given beforehand. 
%

Whenever we speak of limit, we are considering a mapping $f$ of a 
set $E$ into a topological space $F$, and we say that $f$ has a point $a$ of $F$ 
as a limit if the set of elements $x$ of $E$ whose image $f (x)$ belongs to 
a neighbourhood $V$ of $a$ [this set is just the "inverse image" $f^{-1}(V)$] 
belongs, whatever the neighbourhood $V$, to a certain family $\mathfrak F$
 of subsets 
of $E$, given beforehand. 

In terms of $\sFilth$, this is expressed as follows: the mapping 
$f_a: E\lra F\times F, \ x\mapsto (a,f(x))$ is continuous 
with respect to the filter on $E$ defined by the family  $\mathfrak F$,
and the filter of neighbourhoods on $F\times F$ associated 
with the topology on $F$. However %, in terms of $\sFilth$, 
this is a mapping
from $0$-simplices to $1$-simplices, and thus ``shifts'' dimension: 
this is not a problem, as the category $\Dop$ of finite linear orders admits an endofunctor 
$[+1]:\Dop\lra \Dop$ equipped with a natural transformation $[+1]\implies \id$,
and therefore $\sFilth$ admits an endofunctor $[+1]:\sFilth\lra \sFilth$
shifting dimension equipped  with a natural transformation $[+1]\implies \id$.

Considerations above lead to the following reformulation of the notion of a limit 
in terms of $\sFilth$:

\bi
\item To give a mapping of sets $f:E\lra F$ is to give a map of simplicial sets
$$\overrightarrow f: \homm\Sets-E\lra \homm\Sets-F,\ \ \ (x_1,..,x_n)\longmapsto (f(x_1),...,f(x_n))$$ 
\item {\em The mapping $f$ has a point $a$ of $F$
as a limit iff $\overrightarrow f$ factors as 
$$ \homm\Sets-E\xra{\overrightarrow{f_a}} \homm\Sets-F\circ[+1] \lra \Homm\sSets-F$$
$$  \overrightarrow{f_a}: (x_1,..,x_n)\longmapsto (a, f(x_1),...,f(x_n))$$
via a map continuous with respect to appropriate filters.
}\ei

On $\homm\Sets-F$ these are the filters associated with the topology on $F$.
On $ \homm\Sets-E$ these are {\em the $\mathfrak F$-diagonal filters}, i.e.~the finest filters such that the filter on $0$-simplicies is $\mathfrak F$
and the degeneration map from the set of $0$-simplicies to the set of $n$-simplicies is continuous for each $n>0$;
explicitly, a subset of $E^n=\homm\Sets n E$ is a neighbourhood iff it contains
$\{(x,x,..,x): x\in\varepsilon\}$ for some $\varepsilon\in \mathfrak F$. 

Let us check this is indeed a reformulation. As simplicial sets, the source of the morphism is connected 
whereas the target is not and its connected components are parametrised by $a\in F$. 
Hence, to pick such a decomposition  is to pick an $a\in F$. 
Above we saw that the morphism $(\overrightarrow{f_a})_0: E\lra F\times F$ 
of $0$-simplicies is continuous iff $f$ has a point $a$ of $F$ as a limit.
A slight extension of this argument shows this holds for $\overrightarrow{f_a}$ itself.
See details at \S\ref{limits-samples} and \S\ref{limits-defs}.

%%
%%This $\sFilth$-structure associated with a metric 
%%gives  means of comparing the neighbourhoods of two different 
%%points, and
%%allows to 
%%give a precise meaning  to the phrase " phrase "every pair of points x,y which are suffi- 
%%ciently near each other has the property $P\{x,y\}$".%, since a priori 
%%%we have no means of comparing the neighbourhoods of two different 
%%%points. Now the 
%%The 
%%notion of a pair of points near to each other arises fre- 
%%quently in classical analysis (for example, in propositions which involve 
%%uniform continuity)%. It is therefore important that we should be able 
%%, and in fact $\sFilth$-structure associated with a metric is equivalent to 
%%%to give a precise meaning to this notion in full generality, and we are 
%%%thus led to define structures which are richer than topological structures, 
%%%namely 
%%{\em uniform structures} 
%%which are %They are 
%%the subject of Chapter II of \href{http://mishap.sdf.org/tmp/Bourbaki_General_Topology.djvu}{[Bourbaki]}. 
%%

\subsubsection{Axioms of topology as being simplicial}\label{nei-simplicial}
A topology is a collection of (filters of) neighbourhoods of points compatible in some sense. We now show
that it is ``compatible'' in the sense that it is ``functorial'', i.e.~defines a functor from $\Dop$
to the category of filters.
%We now show  that the axioms of topology in terms of neighbourhoods as stated in 

This is almost explicit in the axioms $(\text{V}_\text{I})$-$(\text{V}_\text{IV})$ of \href{http://mishap.sdf.org/tmp/Bourbaki_General_Topology.djvu#page=24}{[Bourbaki,I\S1.2]} 
of topology in terms of neighbourhoods if we rewrite them in terms of $\sFilth$. % the axioms $(\text{V}_\text{I})$-$(\text{V}_\text{IV})$ of [Bourbaki,I\S1.2].
We now quote: 
%say that the neighbourhoods give rise to a 2-dimensional  simplicial object in the category of filters.
% Let us quote  [Bourbaki,I\S1.2]:
\newline\includegraphics[width=\linewidth]{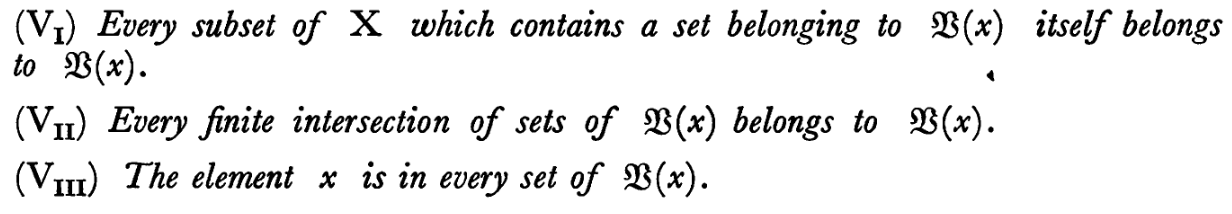}
\newline\includegraphics[width=\linewidth]{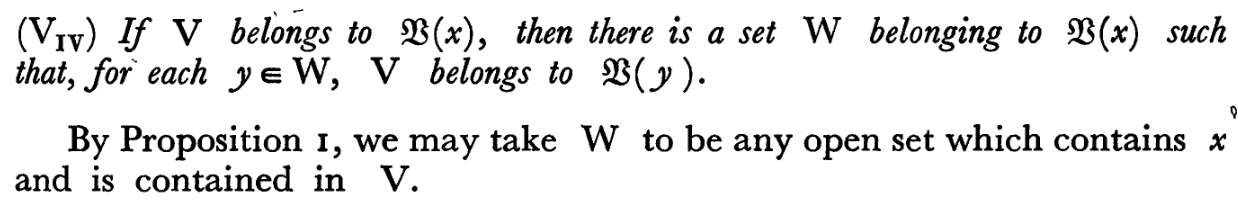}
\newline\includegraphics[width=\linewidth]{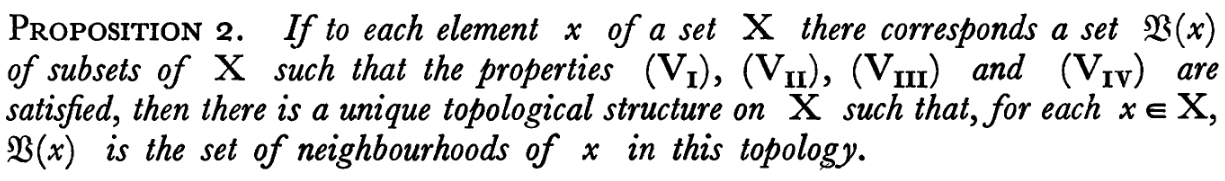}
\newline\includegraphics[width=\linewidth]{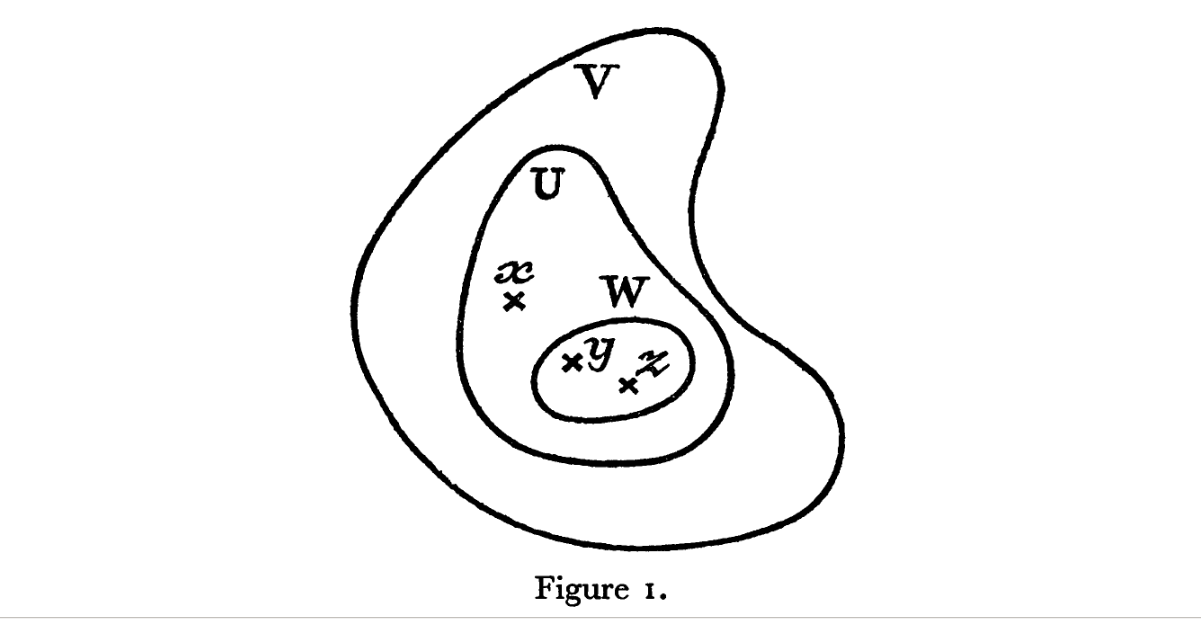}

Call a subset of $X\times X$ a neighbourhood iff it is of the form 
$$\bigcup\limits_{x\in X} \{x\}\times U_x \text{ where } U_x\in \mathfrak B(x)$$
Axiom $(\text{V}_\text{I})$ says that a subset containing a neighbourhood is itself a neighbourhood.
Axiom  $(\text{V}_\text{II})$ says that the neighbourhoods are closed under finite intersection. 
Hence, the first two axioms say that the neighbourhoods so defined form a filter on $E\times E$. 

Axiom  $(\text{V}_\text{III})$ states the continuity of the diagonal map $E\lra E\times E, \ x\longmapsto (x,x)$ 
from the set $E$ equipped with the antidiscrete filter (i.e.~the filter where $E$ itself is the unique neighbourhood). 

Axiom $(\text{V}_\text{IV})$ needs a little argument, as follows.

Equip $E\times E\times E$ with the coarsest filter such that %the two projections on the first 
the following coordinate projections $E\times E\times E\lra E\times E$ are continuous:
$$(x,y,z)\longmapsto (x,y) \text{ and } (x,y,z)\longmapsto (y,z)$$
Axiom $(\text{V}_\text{IV})$ says that the remaining coordinate projection $(x,y,z)\longmapsto (x,z)$ is continuous. 
To see this, consider the preimage of a neighbourhood containing $\{x\}\times V$. By continuity,
there are neighbourhoods $\bigcup\limits_{x'\in X} \{x'\}\times W_{x'}$ and $\bigcup\limits_{y\in X} \{y\}\times V_y$ such that
$$\{x\}\times W_x \times X \,\bigcap\, X\times (\bigcup\limits_{y\in X}   \{y\}\times V_y) \subset   \{x\}\times X\times V$$
Now take $W=W_x$ and see that $V_y\subset V$ for each $y\in W_x$.

Finally, note that the considerations above amount to the following
\newline\includegraphics[width=1.01\linewidth]{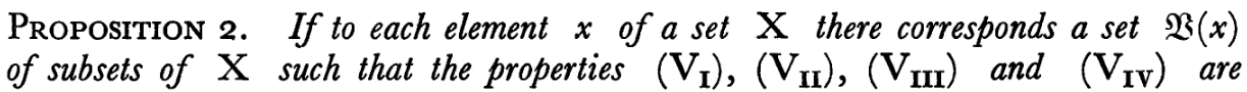}
{\small\em satisfied, then there is a ``unique 2-dimensional $\sFilth$-structure on $X$'' such that
the set of neighbourhoods in $X\times X$ is 
$$\left\{ \bigcup\limits_{x\in X} \{x\}\times U_x \,\,:\,\,  U_x\in \mathfrak {B}(x) \right\} 
.$$}

By this we mean the following:  there is a unique object of $\sFilth$ such that 

%Axioms~$(\text{V}_\text{I})$-$(\text{V}_\text{IV})$ hold iff there is an $\sFilth$-object $X_\bullet$ such that
\bi
\item (``$\sFilth$-structure on $X$'') its underlying simplicial set is $\homm\Sets - X$
\item $(\text{V}_\text{III})$ the set $X_0=X$ of $0$-simplicies carries antidiscrete topology
\item $((\text{V}_\text{I})\&(\text{V}_\text{II}))$ the filter on the set $X_1=X\times X$ of $1$-simplicies is
$$\left\{ \bigcup\limits_{x\in X} \{x\}\times U_x \,\,:\,\,  U_x\in \mathfrak {B}(x) \right\} 
$$
\item $((\text{V}_\text{IV}))$  $X_\bullet$ is 2-dimensional, i.e.~the filter on $X_{n+1}=X^n$ is the coarsest filter such that the face maps 
$$X_{n+1}=X^n\lra X_{1}=X\times X,\ \  (x_1,...,x_n)\longmapsto (x_i,x_{i+1}),\ 0<i<n$$
are continuous.
\ei

A similar reformulation can be given to the axioms of uniform structure \href{http://mishap.sdf.org/tmp/Bourbaki_General_Topology.djvu#page=175}{[Bourbaki,II\S1.1,Def.1]},
cf.~Exercise~\ref{uniform_simplicial}.

\subsubsection{A flexible notion of space}  
%In 
%%...%this way we arrive at a %vast 
%%...%generalization of topology. %Euclidean geometry. 
%%...%It is 
%%...%convenient to continue to use the language of topology:%geometry: 
%%...%thus the $0$-simplicies %elements of a set 
%%...%on which a ``distance" has been defined are called {\em points}, and 
%%...%the object of $\sSets$ %set 
%%...%itself is called a {\em space}.
%%...%
The discussion above suggests this notion of space is somewhat more flexible than the usual notion of 
a topological space. Filters, uniform structures and topological spaces are $\sFilth$-spaces. 
and a limit of a function is an $\sFilth$-morphism.  

This allows to talk in category-theoretic terms about equicontinuous sequences of functions
and their limits, by considering 
$$ \homm\Sets-\NN\times\homm\Sets-X\lra \homm\Sets-Y$$
$$  \homm\Sets-\NN\times\homm\Sets-X\circ [+1]\lra \homm\Sets-Y\circ[+1]$$
where the simplicial sets are equipped with various filters
representing the topology or metric on $X$ and $Y$, and 
the filter of cofinite subsets of $\NN$. See~\S\ref{limits-samples} and \S\ref{ascoli} for a discussion.  

We saw how  the endofunctor $[+1]:\sFilth\lra \sFilth$ is used to talk about limits, a local notion.
In a similar way it can be used to talk about families of functions defined locally. 
Let $X_\bullet$ and $Y_\bullet$ be  $\sFilth$-objects corresponding to  topological spaces $X$ and $Y$.
To give an $\sFilth$-morphism  $X_\bullet\circ[+1]\lra Y_\bullet$ is to give a family of functions $f_x:X\lra Y$, $x\in X$
such that $f_x$ is continuous in a neighbourhood of $x$ (under some assumptions on $X$ and $Y$). 
The definition of local triviality is formulated in terms of $\sFilth$ as follows:
a $X\xra p B$ is a locally trivial bundle with fibre $F$ iff in $\sFilth$ there is 
an isomorphism $B_\bullet\circ[+1]\times_{B_\bullet} X_\bullet \xra{(iso)} B_\bullet\circ[+1]\times F_\bullet$ over $B_\bullet$,
cf.~\S\ref{def-bundle-sample} and \S\ref{def-bundle} for explanation.

\subsection{Algebraic topology}\label{sec:2.2}
Here we offer a couple of vague speculations about a possible intuition in $\sFilth$ originating in homotopy theory.

\subsubsection{Simplicies as $\varepsilon$-discretised homotopies} 
In a metric space, % two points $x$ and $y$ joined by a path can be joined by a sequence of points $x_1,..,x_n$
for $\epsilon>0$ small enough, we may think of a sequence $x_0,..,x_T$, $\dist(x_t,x_{t+1})<\epsilon, t=0,..,T-1$ as an 
{\em $\epsilon$-discretised homotopy} from $x_0$ to $x_T$, and the indices $t,T$ as {\em time}. 
In terms of $\sFilth$, 
this sequence $\overrightarrow x = (x_0,...,x_T)$ is an $T$-simplex %of $M_\bullet$ of 
the $\sFilth$-object $M_\bullet$ corresponding to the metric space,
and the condition $\dist(x_t,x_{t+1})<\epsilon, t=0,..,T-1$ says that the ``consecutive'' faces 
$\overrightarrow x[t<t+1]\in \varepsilon$ where $\varepsilon:=\{(x,y):\dist(x,y)<\epsilon\}$ 
is the neighbourhood associated with distance $\epsilon$.  
This consideration suggests a generalisation of this intuition to an arbitrary object $X:\sFilth$ :
for a small enough neighbourhood $\varepsilon\subset X_\tau$ for $\tau<t'<<T$, think of 
a simplex $s\in X_T$ as {\em $\varepsilon$-discretised homotopy} provided its 
``consecutive'' faces $s[t_0\!<\!...\!<\!t_\tau]\in\varepsilon$ whenever $t_\tau\leq t_0+t'$.
 Note that the essential asymmetry (direction of time) of this notion, or rather intuition, of homotopy, which is apparently
%in a comment by M.Kontsevich probably related to sFilt
a desirable property in the context of $\infty$-categories.

A well-known lemma says that two functions $f,g:A\lra M$ from an arbitrary topological space $A$ to a (sufficiently nice) 
metric space $M$
are homotopic iff for each $\epsilon>0$ there is a $\epsilon$-discretised homotopy 
$f=f_0,..., f_T=g$ such that for any $x\in A$ $\sup_x\dist(f_t(x), f_{t+1}(x))<\epsilon$.
This allows us to think of homotopies of functions as simplicies in a certain function space with $\sup$-metric. 

\subsubsection{Inner Hom and mapping spaces} Think of
the inner hom of the underlying simplicial sets of objects of $\sFilth$ as a space of {\em discontinuous} functions. 
Such spaces of  discontinuous functions are considered in probability theory as 
spaces of functions describing stochastic process, and there are standard metrics
called {\em Levi, Levi-Prokhorov, and Skorokhod metric} on such spaces [Kolmogorov]. We need
the following non-symmetric variant of Skorokhod metric defined on functions between metric spaces:
$$\dist(f,g):=\inf\left\{\epsilon>0: \forall x \exists y ( \dist(x,y)<\epsilon \,\,\&\,\, \dist( f(x),g(y))<\epsilon\right\}
$$ 
This definition generalises to $\sFilth$ as follows (see \S\ref{skorokhod-filter-def} details):
For $N>2n$, $\delta\subset X_N$ and $\varepsilon\subset Y_n$,
a {\em $\varepsilon\delta$-Skorokhod neighbourhood in the Hom-set $\hom X Y$} is the subset %of $\Hom X Y$ 
consisting
of all the function $f:X\lra Y$ with the following property:
\bi\item[] there is a neighbourhood $\delta_0\subset X_n$ such that
each ``$\delta_0$-small'' $x\in \delta_0$ has a ``$\delta$-small''
``continuation'' $x'\in X_N$, $x=x'[1..N]$ such that its ``tail'' maps into something ``$\varepsilon$-small'',
i.e.~$f(x'[N-n+1..N])\in\varepsilon$.
\ei
This defines the {\em Skorokhod filter on the Hom-set $\hom X Y$} 
and thereby {\em Skorokhod neighbourhood structure on the inner Hom $\Homm\sSets X Y $}
of the underlying simplicial sets of $X$ and $Y$,
and in fact a functor 
$$\Homm\sFilth--: \sFilth^\text{op}\times \sFilth \lra \sFilth$$
which we call {\em Skorokhod mapping space}. 

The Skorokhod spaces have the desirable property (cf.~\S\ref{Skorokhod-mapping-spaces-sample} and \S\ref{Skorokhod-mapping-spaces-defs}) 
that there is an evaluation map
$$\Homm\sFilth A {\Homm\sFilth XY} \lra \Homm\sFilth { A\times X} Y $$
and in fact there is a sort of left adjoint functor
$$
\Homm\sFilth { A\ltimes X} Y \xra{(iso)} \Homm\sFilth A {\Homm\sFilth XY}$$
which admits a map $A\ltimes X \lra A\times X$.

In \S\ref{Skorokhod-Grayson} and \S\ref{Skorokhod-paths-space} we show a construction of 
[Besser], [Grayson, Remark 2.4.1-2], and [Drinfeld] can be interpreted in $\sFilth$ as saying
that  the geometric realisation
may be thought as the Skorokhod space of {\em discontinuous} paths, and thus an endofunctor of $\sFilth$.

\subsection{Ramsey theory and model theory}\label{sec:2.3}

Ramsey theory suggests $\sFilt$-spaces which do not come from topology.
Given an sset $X$ and a colouring $c$ of simplices, call a simplex
$c$-homogeneous iff all of its non-degenerate faces have the same $c$-colour. 
A collection of colourings allows one to define a notion of neighbourhood:
a neighbourhood (of the diagonal) is a subset containing all the simplices homogeneous with respect to 
finitely many colours. This construction allows to generalise the notion of the
Stone space of types in model theory if we consider formulae as 
colourings on the sset $n_\leq\longmapsto \homm\Sets {n} M$ of tuples of elements of a model $M$ : 
the neighbourhoods in the $\sFilt$-Stone space consists of  sequences (with repetitions) 
indiscernible with respect to finitely many formulae. See \S\ref{ramsey}-\ref{model} for details.

\section{Definitions and constructions}\label{sec:3}

\subsection{The main category: the definition
}\label{def:filt}

We now state formally the definition of the key categories of the paper. 

\subsubsection{The two categories of filters (sets with a notion of
smallness.)
}
First we quote the definition of a filter by \href{http://mishap.sdf.org/tmp/Bourbaki_General_Topology.djvu#page=63}{[Bourbaki, I\S6.1, Def.I]}, again. 
\newline\includegraphics[width=\linewidth]{Bourbaki_filter_def.png}
%%A {\em filter} on a set $X$ is a collection of subsets called {\em neighbourhoods}
% topology 
%such that a subset containing a non-empty open subset is necessary open; 
% non-empty open subsets are called {\em big}. 
%%closed under finite intersection and 
%%such that a subset containing a neighbourhood is necessary a neighbourhood as well.
\newline\includegraphics[width=\linewidth]{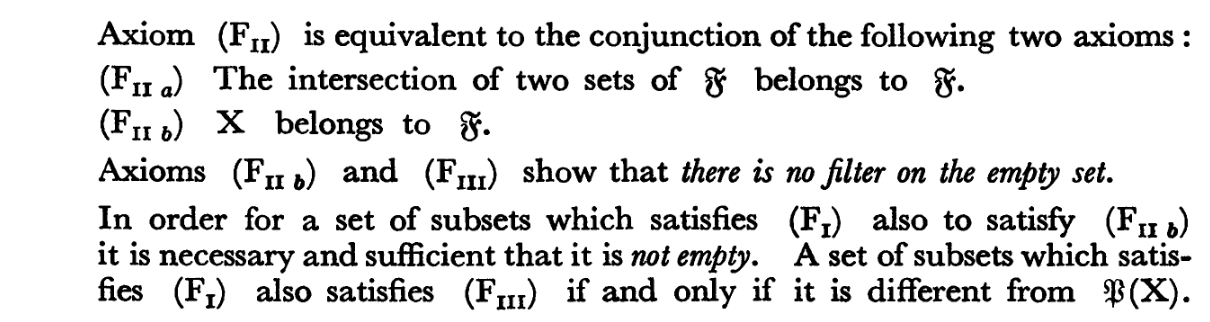}
Subsets in $\mathfrak F$ are called {\em neighbourhoods} or {\em $\mathfrak F$-big}.
Unlike \href{http://mishap.sdf.org/tmp/Bourbaki_General_Topology.djvu}{[Bourbaki]}, we do {\em not} require $(\text{F}_\text{III})$ and allow both $X=\emptyset$ and  $\emptyset\in\frak F$;
necessarily $X\in \mathfrak F$..

\begin{defi}A {\em filter} on a set $X$ is a collection of subsets of $X$ called {\em neighbourhoods}
satisfying  $(\text{F}_\text{I})$  and  $(\text{F}_\text{II})$ above. 
% topology 
%such that a subset containing a non-empty open subset is necessary open; 
% non-empty open subsets are called {\em big}. 
%closed under finite intersection and 
%such that a subset containing a neighbourhood is necessary a neighbourhood as well. 
%Unlike \href{http://mishap.sdf.org/tmp/Bourbaki_General_Topology.djvu}{[Bourbaki]}, we allow $\emptyset$ to be a neighbourhood and we also allow $X$ to be empty.

A {\em morphism of filters} is a function of the underlying sets 
such that the preimage of a neighbourhood is necessarily a neighbourhood;
we call such maps of filters {\em continuous}.  
%A continuous map is a map of filters iff its image intersects each neighbourhood, hence 
%not each map sending everything to a single point is a morphism of filters. 
%Note that requiring non-emptiness condition is stronger than being continuous: a map sending everying everything to a single point 
%is necessarily continuous and is  a morphism of filters iff its image lies in every neighbourhood.

Let $\Filt$ denote the category of filters. Let $\FFilt$ denote the category of filters
where we identify the maps equal almost everywhere: 
$\Ob \FFilt=\Ob \Filt$ and $\homm\FFilt X Y =\homm\Filt X Y /\approx_\FFilt$
where $f\approx_\FFilt g$ iff $\{x:f(x)=g(x)$ is a neighbourhood in $X$.
\end{defi}

We sometimes refer to $\Filt$ as {\em the category of neighbourhoods}.

The {\em category $\Filt$ of filters} is  
equivalent\footnote{The equivalence is given by adjoining/removing the base-point:
a filter $F$ on a set $X$ corresponds to the topological space with points $X\sqcup \{x_0\}$
where a subset is open iff it is of form $\{x_0\}\sqcup U$, $U\in F$ is an $F$-neighbourhood; 
a topological space $X$ with a base-point $x_0$ corresponds to the filter on $X\setminus \{x_0\}$
induced by the neighbourhood filter of $x_0$. The requirement 
``no point other than the base-point
goes into the base-point'' means the map of sets $X\setminus \{x_0\} \lra Y\setminus \{y_0\}$
is well-defined, and continuity at the base-point means precisely that it is a morphism of filters.}
 to the category of topological spaces with a base-point
and morphisms being functions continuous {\em at} the base-point 
(not in a neighbourhood of the base-point), 
and---not a natural requirement---no point other than the base-point
goes into the base-point.

%\begin{enonce}{Examples}[Filters]
\begin{exem}[Filters] \bi\item In a topological space, (not necessary open) neighbourhoods of a point form a filter.
\item The filter consisting of subsets containing a given subset. 
For a set $X$, the filter of diagonals consists of subsets of $X\times ..\times X$ containing $\{(x,...,x): x\in X\}$.
\item Given a measure $\mu$ on a set $X$, subsets of full measure 
%(i.e. those whose complements are (contained in some subset) of measure $0$ 
form  a filter; so do subsets of positive measure. %({\em almost everywhere)}
\item An {\em ultrafilter} on a set $X$ is a filter such that either $A$ or $X \setminus A$ is a neighbourhood, i.e. is either open or closed.
In other words, a not necessarily countably additive ``measure'' such that each subset is measurable and has measure either $0$ or $1$.
\item The filter of cofinite subsets of a set: a subset is a neighbourhood iff its complement is finite.
\item Sets of points of a metric space whose complement has finite diameter form a filter.

\item Sets of pairs of points of a metric space which contain all pairs sufficiently close to each other, i.e.~subsets of $M\times M$ which contain 
the following set for some $\epsilon>0$:
$$\{(x,y)\in M\times M \,:\, \dist(x,y)<\epsilon \} $$

\item Sets of pairs of points of a metric space which contain all pairs sufficiently far apart, i.e.~subsets of $M\times M$ which contain 
the following set for some $d>0$:
$$\{(x,y)\in M\times M \,:\, \dist(x,y)>d \} $$

\item fixme: other examples..
\ei
\end{exem}%\end{enonce}

\subsubsection{Simplicial sets: notations and first definitions}\label{notations}
Recall a preorder $\leq\,\subseteq P\times P$, i.e.~a reflexive transitive binary relation on a set $P$, 
defines a topology on $P$: a subset $U\subseteq P$ is open iff $x\in U$ whenever $x\leq y$ for some $y\in U$. 
Recall also that the preorder also defines a category whose objects are elements of $P$ 
and where all diagrams commute: there is a necessarily unique morphism $x\lra y$ iff $x\leq y$.
Monotone non-decreasing maps correspond to continuous maps and, resp., functors. 

In fact, a preorder on a finite set can be viewed in three equivalent ways: 
as a reflexive transitive relation,  a finite topological space, and a category where all diagrams commute. 
Maps (morphisms) of preorders are then viewed as monotone maps, continuous maps, and functors.

Let $\Delta$ be the category whose objects are finite linear orders $n_\leq:=\{0,..,n-1\}$, $n\in\Bbb N$, and whose morphisms are non-decreasing maps;
note that in $\Delta$ each isomorphism is an identity, in other words, there is a unique object in each isomorphism class.
Let $\Dop$ be its opposite category.
Let $\Delta_\text{big}$ denote the equivalent category of all finite linear orders and non-decreasing maps.
 A {\em simplicial object of a category} $C$ is a  functor $X_\bullet:\Dop\lra C$. 
A {\em simplex} is an element of $X(n_\leq)$ for some $n>0$. An $n$-simplex is an element of $X((n+1)_\leq)$.

An increasing sequence $1\leq t_1\leq ... \leq t_n \leq m$ determines a morphism $m_\leq\lra n_\leq$ in $\Dop$. 
We denote the corresponding {\em faces and degenerations} of a simplex $s$ by $s[t_1\!\leq\! ... \!\leq\! t_n]$ or 
$s[t_1,...,t_n]$. 

In a category $C$, the set of all maps from an object $X$ to $Y$ is denoted by $\homm C X Y$.
The space of maps from $X$ to $Y$, if defined, is denoted by $\Homm C X Y$; typically it is an object of 
the category $C$ itself.

Sometimes we borrow notation from type theory and write $X:C$ to indicate that $X$ is an object of $X$.

%%outdated-rewritten%Finally,  $\sFilt=Func(\Dop, \Filt)$ denote the category of functors (natural transformations)
%%outdated-rewritten%from the opposite of the category $\Delta$ of finite linear orders to the category $\Filt$ of filters. 
%%outdated-rewritten%We refer to $\sFilt$ as  the {\em category of simplicial filters} 
%%outdated-rewritten%or of {\em simplicial neighbourhoods}, for lack of better name. 
%%outdated-rewritten%
%%outdated-rewritten%
%%outdated-rewritten%%%With a set $M$ associate a simplicial set ``co-representable'' functor $\Dop\lra Sets$, $n\mapsto
%%outdated-rewritten%%%Hom_{Sets}(n_{\text{as Set}},M)$, otherwise known as a simplicial set. 
%%outdated-rewritten%
%%outdated-rewritten%A fully faithful embedding of $\Sets$ into $\sSets$ is given by associating to each set 
%%outdated-rewritten%a  ``co-representable'' functor 
%%outdated-rewritten%$$h_X:\Dop\lra \Sets, \ \ \ n\mapsto Hom_{Sets}(n_{\text{as Set}},X)$$
%%outdated-rewritten%
%%outdated-rewritten%In other words, with a set one associates {\em the simplicial object of its Cartesian powers $X^n,n>0$ 
%%outdated-rewritten%and coordinate projections $X^n\lra X^m, (x_1,...,x_n) \mapsto (x_{i_1},...,x_{i_m})$ where $1\leq i_1\leq ...\leq i_m\leq n$}. 
%%outdated-rewritten%

\subsubsection{The definition of the categories of simplicial sets with a notion of
smallness.}

\begin{defi}
Let  $\sFilt=Func(\Dop, \Filt)$ be the category of functors
from 
$\Dop$, the category opposite to the category $\Delta$ of finite linear orders,
to the category $\Filt$ of filters.

Let $\sFFilt$ be the category $\sFilt$ localised as follows:
$\Ob \sFFilt=\Ob\sFilt$, $\homm \sFFilt X Y =\homm\sFilt X Y/\approx$ where 
$f\approx_\sFFilt g$ iff there is $N>0$ such that for every $n>N$ there is a neighbourhood $\varepsilon\subset X_n$ such that
for all $x\in \varepsilon$ it holds $f(x)=g(x)$.
\end{defi}

We think of an object of $\sFilt$ as a simplicial set equipped with a notion of smallness, and that it provides us with 
a notion of {\em a space} which is more flexible than
the notion of a topological space. We suggest no good name for these spaces and refer to such a space
as either a simplicial neighbourhood, 
a neighbourhood structure, 
a simplicial filter; the reader preferring a short word may want to call it a situs. 

We think of an object of $\sFFilt$ as a space where we only care about {\em local} properties.

\subsection{Topological and metric spaces as simplicial filters}
Here we show that $\sFilth$ contains both categories of topological and of metric spaces (with uniformly continuous maps) 
as full subcategories. 

\subsubsection{Topological and metric spaces as neighbourhood structures%simplicial filters
}\label{metric-filth} \label{top-metr-defs} See \S\ref{top-sf-samples} and \S\ref{top-sf-intuit} for examples and intuition. 

%\begin{enonce}{Example}[Simplicial filters]
\begin{defi}[Topological and metric spaces in $\sFilth$]
\label{defi:topoic_and_e-neighbs}
A metric $\dist:X\times X\lra \RR$ on a set $X$ defines {\em a filter of $\varepsilon$-neighbourhoods of the diagonal} 
on $X^n=  \homm{\Sets}{n_{\text{as Set}}} X$: 
a subset $\varepsilon \subset X^n$ is a neighbourhood iff there is an $\epsilon>0$ such that 
$(x_1,..,x_n)\in \varepsilon$ whenever $\dist(x_i,x_j)<\epsilon$ for all $1\leq i,j\leq n$. 

A topology on a set $X$ defines {\em a filter of topoic neighbourhoods of the diagonal} 
on $X^n=  \homm{\Sets}{n_{\text{as Set}}} X$ as follows:
\bee\item for $n=1$, $X$ itself is the unique neighbourhood of $X=X^1$.
\item for $n=2$, $U\subseteq X\times X$ is a neighbourhood iff
for all $x\in X$ there is an open neighbourhood $U_x\ni x$ such that $\{x\}\times U_x \subseteq U$.

\item the filter on $X^n$ is the coarsest filter such that maps $X^n\lra X\times X, (x_1,...,x_n) \mapsto  (x_i,x_{i+1})$ are 
maps of filters for each $0< i<n$. Explicitly, a subset $\varepsilon$ of $X^n$ is {\em a neighbourhood} 
iff  either of the following equivalent conditions holds:
\bi 
\item 
there exist neighbourhoods $\varepsilon_i\subset X\times X$,$0<i<n$, such that $(x_1,...,x_n)\in \varepsilon$ 
whenever for each $0< i<n$, $(x_i,x_{i+1})\in \varepsilon_i$.
\item for each $x\in X$ there a neighbourhood $\varepsilon_x\subseteq X^{n-1}$ such that $\{x\}\times \varepsilon_x\subseteq \varepsilon$
\ei
\eee 
\end{defi}

Denote the simplicial filters associated with a metric and a topology on $X$ by $X_\mU$ and $X_\ttt$, resp. 
\begin{exer} A verification shows that $-_\ttt:\Topp\lra \sFilt$ and $-_\mU:\text{(Metric Spaces, Uniformly Continuous Maps)}\lra \sFilt$
define fully faithful embeddings of the category of topological spaces and the category of metric spaces and uniformly continuous maps.
\end{exer}

\begin{exer} 
%Define the functors $-_\ttt:\sFilt\lra\Topp$ and $-_\mU:\sFilt \lra UniformStructures$ as follows; we slighly abuse notation
%by denoting them by the same letters as the corresponding embeddings.
Check these two embeddings have inverses.

\bi\item
Check that the embedding of topological spaces admits an inverse $-_{\ttt\inv}:\sFilth\lra\Topp$ 
defined as follows. See~\S\ref{sFtoTop} for another brief exposition, and note the similarity 
to the definition of the topology associated with a uniform structure  \href{http://mishap.sdf.org/tmp/Bourbaki_General_Topology.djvu#page=180}{[Bourbaki,II\S1.2,Prop.1,Def.3]}.

The set of points of $X_{\ttt\inv}$ is the set of $0$-simplicies which are $\varepsilon$-small for each neighbourhood $\varepsilon\subset X_0$, i.e.~$X_\text{points}:=\bigcup\limits_{\varepsilon\subset X_0\text{ is a neighbourhood}} \varepsilon$.
The topology is generated by the subsets that together with each point contain all $\varepsilon$-near points
for some $\varepsilon\subset X_1$, i.e.~subsets $U$ with the following property:
there is a neighbourhood $\varepsilon\subset X_1$ such that %neighbourhood base consisting of the following subsets: 
for every $1$-simplex $x\in \varepsilon$ such that both $x[0],x[1]\in X_\text{points}$, it holds that $x[0]\in U$ implies $x[1]\in U$.  

\item Check that for a topological space $X$, $(X_\ttt)_{\ttt\inv}=X$, and that 
$([0.1]_\leq)_{\ttt\inv}=[0,1]$ is the usual unit interval. 

%%incorrect_construction%
%%incorrect_construction%The set of points of the topological space $X$ corresponding to a simplicial neighbourhood $X_\bullet$ 
%%incorrect_construction%is the set $X_0$ of $0$-simplicies equipped with 
%%incorrect_construction%the coarsest topology such that the obvious map $X_\bullet\lra X_\ttt$ is continuous.
%%incorrect_construction%\item Explicitly the topology is described as follows:
%%incorrect_construction%a subset $U_x\ni x$ is a neighbourhood 
%%incorrect_construction%iff there is a neighbourhood $\varepsilon \subset X_1$ such that $s\in U$ whenever $x=s[0]$ and $s[1]\in U_x$.   
%%incorrect_construction%
\item Check whether  the embedding of uniform spaces admits an inverse $-_{\mU\inv}:\sFilth\lra\text{UniformSpaces}$
 defined as follows.

The set of points of the uniform  space $X_{\mU\inv}$ corresponding to a simplicial neighbourhood $X$
is the set $X_0$ of $0$-simplicies. 
The uniform structure on $X_0$ is the coarsest such that the obvious map $X \lra X_{\mU\inv}$ is continuous.
%incorrect% \item Describe the adjoint to the embedding of uniform structures in terms of the localisation $\Dop\lra FiniteSets^\text{op}$.
\ei
\end{exer}

\begin{rema} %(fixme: improve exposition) 
The neighbourhood structure  of topoic subsets associated with a topology 
lacks symmetry of the filter associated with the neighbourhood structure %filter 
of $\varepsilon$-neighbourhoods of the diagonal associated with a metric. 
This accords well with a remark of in the introduction of (Bourbaki, General Topology): 
%\endnote{We quote in full:
%  \href{http://mishap.sdf.org/mints-lifting-property-as-negation/tmp/Bourbaki_General_Topology.djvu}{[Bourbaki, Introduction, p.13]}\newline\noindent 
%\includegraphics[width=\linewidth]{Bourbaki-top.png}
%\includegraphics[width=\linewidth]{Bourbaki-top-uniform-intuition.png}
% \begin{quote}
% `a topological structure now enables us to give precise meaning to the
%  phrase ``such and such a property holds for all points sufficiently near $a$'':
%  by definition this means that the set of points which have this property is
%  a neighbourhood of $a$ for the topological structure in question.' 
% ....
% `As we have already said, a topological structure on a set enables one
% to give an exact meaning to the phrase ``whenever $x$ is sufficiently
% near $a$, $x$ has the property $P\{x\}$''. But, apart from the situation
% in which a ``distance'' has been defined, it is not clear what meaning
% ought to be given to the phrase ``every pair of points $x,y$ which are sufficiently near each other has 
% the property $P \{x, y\}$'', since a priori																																																																																																																																								
% we have no means of comparing the neighbourhoods of two different
% points. Now the notion of a pair of points near to each other arises fre-
% quently in classical analysis (for example, in propositions which involve
% uniform continuity). It is therefore important that we should be able
% to give a precise meaning to this notion in full generality, and we are
% thus led to define structures which are richer than topological structures,
% namely {\em uniform structures}.'
%\end{quote} 
%}
\begin{quote}
a topological structure on a set enables one
 to give an exact meaning to the phrase ``whenever $x$ is sufficiently
 near $a$, $x$ has the property $P\{x\}$''. But, apart from the situation
 in which a ``distance'' has been defined, it is not clear what meaning
 ought to be given to the phrase ``every pair of points $x,y$ which are sufficiently near each other has
 the property $P \{x, y\}$'', since a priori
 we have no means of comparing the neighbourhoods of two different
 points. %Now the notion of a pair of points near to each other arises fre-
% quently in classical analysis (for example, in propositions which involve
% uniform continuity). It is therefore important that we should be able
% to give a precise meaning to this notion in full generality, and we are
% thus led to define structures which are richer than topological structures,
% namely {\em uniform structures}.'
\end{quote}
In fact, a metric gives rise to a functor $FiniteSets\lra \Filt$, and the notion of a uniform structure 
is equivalent to such a functor %with an extra condition FIXME being kinda 2-dimensional.
satisfying the ``2-dimensionality'' Condition of Def.~\ref{defi:topoic_and_e-neighbs}(3). 
\end{rema}

\begin{exer}\label{uniform_simplicial} Show that  the definition of a uniform structure \href{http://mishap.sdf.org/tmp/Bourbaki_General_Topology.djvu#page=174}{[Bourbaki,General Topology, II\S1.1,Def.1]} 
in fact describes a 2-dimensional (i.e.~satisfying Def.~\ref{defi:topoic_and_e-neighbs}(3))
symmetric simplicial filter, where symmetric means that it factors as 
$\Dop \lra FiniteSets^{\text{op}} \lra \Filt$.
\end{exer}

\subsubsection{Simplicial neighbourhoods associated with metric maps on the large scale}

We shall now show that $\sFilt$ contains the opposite to the category of geodesic (on the large scale) metric spaces
with surjective quasi-Lipschitz maps.

Call a map $f:X\lra Y$ of metric spaces {\em injective on the large scale} [Gromov,\S0.2.D] iff 
there is a monotone unbounded real function $\lambda(d), d>0$, such that either of the %following 
equivalent conditions hold:
\bi\item
$\dist_X(x,y)\leq \lambda(\dist_Y(f(x),f(y)))$  for all $x,y\in X$
\item $\dist_Y(f(x),f(y))\geq \lambda^{-1} (\dist_X(x,y))$ for all $x,y\in X$
\ei

Associate with a metric space $X$ the sset $\Dop\lra Sets,\ n\longmapsto \homm{\Sets}n X$ represented
by the set of its points. Call a subset $U$  of $X^n$  a neighbourhood iff there is $d>0$ and a set $B$
of diameter $<d$ such that 
\bi\item $(x_0,...,x_n)\in U$ 
whenever $x_0,...,x_n\not\in B$ and \\ \ \ \ \ \ ${\tiny{}\ \ \ \ \ \ \ \ \ \ \ \ \ \ \ \ \ \ \ \ \ \ \ } \dist(x_i,x_j)>d$ for all $0\leq i<j\leq n$ such that $x_i\neq x_j$
\ei

Let $X_\LL:\sFilt$ be the simplicial neighbourhood obtained. 

\begin{exer} Check that $-_\LL$ defines a contravariant fully faithful embedding of the category of 
geodesic (on the large scale) metric spaces with surjective quasi-Lipschitz maps, 
cf.~[Gromov,\S0.2.D,\S0.2.A${}_2'$].
\bi \item A morphism $X_\LL\lra Y_\LL$ in $\sFilt$ is a map of metric spaces $X\lra Y$
injective on the large scale.
%such that 
%there is a real function $\lambda(d), d>0$, such that
%\bi
%\item $\dist_X(x,y) < \lambda (\dist_Y(f(x),f(y)))$ for all $x,y\in X$
%\ei
%\item This is equivalent to saying that $f^{-1}:Y\lra X$ is uniform on the large scale 
%whenever the map $f:X\lra Y$ 
%of sets is surjective and $f^{-1}:Y\lra X$ is well-defined. More generally, it is enough that the image $f(X)$
%forms a $d$-net for some $d>0$.  
\item If $X$ and $Y$ are geodesic and $f:X\lra Y$ is surjective, then $f^{-1}$ is quasi-Lipschitz 
iff $f:X_\LL\lra Y_\LL$ is well-defined, i.e.~iff $f$ is injective on the large scale.  
\item Conclude that  $-_\LL$ defines a contravariant fully faithful embedding of the category of 
geodesic (on the large scale) metric spaces with surjective quasi-Lipschitz maps.
\item Work out the geometric meaning of %Reformulate geometrically 
the lifting property defining connectedness %and compare to [Gromov,0.2$A_2'$]
(see Exercise~\ref{connected} for notation and explanations): 
\bi
\item  $X\lra \{0=1\}_\ttt \rtt  \{0,1\}_\ttt \lra\{0=1\}_\ttt$ 
\ei
\ei
\end{exer}

\begin{exertodo} In \S\ref{Skorokhod-Grayson}%{grayson} 
we define the notion of geometric realisation as an endofunctor of $\sFilt$. 
Define in $\sFilth$, and then work out the geometric meaning of, $X_\LL \lra B(G_\LL)$: does $B(G_\LL)$ classify something which 
may be called a $G$-bundle on the large scale ? 
\end{exertodo}

\subsection{The subdivision neighbourhood structure on a simplicial set}\label{subdivision-def}
%## Fine ... decomposition ... bacycentric ...

We give precise meaning to ``subdivide a simplex
into simplexes small enough'', as follows.

\begin{defi}%_Definition_. 
%Equip an X:sSet with simplicial filters
%  formed by all e-neighbourhoods for e:X_n,n>=0,
%  where _e-big_ means one of the following:
%
Let $X:\sSet$ and let $\epsilon\in X_k$ be a simplex, $m\geq 0$. 
A subset $\varepsilon\subset X_n$ is {\em  $(\epsilon,m)$-neighbourhood}, 
resp.~{\em  $(\epsilon,m)_>$-neighbourhood} or {\em  $(\epsilon,m)_<$-neighbourhood}, iff 
for each $N>0$ and each simplex $\epsilon':X_N$ 
such that $\epsilon$ is a face of $\epsilon'$ 
it holds:
\bi
\item $\epsilon'[t_1\!\leq...\leq\!t_n]\in \varepsilon$ whenever 
% $\varepsilon$ is a face of a simplex $\varepsilon':X_N$ and 
 $0\leq t_1\leq ... \leq t_n\leq t_1+m\leq N$ 
\item $\epsilon'[t_1\!\leq...\leq\!t_n]\in \varepsilon$ whenever 
% $\varepsilon$ is a face of a simplex $\varepsilon':X_N$ and 
 $0\leq t_1\leq ... \leq t_n\leq m\leq N$ 
\item $\epsilon'[t_1\!\leq...\leq\!t_n]\in \varepsilon$ whenever 
% $\varepsilon$ is a face of a simplex $\varepsilon':X_N$ and 
 $N-m<t_1\leq ... \leq t_n\leq N$ 
\ei

The {\em subdivision}, resp.~{\em $>$-subdivision} or {\em $<$-subdivision},  
filter on $X_n$ is generated by 
  the $(\epsilon,m)$-neighbourhoods where $\epsilon:X_k$ varies through simplices of arbitrary dimension  $k\geq 0$, $m>0$.

Let $X_{\text{subd}}$ denote the simplicial set equipped with the subdivision neighbourhood structure,
and similarly for $<$-subdivision and $>$-subdivision neighbourhood structures.
\end{defi}
%
%Equip an $X:\sSet$ with simplicial filters
%  generated by all $(\varepsilon,m)$-fine subsets for simplices $\varepsilon:X_k,k\geq 0$, $m>0$.
%
%  where {\em $\varepsilon$-fine} subset $U\subseteq X_m$ means one of the following:
%%
%%
%%
%% (m)   e'[i,...,i+m] in U whenever
%%     there is a map
%%      n--->k,   X_k--->X_n , e'|--->e
%%
%%
%% (t_<) e'[0,...,m] in U whenever
%%     there is a map
%%      n--->k,   X_k--->X_n , e'|--->e
%%
%% (t_>)  e'[n-m,...,k] in U whenever
%%     there is a map
%%      n--->k,   X_k--->X_n , e'|--->e
%%
%% (t) e'[0,...,m] and e'[n-m,...,k] in U whenever
%%     there is a map
%%      n--->k,   X_k--->X_n , e'|--->e
%%
%for each simplex $\varepsilon'$ such that $\varepsilon$ is a face of $\varepsilon'$, it holds
%\bi
%\item[ (m)]   $\varepsilon'[i,...,i+m] \in U$ whenever $0<i\leq i+m\leq n$ 
%\item[$(t_<)$] $\varepsilon'[1,...,m] \in U$ whenever $m\leq n$ 
%\item[$(t_>)$]
%$\varepsilon'[n-m,...,n] \in U$ whenever $m<n$ 
%\item[$(t_{<>})$] 
% $\varepsilon'[1,...,m] \in U$ whenever $m\leq n$ 
%and 
%$\varepsilon'[n-m,...,n] \in U$ whenever $m<n$ 
%\ei
%%and 
%%     $\varepsilon$ is a face of $\varepsilon'$, 
%%      i.e.~there is a map $k\lra n$ inducing $X_n\lra X_k$ sending $\varepsilon'\mapsto\varepsilon$
%\end{defi}
%
%
%\begin{exer}%\begin{exer}%_Exercise_.?_ 
%\bi\item Check that a surjective map of simplicial sets induces a map of the  corresponding simplicial filters.
%\item If the simplicial set is connected, then the filter defined above is formed, rather than generated, 
%by the $\varepsilon$-neighbourhoods.
%\ei 
%\end{exer}
%
%

\begin{exem}%\begin{exer}%\begin{exer}%_Exercise_.?_ 
\bi\item For the simplicial set of Cartesian powers $\homm\Sets-M,\ n \mapsto M^n$ of a set, 
the subdivision filter is trivial: $M^n$ is the only neighbourhood of $M^n$. 
\item (see details in the next subsection) For a linear order $[0,1]_\leq$, 
in the simplicial set consider ``co-represented'' by $I$
$$ n_\leq\longmapsto \homm{{\text{preorders}}}{n_\leq}{[0,1]_\leq},$$ the subdivision filter is generated by the sets
$\{(t_1\leq ...\leq t_n): \dist(t_i,t_{i+1})<\epsilon\}$ where $\epsilon>0$. Equivalently, this is  the filter is generated by the subsets of simplices of 
diameter $\epsilon$, for $\epsilon>0$.

%%not helpful%\item More generally, let $I$ be a partial preorder equipped with a metric such
%%not helpful%that for each $\varepsilon>0$ there is an increasing sequence 
%%not helpful% $t_1<t_2<...<t_n$ forming an $\varepsilon$-net, i.e.~such that for any $t\in I$ there is some $t_i,0<i\leq n$ such that 
%%not helpful%$\dist(t,t_i)<\varepsilon$. Then the notion of $(\varepsilon,m)$-fine sets is  generated by the metric, i.e.~a subset of 
%%not helpful%$Hom(n_\leq,I_\leq)$ is $(\varepsilon,m)$-fine for some simplex $\varepsilon$ and $m>0$ iff there is $\varepsilon>0$ such that it contains 
%%not helpful%all the simplices of diameter $<\varepsilon$. 
%%not helpful%
%%not helpful%\item For the partial order $I=[0,1]_\leq \times [0,1]_\leq$, 
%%not helpful%in the simplicial set ``co-represented'' by $I$
%%not helpful%%$$n\longmapsto Hom_{{\text{preorders}}}(n_\text{as preorder},(I,\leq))$$, 
%%not helpful%the filter is generated by the sets
%%not helpful%$\{(t_1\leq ...\leq t_n;s_1\leq ...\leq s_n): \dist(t_i,t_{i+1})<\varepsilon, \dist(s_i,s_{i+1})<\varepsilon \text{ for all }i\,\}$.
%%not helpful%
%%--%
%%--%\item For the alphabetical order $I=[0,1]\times [0,1]$ where $(a,b)\leq (a',b,)$ iff either $a\leq a'$ or $a=a'$ and $b\leq b'$, 
%%--%in the simplicial set ``co-represented'' by $I$
%%--%%$$n\longmapsto Hom_{{\text{preorders}}}(n_\text{as preorder},(I,\leq))$$, 
%%--%the filter is generated by the sets
%%--%$\{(t_1\leq ...\leq t_n;s_1\leq ...\leq s_n): \dist(t_i,t_{i+1})<\varepsilon, \dist(s_i,s_{i+1})<\varepsilon \text{ for all }i\,\}$.
%%--%???
%%--%

\ei 
\end{exem}%\end{exer}

\subsection{%# 
The real line interval $[0,1]$}\label{unit:sF} 
\def\Sing{\text{Sing\,}}
We now may define the interval object in $\sFilt$ by equipping the sset corepresented by a linear order with 
the subdivision neighbourhood structure.

\begin{defi} Let {\em the interval object $[0,1]_\leq$ in $\sFilt$} be the sset 
%Let $I$ be a preorder; the interesting case is $I = ( [0,1],\leq )$ % I = ( [0,1],<= )
%or another dense linear order. Associate with $I$ a simplicial set
%$\Sing I$ by taking the functor "almost" co-represented by $I$
%in the category of preorders
%
$$\Dop\lra \Sets$$ %  /_\^op  ---> Sets
$$n_\leq \longmapsto \homm{\leq}{n_\leq}{[0,1]_\leq}$$ %        n |---> Hom_preorders ( n, I )
equipped with the subdivision neighbourhood structure. 
Let  $[0,1]_{\leq^-}$,  $[0,1]_{\leq^+}$, and  $[0,1]_{\leq^\pm}$ a denote the same sset equipped with the
 $<$-subdivision, $>$-subdivision, and the union of  $<$-subdivision and $>$-subdivision neighbourhood structures, resp.
%%
%%Let $\Sing^mI, \Sing^{t_<} I, \Sing^{t_>} I, \Sing^{t_{<>}} I$ %Sing^mI, Sing^ts I, Sing^tt I, Sing^t I 
%%be  the sset equipped with the simplicial filters as
%%defined above.
%%
\end{defi}

We also use similar notation for any preorder $(I,\leq)$.

\begin{exer} The subdivision neighbourhood structure on  $[0,1]_\leq$ is induced by the metric in the following sense:
it is generated by the $\epsilon$-neighbourhoods of the diagonal $\{0\leq t_0\leq ... \leq t_n\leq 1: \dist(t_i,t_{i+1})<\epsilon\}$,
for $\epsilon>0$.
\end{exer}

\begin{exer} %\begin{exer}%_Exercise_.._ 
 Let $X$ be a topological space, $M$ a metric space.
\bi\item
 A continuous function $[0,1]\lra M$ %$[0,1]-->M$ 
is the same as a morphism $[0,1]_\leq \lra M_\mmU$ in $\sFilt$.
\item

 A continuous function $[0,1]\lra X$ %[0,1]-->X$[0,1]-->X$ 
is the same as a morphism $ [0,1]_{\leq^\pm} \lra X_\ttt$ in $\sFilt$.

% A continuous function $[0,1]-->M$ is the same as $Sing^t I --> M_\mmU$ in sFilth

\item
 A upper semi-continuous function $[0,1]\lra X$ %[0,1]-->X $[0,1]\lra X$ %[0,1]-->X 
 is the same as $[0,1]_{\leq^-} \lra X_\ttt$ in $\sFilt$.
\item
 A lower semi-continuous function $[0,1]\lra X$ % [0,1]-->X 
is the same as $[0,1]_{\leq^+} \lra  X_\ttt$ in $\sFilt$.
\ei
\end{exer}

\subsection{A notion of homotopy based on the interval $[0,1]_{\leq}$} 
In the standard way this notion of an interval leads to a notion of homotopy.
Note that we later define a notion of homotopy based on a notion of the mapping space,
which we feel is more appropriate.

\begin{rema}[Homotopy on $\sFilt$] The definition above lets us define a notion of homotopy in $\sFilt$ in the usual way: 
two maps $f,g:X\lra Y$ are {\em homotopic} 
iff %they both factor through the same map
there is a linear order $I$, elements $i_f$ and $i_g$, a morphism $ X\times I_\leq \xra h Y$ such that 
$f$ factors as  $X\times \{i_f\} \lra X\times I_\leq \xra h Y$  
and 
$g$ factors as  $X\times \{i_g\} \lra X\times I_\leq \xra h Y$.
\end{rema}

\begin{exertodo} Compare this to the notion of homotopy defined in \S\ref{homotopy-mapping}.
\end{exertodo}

\subsection{Mapping spaces, geometric realisation and a notion of homotopy}

\subsubsection{Discretised homotopies as Archimedean simplices}\label{discr_homotpies}

\begin{defi}
A simplex $s$ is {\em $\varepsilon/n$-fine} iff $s[t_1\!\leq\! ... \!\leq\! t_k]\in \varepsilon$ whenever $t_1\leq ... \leq t_k\leq t_1+n$.
A simplex $s$ is {\em Archimedean} iff it can be split into finitely many of arbitrarily small parts, 
i.e.~is a face of an $\varepsilon/n$-fine simplex for every neighbourhood $\varepsilon\subset X_k$ and every $k,n>0$.
Call such an  $\varepsilon/n$-fine simplex an  $\varepsilon/n$-refinement of $s$.\footnote
{The geometric intuition suggests this definition should possibly be modified:  $\varepsilon/n$-refinement of a simplex may ``go off to infinity''
as $\epsilon$ and $n$ vary; in the case of a metric space, the $\epsilon$-chains (i.e.~$\epsilon$-discretised homotopies) connecting two 
points go off to infinity rather than converge on an actual path). One may want to require something like that 
a simplex $s$ is {\em Archimedean} iff for  every neighbourhood $\varepsilon\subset X_k$ and every $k,n>0$ 
it has an  $\varepsilon/n$-refinement which in its $\varepsilon$-neighbourhood has 
$\delta/m$-refinement with the same property, for every  neighbourhood $\delta\subset X_{l}$ and every $l,m>0$. 
}

Call a set of simplices {\em bounded} iff its simplices are Archimedean and there is an upper bound on the dimension
of their $\varepsilon/n$-refinements for each $n$ and neighbourhood $\varepsilon$.
\end{defi}

\begin{exer} Check that for a metric space $M$, 
a pair of points $s=(x,y)\in M\times M$ form an Archimedean simplex in $M_\mU$ iff for each $\epsilon>0$ there is 
an $\epsilon$-{\em discretised path} 
$x=x_0,x_1,...,x_l=y$, $\dist(x_t,x_{t+1})<\epsilon$ for $0\leq t<l$, {\em from $x=(x,y)[0]$ to $y=(x,y)[1]$}. 

Note we do not require that these $\epsilon$-discretised path converge on an actual path; perhaps 
this hints the definition of an Archimedean simplex should be modified.

Check that in a metric space $M$ with an inner metric a subset $B\subset M\times M$ is bounded iff 
there is a bound on the distance between points for $(x,y)\in B$, i.e.~there is $d>0$ such that 
$\dist(x,y)\leq d$ whenever $(x,y)\in B$.
\end{exer}

%In other words, 
%an Archimedean simplex in $M\times M$ is a pair of homotopic points, or rather a pair of points
%which can be connected by a {\em discretised homotopy} of arbitrary finess.

%Note that if $M$ is the space of functions, 

\begin{exer} 
Archimedean simplices form a subobject (subfunctor) $X_{\Arch}$ of $X$, for $X:\sFilt$. 
Check that a surjective map $X\lra Y$ induces a map $X_\Arch\lra Y_\Arch$.
%Two functions $f,g:A\lra M$ are homotopic iff $(f,g)$ is an Archimedean simplex 
\end{exer}

\begin{exer}\label{e-chain-of-functions}
A well-known lemma says that two functions $f,g:A\lra M$ from an arbitrary topological space $A$ to a metric space $M$
are homotopic iff there is a $\epsilon$-discretised homotopy 
$f=f_0,..., f_n=g$ such that for any $x\in A$ $\dist(f_t(x), f_{t+1}(x)<\epsilon$, under some assumptions on 
the metric space $M$; it is enough to assume that for every $\epsilon>0$ there is $\delta>0$ such that 
every $\epsilon$-ball contains a contractible subset containing a $\delta$-ball with the same centre. (Ref!) 

Check that the lemma says that two functions $f,g:A\lra M$ are homotopic iff $(f,g)$ is an Archimedean simplex of the mapping space $Func(A,M)$ 
with the sup-metric, or, equivalently, iff $(f,g)$ is an Archimedean simplex of the mapping space $(Func(A,M)_\mU)_{\Arch}$.
\end{exer}

\subsubsection{Skorokhod mapping spaces: the definition}\label{mapping_def}\label{Skorokhod-mapping-spaces-defs} 
We now repeat somewhat more formally the definition of the Skorokhod filters given in \S\ref{skorokhod-filter-def}.
\begin{defi}
 Let $X,Y:\sFilth$ be objects of $\sFilth$. For $N>2n$, $\delta\subset X_N$ and $\varepsilon\subset Y_n$, 
a {\em $\varepsilon\delta$-Skorokhod neighbourhood of the Hom-set $\homm\sSets  {X_\text{as sSet}} {Y_\text{as sSet}}$} of underlying simplicial sets
of $X$ and $Y$ is the subset %of $\Hom X Y$ 
consisting 
of all the functions $f:X_\text{as sSet}\lra Y_\text{as sSet}$ with the following property:
\bi\item[] there is a neighbourhood $\delta_0\subset X_n$ such that 
each ``$\delta_0$-small'' $x\in \delta_0$ has a ``$\delta$-small''
``continuation'' $x'\in X_N$, $x=x'[1..N]$ such that its ``tail'' maps into something ``$\varepsilon$-small'', 
i.e.~$f(x'[N-n+1..N])\in\varepsilon$.
\item[] As a formula, this is 
$$\{ f:X\lra Y\,:\, \exists \delta_0 \subset X_n\, \forall x\in\delta_0 \, \exists x'\in \delta 
( x=x'[1...n] \,\,\,\&\,\,\, f(x'[N-n+1,...,N])\in \varepsilon ) \}$$

\ei
The {\em Skorokhod filter on $\homm\sSets {X_\text{as sSet}} {Y_\text{as sSet}}$} is the filter generated by all the Skorokhod $\varepsilon\delta$-neighbourhoods for 
$N\geq 2n>0$ (sic!), neighbourhoods $\delta\subset X_N$ and $\varepsilon\subset Y_n$. 

As a formula, it is 
$$\{ f:X\lra Y\,:\, \exists \delta_0 \subset X_n\, \forall x\in\delta_0 \, \exists x'\in \delta 
( x=x'[1...n] \,\,\,\&\,\,\, f(x'[N-n+1,...,N])\in \varepsilon ) \}$$

Let $\Homm\Filt X Y$ denote $\homm\sSets {X_\text{as sSet}} {Y_\text{as sSet}}$ equipped with the Skorokhod neighbourhood structure.
\end{defi}

This allows to  define mapping spaces in $\sFilt$ by equipping the inner Hom of ssets with Skorokhod filters.
%\begin{defi}[Skorokhod mapping space]ooo The {\em Skorokhod mapping space $\Homm\sFilth X Y$} is the simplicial set $\Homm\sSets XY$ equipped 
\begin{defi}[Mapping space]\label{mapping_space} 
 The {\em Skorokhod mapping space $\Homm\sFilth X Y$} is the inner Hom $\Homm\sSets {X_\text{as sSet}} {Y_\text{as sSet}}$ of the underlying
simplicial sets of $X$ and $Y$ equipped
with the neighbourhood structure as follows. Equip the $(n-1)$-simplex $\Delta_{n-1}=\homm\leq-{n_\leq}$ with the antidiscrete filter, 
equip  $X\times \Delta_{n-1}=X\times \homm\leq-{n_\leq}$ with the product filter, and, finally, equip the set of $(n-1)$-simplicies  
$\homm\sSets{{X_\text{as sSet}} \times \homm{{\text{preorders}}}-{n_\leq}}{{Y_\text{as sSet}}}$ with the resulting Skorokhod neighbourhood structure. %filter. 
\end{defi}

\begin{rema}\label{rema:mapping}
%Now recall the construction of inner hom in ssets. The following definition
%is made in a way so that face and degeneration maps are continuous because
%of simplicial combinatorics:
%define a neighbourhood $ U=U_{\varepsilon,\delta, t_0\leq ... \leq t_n, s_0\leq ... \leq s_n} \subset Hom(X\times Hom(-,n), Y)$ 
It is possible to define a different neighbourhood structure.
We give this definition to demonstrate that there are implicit choices made in formalisation,
and care needs to be taken. Another reason is that we are not sure our definitions are most appropriate.

Call $U\subset \homm\sSets{X\times \homm\leq-{n+1_\leq} } Y$
a neighbourhood iff 
for some $m\geq 0$, neighbourhoods $\delta\subset X_m$ and $\varepsilon\subset Y_m$,
some (equivalently, all) sequences $\overrightarrow t:= (0\leq t_0\leq ... \leq t_m\leq n)$, $\overrightarrow s := (0\leq s_0\leq ... \leq s_m \leq n)$ such that
$t_0\leq s_0$, $t_1\leq s_1$,..., $t_m\leq s_m$, 
it holds 
\bi\item[] $f\in U$ whenever 
 the following implication holds
$$
f_n[\overrightarrow t]\big( (\delta \times \theta) [\overrightarrow t] \big) \subset \varepsilon[\overrightarrow t]  \implies 
f_n[\overrightarrow s]\big( (\delta \times \theta) [\overrightarrow s] \big) \subset \varepsilon[\overrightarrow s] $$
\ei
For example, if $n=0$ the only neighbourhood is $U=\homm\sSets{X\times \homm\leq-{1_\leq} } Y$ itself.
\end{rema}

\begin{exertodo}\label{mapping_top}
Let $X$ be a locally compact space, and let $M$ be a metric space such that
each ball contains a contractible subset containing an open ball around the same point.
\bi\item Calculate the Skorokhod mapping space $\Homm\sFilth {X_\ttt}{M_\mU}$.
%%incorrect...%\bi
%%incorrect...%\item Check whether the two definitions of the mapping space  $\HHom_{\sFilt}(X_\ttt,M_\mU)$ coincide in this case. 
%%incorrect...%\item 
%%incorrect...%Check whether  $\HHom_{\sFilt}(X_\ttt,M_\mU)_{\ttt^*}$ is the space of continuous maps from $X$ to $Y$
%%incorrect...%equipped with compact-open topology. 
%%incorrect...%\item Let $X$ be compact. Check whether  $\HHom_{\sFilt}(X_\ttt,M_\mU)_{\mU^*}$ is the space of continuous maps from $X$ to $Y$
%%incorrect...%equipped with sup-metric $\dist(f,g)=\sup_{x\in X} \dist(f(x),g(x))$. 
%%incorrect...
\item Check whether $f:X\lra M$ is homotopic to $g:X\lra M$ iff the simplex
$(f,g)$ is Archimedean in $\Homm{\sFilt}{X_\ttt}{M_\mU}$.%_{\ttt^\inv}$.
\ei
\end{exertodo}

\begin{exertodo} Check the following basic properties of the evaluation map for Skorokhod spaces.
\bi\item Check whether 
%the evaluation map in $sSets$ gives rise to a isomorphism %bijection 
%natural in $K,X,Y$:
the isomorphism of the underlying ssets (natural in $A,X,Y$) is necessarily continuous:
\begin{equation*}
\begin{split}
\underline{ev}_*:\Homm\sFilt A{\Homm\sFilt X Y}  \lra \Homm\sFilt {A\times X} Y
%\\ X\times K \xra { \id_X\times f} X \times  \Hom_\sSet(X_.,Y_.) \xra {ev} Y_.
\end{split}
\end{equation*}

\item Check whether for the initial object  $\Delta_0= \homm\leq{n_\leq}{1_\leq}$ of $\sFilth$, 
the mapping space from $\Delta_0$ to an $Y:\sFilth$ is $Y$ itself:
\begin{equation*}
\begin{split}
\underline{ev}_*:\Homm\sFilt {\homm\leq{n_\leq}{1_\leq}} Y   \xra{(iso)} Y 
%\\ X\times K \xra { \id_X\times f} X \times  \Hom_\sSet(X_.,Y_.) \xra {ev} Y_.
\end{split}
\end{equation*}
% seems it is! Find sufficient conditions if not. 
\ei\end{exertodo} 

\begin{exer} Check whether $\Homm\sFilt - Y$ has an ``inner'' analogue of the left adjoint in the following sense.
\bi
\item Define a semi-direct product $A\ltimes X$ in $\sFilth$ as follows. 
The underlying simplicial set is that of the direct product in $\sSets$, or, equivalently, in $\sFilth$: 
$(A\ltimes X)_n:=(A\times X)_n$. The filter is defined as follows: 
$\delta\subset A_n\times X_n$ is a neighbourhood 
iff there is a neighbourhood $\alpha\subset A_n$ such that for each ``$\alpha$-small'' $a\in \alpha$ 
there is a neighbourhood $\delta_a\subset X_n$ such that 
$(a,x)\in \delta$ whenever $a\in\alpha$ and $x\in\delta_a$. 
\item Check whether there is an isomorphism in $\sFilth$ natural in $A,X,Y$:
\begin{equation*}
\begin{split}
\underline{ev}_*:\Homm\sFilt A{\Homm\sFilt X Y}  \xra{(iso)} \Homm\sFilt {A\ltimes X} Y
%\\ X\times K \xra { \id_X\times f} X \times  \Hom_\sSet(X_.,Y_.) \xra {ev} Y_.
\end{split}
\end{equation*}
\item Ponder the similarity to the definition of the topoic filter in the
definition of the $\sFilth$-neighbourhood structure
associated with topological spaces.

\ei

\end{exer}

%%b4_realised_connection_to_Skorokhod%\begin{exertodo}  Does the product in $\sFilt$ have a right adjoint? 
%%b4_realised_connection_to_Skorokhod%\bi\item Check whether it is possible to replace in Definition~\ref{mapping_space} 
%%b4_realised_connection_to_Skorokhod%by arbitrary continuous function 
%%b4_realised_connection_to_Skorokhod%the ``forgetful'' coordinate projections $x=x'[\overrightarrow t]$ 
%%b4_realised_connection_to_Skorokhod%and the function $f$ in the assumption  $f(x_l, l\! \xra{\theta}\! m) \in\delta[\overrightarrow t]$.
%%b4_realised_connection_to_Skorokhod%\ei\end{exertodo}
%%b4_realised_connection_to_Skorokhod%

%%b4_realised_connection_to_Skorokhod%
%%b4_realised_connection_to_Skorokhod%\begin{exertodo} 
%%b4_realised_connection_to_Skorokhod%Does  the mapping space functor as defined in Definition~\ref{mapping_space} have a left adjoint?
%%b4_realised_connection_to_Skorokhod%\end{exertodo}
%%b4_realised_connection_to_Skorokhod%

\subsubsection{A notion of homotopy based on the mapping space}\label{homotopy-mapping}
Now that we have a notion of a mapping space, we define a notion of homotopy using Archimedean simplices. %\S\ref{discr_homotpies}
Recall that a pair of continuous functions $f_0,f_1:X\lra Y$ defines a map 
$X\times \Delta_1\lra Y$, $(s,\theta)\mapsto (f_{\theta(0)}(s),...,f_{\theta(n)}(s))\in Y^{n+1}=Y_n$ and thus 
 can be viewed as a $1$-simplex in $\Homm{\sFilt}XY$.
\begin{defi}[Skorokhod homotopic]
Say that {\em a map $f:X\lra Y$ is {\em Skorokhod homotopic} to $g:X\lra Y$} iff
%Check whether $f:X\lra M$ is homotopic to $g:X\lra M$ iff
 the simplex $(f,g)$ is Archimedean in $\Homm{\sFilt}XY$.
\end{defi}

\begin{exertodo} Is this notion symmetric? transitive? Study this notion of homotopy.
For example:
\bi\item Let $X,Y:sSets$ be ssets equipped with antidiscrete filters: the only neighbourhood is the whole set. 
Check that this is equivalent to the definition of simplicial homotopy in [Goerss-Jardine, I\S6].
\item Check whether this gives the standard notion of homotopy of maps from topological spaces to metric spaces.
Compare~Exercise~\ref{e-chain-of-functions} and Exercise~\ref{mapping_top}.
\ei\end{exertodo}

\subsubsection{The geometric realisation via the approach of Grayson}\label{grayson} % and Besser and Drinfeld}
We now repeat more formally the construction of the geometric realisation due to [Besser], [Grayson, Remark 2.4.1-2], and [Drinfeld]
sketched in \S\ref{Skorokhod-paths-space}.

[Besser, Def.3.3] and [Grayson, Remark 2.4.1-2] give a construction of the geometric realisation 
of a simplicial set based on the observation 
 that the standard simplex $\Delta_n=\{ (s_1,..,s_n) \in [0,1]^n : 0\leq s_1\leq ... \leq s_n\leq 1\}$
is the space of maps $[0,1]_\leq \lra (n+1)_\leq$ of preorders modulo some identifications,
i.e.
$$\Homm{\sSets}{\homm{\leq}-{[0,1]_\leq}} { \homm{\leq}-{(n+1)_\leq}}.$$
The notion of a mapping space in $\sFilth$ suggests we should try to define the geometric realisation of a simplicial set $X$
as the Skorokhod  space $\Homm\sFilt{\Homm{preorders}-{[0,1]_\leq}} X$ of (discontinuous) paths in $X$ equipped 
with an appropriate neighbourhood structure. 

\begin{defi} The {\em Besser %-Drinfeld-Grayson 
geometric realisation} of $X:\sFilt$ is 
the endofunctor $$\sFilt\lra\sFilt,\ \ X\longmapsto \Homm{\sFilt}{[0,1]_\leq} X$$

The {\em Grayson subdivision} is the endofunctor 
$$\sFilt\lra\sFilt,\ \ X\longmapsto X\circ e$$
where $e:\Dop\lra\Dop$ is the endofunctor defined following [Grayson,\S3.1, esp.~Def.3.1.4, Def.3.1.8]): 
$$n\longmapsto 2n, \ f:m\ra n \longmapsto \{ n+i\mapsto n+f(i),\ n-i \mapsto n-f(i) , \text{ for }i=0,...,n-1\}$$ 
\end{defi}

\begin{exertodo} Verify details of the argument in  \S\ref{Skorokhod-paths-space} and prove the following.
\bi
\item Verify \S\ref{Skorokhod-paths-space} gives a well-defined map of ssets
 $$ |X|\lra \homm{sSets}{[0,1]_\leq}{X} 
 $$
\item Let $X_\diag$ denote the simplicial set $X$ equipped with the finest neighbourhood structure
such that the set $X_0$ of $0$-simplicies is antidiscrete. 
Explicitly, a subset of $X_n$ is a neighbourhood iff it contains the diagonal, i.e.~the image of 
 $X_0$ in $X_n$ under the unique degeneracy map. Verify that $X_\diag:\sFilth$ is well-defined.

\item Verify that for $X=\Delta_{n-1}=\homm{\leq}{\cdot}{n_\leq}$, $n>0$,  
the Hausdorffisation of the topologisation of the  Skorokhod paths space %$\sFilt$-geometric realisation 
is the standard simplex:
$$ \left(\Homm{\sFilt}{[0,1]_\leq}{(\Delta_n)_\diag}_{\ttt\inv}\right)_\text{Hausdorff} \xra{(iso)} |\Delta_n|$$

\item Verify that $\Homm{\sSets}{[0,1]_\leq} -:\sSets\lra\sSets$ 
preserves finite directed limits, and is also compatible with Skorokhod neighbourhood structure, 
i.e.~that $\Homm{\sFilth}{[0,1]_\leq} -:\sFilth\lra\sFilth$ also preserves finite directed limits.

\item (todo) Conclude that for a finite simplicial set $X$ there is an isomorphism of the geometric realisation and 
the Hausdorffisation of the topologisation of the Skorokhod paths space: 
$$ |X| \xra{(iso)} \left(\Homm{\sFilt}{[0,1]_\leq}{X_\diag}_{\ttt\inv}\right)_\text{Hausdorff} $$

\item (todo) Give a precise meaning to the following argument. For every 
$n>0$ 
%neighbourhood $\varepsilon\subset \Hom{n_\leq}{[0,1]_\leq} $
the sequence of $x_\theta$ is 
determined, up to $\varepsilon=1/n$, by $x_\theta$ where $\theta:n_\leq \lra [0,1], \theta=(0<1/n<...<(n-1)/n)$. 
Hence, the topological geometric realisation of a simplicial set  $X$ 
is dense in the topologisation of its Skorokhod paths space.
\ei
\end{exertodo}

\begin{rema} A.Smirnov suggested it maybe worthwhile to see whether the use of geometric realisation by [Suslin, On the K-theory of local fields]
%Suslin, A. A. (1984). On the K-theory of local fields. Journal of Pure and Applied Algebra, 34(2-3), 301–318. doi:10.1016/0022-4049(84)90043-4 
can be interpreted in terms of $\sFilt$.
\end{rema}

%\begin{exertodo} A $G$-bundle $\tilde X\lra X$ over a space $X$ 
%\end{exertodo}

\subsection{%### 
The connected components functor  $\pi_0$ as M2(l-lr)-replacement.} 
We observe that the connected components functor  $\pi_0$  is analogues to the (co)fibrant replacement 
postulated by Axiom M2 of model categories where the (co)fibrant replacement is taken with respect
to a morphism implicitly appearing in the definition of connectivity.

Recall that $ \{0,1\}\lra\{0=1\}$ %{0,1}-->{0=1}
 denotes the map of topological spaces gluing together the points of 
the discrete space with two points. As usual, we denote by 
$ \{0,1\}_\ttt\lra \{0=1\}_\ttt$ the corresponding map in $\sFilth$. 

Recall that $ \{0,1\}_\ttt$ can be explicitly described as follows: $n \lra \{0,1\}^n$,
and a subset of $ \{0,1\}^n$ is a neighbourhood iff it contains the diagonal $\{(0,..,0),(1,...,1)\}$.

\begin{exer}\label{connected} Check the following. 
\bi\item A topological space $X$ is connected iff $X\lra \{0=1\} \rtt  \{0,1\}\lra\{0=1\}$. 
\item A simplicial set is connected iff $X\lra \{0=1\}_\ttt \rtt  \{0,1\}_\ttt \lra\{0=1\}_\ttt$.
\ei
\end{exer}

Denote by $P^\lrl$ and $P^\rlr$ the classes (properties) of morphisms defined with respect to the left, resp.~right, lifting property:
$$P^\lrl:=\{ f\rtt g: g\in P\} \ \ \ \ P^\rlr := \{ f\rtt g: f\in P\}$$
It is convenient to refer to $P^\lrl$ and $P^\rlr$ as the property of {\em left, resp.~right, Quillen negation of property $P$}.

\begin{defi}%_Definition_. 

        In $\Top$, let $\pi_0$ be the functor defined by 
         the following $M2(l-lr)$ decomposition:
$$        
X\xra{(\{0,1\}\lra \{0=1\})^\lrl} \pi_0(X)\xra{(\{0,1\}\lra \{0=1\})^\lr} \{0=1\} $$
 %       X---({0,1}-->{0=1})^l--->\pi_0(X)---({0,1}-->{0=1})^lr-->{0=1}    
 
        In $\sFilth$, let $\pi_0$ be the functor defined 
         the following $M2(l-lr)$ decomposition:
$$        
X\xra{(\{0,1\}_\ttt\lra \{0=1\}_\ttt)^\lrl} \pi_0(X)\xra{(\{0,1\}_\ttt\lra \{0=1\}_\ttt)^\lr} \{0=1\}_\ttt $$        
%        X---({0,1}.-->{0=1}.)^l--->\pi_0(X)---({0,1}.-->{0=1}.)^lr-->{0=1}. 
\end{defi}

% 
%\begin{exer}%_Exercise._ 
%Check whether in $\sFilth$ for a topological space $X$ it holds that
%$\pi_0^{\sFilth}(X_\ttt)=(\pi_0^\Top(X))_\ttt.
%$
%Is this functor well-defined, or only up to non-canonical 
%        isomorphism or something ? Probably there is 
%        a universal functor among those of this form...
%\end{exer}
%

\begin{exer}%_Exercise._ 
     Check this definition is consistent with the usual definition of $\pi_0$ in $\sSets$. 
     Recall there is an embedding $\sSets\lra\sFilth, n\longmapsto (X_n)_{\text{antidiscrete}}$ which sends  a simplicial set $X$
into itself equipped the antidiscrete filters, i.e.~the filter such that the only neighbourhood is the whole set. 

Check that the set of 1-simplices $\pi_0^\sFilth(X)_1$ is the set of connected components of $X_{\text{antidiscrete}}$, and, in fact,  
$\pi_0^\sFilth(X)=\diag(\pi_0^\sSets(X)_{\text{antidiscrete}})$.
\end{exer}

\begin{exer}%_Exercise._ 
     Check this is consistent with the usual definition (notation) 
     of $\pi_0$ on $\Top$. For example, check the following.

     Let $X$ be a topological space such that the $\pi_0(X)$ is 
     well-defined (behaved), e.g.~$X$ has finitely many connected components.
     Then in $\Top$ it holds $\pi_0^{\sFilth}(X_\ttt)=(\pi_0^\Top(X))_\ttt$, or, in other words,
$$        
X_\ttt \xra{(\{0,1\}_\ttt\lra \{0=1\}_\ttt)^\lrl} \pi_0^\Top(X)_\ttt \xra{(\{0,1\}_\ttt\lra \{0=1\}_\ttt)^\lr} \{0=1\}_\ttt $$        
%        X---({0,1}.-->{0=1}.)^l--->\pi_0(X)---({0,1}.-->{0=1}.)^lr-->{0=1}. 
%     X---({0,1}-->{0=1})^l--->\pi_0(X)---({0,1}-->{0=1})^lr-->{0=1} 
\end{exer}

\subsection{Locally trivial bundles}\label{def-bundle} Here we repeat somewhat more formally \S\ref{def-bundle-sample} 
about local trivial bundles. 

It is said that being locally trivial means being locally a direct product. The precise meaning of this phrase
in terms of $\sFilth$ is straightforward: a map over a base $B$ 
is locally trivial iff it becomes a direct product after pullback along  $B[+1]\lra B$.  
We state this in the next $\S$ and then speculate whether this observation can be used to define a model structure on $\sFilth$. 
%This section uses notation introduced later  in \S\ref{embeddings}-\ref{shift}. 

\subsubsection{Local triviality as being a product after pullback along $B[+1]\lra B$}
%\begin{minipage}[c]{0.1458980337\textwidth}
%{\small $ \ \ \xymatrix{ X[+1] \ar[r]|--f \ar@{->}[d]|p &  B[+1] \times\{\NN\}\ar[d]|{} \\ B[+1] \ar[r]|-{\id} & B[+1]}$
%$ \ \ \xymatrix{  B[+1]\times_B X  \ar[r]|--f \ar@{->}[d]|p &  B[+1]\times F \ar[d]|{} \\ B[+1] \ar[r]|-{\id} & B[+1]}$
%}%\end{minipage}
%\begin{minipage}[c]{0.85841019663\textwidth}
\begin{exer} A map $p:X\lra B$ of topological spaces is a locally trivial bundle with fibre $F$ 
iff there is an $\sFilt$-isomorphism 
$$\tau:B_\ttt[+1]\times_{B_\ttt} X_\ttt \xra{(iso)} B_\ttt[+1]\times F_\ttt\text{ over } B_\ttt[+1]$$ 
%i.e.~such that $ p_\ttt[+1] = \tau\circ \pr_{B_\ttt}[+1] : B_\ttt[+1]\times_{B_\ttt} X_\ttt \lra B_\ttt[+1]$. 
i.e.~there is a commutative diagram as shown
$$ \ \ \xymatrix{   B_\ttt[+1]\times F_\ttt \ar[d]|{}  \ar[r]|--{(iso)} &  B_\ttt[+1]\times_{B_\ttt} X_\ttt \ar[d] \ar[r]|--{} & X_\ttt\ar@{->}[d]|p 
\\ B_\ttt[+1] \ar[r]|-{\id} & B_\ttt[+1] \ar[r] & B_\ttt 
} $$

%false,it_seems%The same holds for metric spaces: a map $p:X\lra B$ of metric spaces is a locally trivial bundle with fibre $f$ 
%false,it_seems%iff there is an $\sFilt$-isomorphism $$\tau:B_\mU[+1]\times_{B_\mU} X_\mU \xra{(iso)} B_\mU[+1]\times F_\mU\text{ over } B_\mU[+1]$$
%false,it_seems%i.e.~such that $ p_\mU[+1] = \tau\circ \pr_{B_\mU[+1]} : B_\mU[+1]\times_{B_\mU} X \lra B_\mU[+1]$. 
%false,it_seems%

Verify the argument in \S\ref{def-bundle-sample} using the following steps.
\bi
\item As simplicial sets, $B_\ttt[+1]\times_{B_\ttt} X_\ttt=\sqcup_{b\in B} X_\ttt$ and $ B_\ttt[+1]\times F_\ttt=\sqcup_{b\in B} B_\ttt\times F_\ttt$
where $\sqcup$ denotes disjoint union.
\item  To give a map of these simplicial sets over $B_\ttt[+1]$ is to give for each $b\in B$ a map of sets $f_b:X\lra  B\times F$ over $B$
which extends to a map of the corresponding simplicial sets $X_\ttt\lra B_\ttt\times F_\ttt$.
% as a connected component has to to a connected component.
\item The maps $f_b:X\lra  B\times F$, $b\in B$ represent a continuous $\sFilt$-isomorphism 
$B[+1]\times_B X \xra{(iso)} B[+1]\times F$ in $\sFilt$ iff for every $b\in B$ there is a neighbourhood $U_b\ni B$
such that %$f_b_{|p^{-1}(U_b)}:p^{-1}(U_b)\xra{(iso)}
$f_b$ defines a homeomorphism between $p^{-1}(U_b)$ and $U_b\times F$.
%this is false% \item Verify the same steps for the embedding $\mU$ of metric spaces  instead of the embedding $\ttt$ of topological spaces.
\item Check whether the above holds for the $\mU$-embedding of metric spaces.
\item Check whether  the above also holds for the localised category $\sFFilt$. 
\ei

\end{exer}

\begin{exertodo} Use the reformulation above to rewrite for $\sFilth$ a definition of the long exact sequence of a (co)fibration
in terms of the endofunctor $[+1]:\sFilth\lra \sFilth$ and base change $B[+1]\lra B$.
\end{exertodo}
\subsubsection{Suggestions towards a model structure on the category of simplicial filters}

\begin{exertodo} Can this notion of local triviality be used to define a model structure on $\sFilt$? For example, do the following. 
\bi
\item Calculate in $\sFilt$ 
$$(f):=\{X\xra p B: \text{in }\sFFilt\text{ }B[+1]\times_B X\lra B[+1]\text{ is of form } B[+1]\times F\lra B[+1]\}^\lr$$
\item Is it true that maps in $(f)$ have the homotopy extension property, i.e.~for any $A$ $A \lra A\times [0,1]_\leq \rtt (f)$ ? 
\item Does this define the class of fibrations in the category of topological spaces, i.e.~is it true under suitable assumptions 
that a map $p$ of topological spaces is a fibration 
iff $p_\ttt\in (f)$ ? 
\item(todo) Is there a model structure on $\sFilt$ where $(wc):=(f)^\lrl$ and  $(f)$ are the classes of weak cofibrations and fibrations, resp.?  
\item(todo)
More generally, define a model structure on $\sFilt$ or $\sFFilt$.
\ei
\end{exertodo}

\subsection{Taking limits of sequences and filters.}\label{limits}\label{limits-defs} %, compactness and completeness in terms of orthogonality}
Here we show how to reformulate about taking limits, convergence of sequences, equicontinuous families of functions, 
Arzela-Ascoli theorems, compactness and completeness,
with help of the  endofunctor $[+1]:\sFilth\lra \sFilth$ ``shifting dimension'' and Quillen lifting properties.
 
%\subsection{% ## Three embeddings const, diag, cart: Filth ---> sFilth
\subsubsection{Three embeddings $const, diag, cart$ of filters into simplicial filters} %: $\Filth \lra \sFilth$}
\label{embeddings}
Let $F$ be a filter. There are two natural ways to equip $\homm\Sets  n  F $ with a filter:
\bi
\item[(a)] the finest filter such that degeneracy (diagonal) map 
   $$\homm\Sets{1}{F}=F \lra \homm{\Sets}{ n}{ F }=F^n , x \longmapsto (x,x,...,x)$$
   %Hom(1,F)=F ---> Hom_Sets ( n, F )=F^n , x |--> (x,x,...,x)
   is continuous. 

\item[(b)] the coarsest filter such that all the face ``coordinate projection'' maps 
  $$\homm{\Sets}{ n }{F }=F^n \lra  \homm\Sets{1}{F}=F,\ (x_1,...,x_n)\longmapsto x_i, 0<i\leq n$$ 
%  Hom_Sets ( n, F )=F^n ---> Hom(1,F)=F 
  are continuous. 
\ei

\begin{exer}%\begin{exer}%_Exercise_.._ 
Explicitly, these filters can be defined as:
\bi
\item[ (a')] a subset of $F^n$ is a $diag$-neighbourhood iff it contains $\{(x,x,..,x) : x \in U \}$ for some $U$ an $F$-neighbourhood

\item[ (b')] a subset of $F^n$ is a $cart$-neighbourhood iff it contains $U_1 \times U_2 \times .. \times U_n$, for some $U_1,..,U_n$ $F$-neighbourhoods
\ei
\end{exer}
There is also the constant functor $\const:\Filt\lra \sFilt, n\mapsto F$.

\begin{exer}%\begin{exer}%_Exercise_.._ 
Check these define three  fully faithful embeddings  
$\diag,\cart,\const:\Filth\lra \sFilth$
%diag, cart : Filth ---> sFilth
\end{exer}

%__todo: define const: Filth ---> sFilth

%\begin{exer}%_Exercise_.._ Check const: Filth ---> sFilth is a fully faithful functor

%## Shift [+1]: /_\ ---> /_\ 
\subsection{Shift endofunctor $[+1]: \Delta\lra\Delta$}
\label{shift}

%__todo: define [+1]
\def\cof{\text{cof}}
\def\lr{\text{lr}}
\def\rl{\text{rl}}
The reformulation of the definition of limit in $\sFilth$ uses
the following ``shift'' endofunctor of $\Dop$ ``forgetting the first coordinate''.
Let $[+1]:\Delta\lra\Delta$ denote the shift by $1$ adding a new minimal element:
 $$[+1]: (1\!<\! ... \!<\!n)\longmapsto (-\infty\!<\!1\!<\!...\!<\!n) ; \ \  ([+1] f) (-\infty):=-\infty;\ \, ([+1]f)(i):=i$$
% $$[+1]: n\longmapsto n+1;\ \  n\xra\theta m \longmapsto ( \theta[+1] : 1\mapsto 1, i+1\mapsto \theta(i)+1). $$
The endofunctor is equipped with a natural transformation $[-1]\,:\,[+1]\implies\id:\Dop\lra\Dop$,
and its morphisms $X[+1]\xra{[-1]} X$ are particularly useful to us.

Now we recast in these terms a number of familiar notions in analysis: limit, Cauchy sequence, etc.

\subsubsection{%## 
Cauchy sequences and their limits.
}
A sequence $(a_n)$ of points of a metric space may be viewed as a map $\NN\lra M$. 
Let $\NN_\cof$ be $\NN$ equipped with the filter of cofinite subsets: a subset $U$ of $\NN$ is {\em a neighbourhood} 
iff 
 there is $N>0$ such that $m \in U$ whenever $m>N$. 

\begin{exer}%_Exercise_.._ 
The sequence $(a_n) \in M$ is Cauchy iff it determines a continuous map 
$$
  \homm\Sets-{\NN_\cof}_\cart \xra{} M_\mmU ,\ \ (i_1,..,i_n)\longmapsto (a_{i_1},...,a_{i_n})
$$
   The sequence {\em converges} iff this map factors as $$ \homm\Sets-{\NN_\cof}_\cart \lra M_\mmU[+1]\xra{[-1]} M_\mmU,$$
and the map is necessarily of form $ (i_1,..,i_n)\longmapsto (a_\infty, a_{i_1},...,a_{i_n}) $
where $a_\infty$ is the limit of the sequence.
%The limit of the sequence is given by the first coordinate: $$\homm\Sets-{\NN_\cof}_\cart \lra  M_\mmU[+1]\xra{(pr_0)}M$$. 
\end{exer}

\begin{exer}%_Exercise_.._ 
   A filter $F$ on $M$ is Cauchy iff $\homm\Sets-{F}_\cart \lra  M_\mmU$ is well-defined in $\sFilth$. 
  % (Check [Bourbaki,...] for definitions, also end-note??). 

   A Cauchy filter converges on $M$ iff the morphism $\homm\Sets-{F}_\cart\lra M_\mmU$ factors as 
    $F_\cart \lra  M_\mmU[+1]\lra M_\mmU $.

   A filter $F$ converges on $M$  iff the morphism $ \homm\Sets-{F}_\diag\lra M_\mmU$ factors as 
     $F_\diag \lra  M_\mmU[+1]\lra M_\mmU$. 
\end{exer}

\subsubsection{%## 
Limits on topological spaces. 
}

The same construction works for topological spaces.

Let $F$ be a filter on the set of points of a topological space $T$. 

The inclusion (equality) of underlying subsets $F \subset T$  defines a morphism of sSets
$\homm\Sets n F  \lra  \homm\Sets n T$.

\begin{exer}%_Exercise_.._   

  A filter $F$ converges on a topological space $T$ iff the morphism $\homm\Sets - F_\diag\lra T_\ttt$ factors as 
    $\homm\Sets - F_\diag \lra  T_\ttt[+1]\lra T_\ttt$. 
\end{exer}

\subsection{%## 
Compactness and completeness.
}\label{Compactness-as-lrl}
%....+++ todo: introductory remarks
We reformulate compactness and completeness in terms of iterated orthogonals/Quillen negation and 
morphisms representing typical examples of these notions. 
See the footnote in \S\ref{complete_as_negation} for the definition of the Quillen negation/orthogonals.

%\begin{defi}%Definition. 
%
%   Call X:sFilth  _complete_ iff for each filter {}\lra F^cart /_ X[+1]-->X
%   
%   Call X:sFilth _precompact_ iff for each filter F_const --> F^cart /_ X --> _|_
%   
%   Call X:sFilth  _compact_ iff for each ultrafilter {}\lra F^diag /_ X[+1]-->X
%   
%   Call X\lra Y in sFilth  _proper_  iff {}\lra F^diag /_ X[+1] \lra  X /\_X[+1] Y[+1]
%\end{defi}
%
%\begin{exer}%_Exercise_.._ 
%
%  A topological space $K$ is compact iff 
%
%    K.[+1]\lra K  (-  { K.[+1]-->K. : K a compact topological space }^lr 
% 
% A metric space M is complete iff 
%
%
%    M_\mmU[+1]\lra M  (-  { M_\mmU[+1]-->M_\mmU : M a complete metric space }^lr 
% 
%
%(Hint: the duals {..}^l contain the maps associated with filters, as above)
%
%Let  %{o<1}^cart |_| {o>1}^cart \lra  {o<->1}^cart

\subsubsection{Compactness and completeness as lifting properties/Quillen negation.}

\begin{exer}%Definition. 
\bi
\item A metric space $X$ is complete iff each Cauchy filter converges, i.e.~for each filter $F$ it holds 
$$\emptyset \lra \homm\Sets - F_\cart \rtt X_\mmU[+1]\lra X_\mmU$$

\item A topological space $X$ is quasi-compact iff each ultrafilter converges, i.e.~for each ultrafilter $\mathfrak U$ 
it holds $$\emptyset \lra \homm\Sets - {\mathfrak U}_\diag \rtt X_\ttt[+1]\lra X_\ttt$$

\item A metric space $X$ is compact iff each ultrafilter converges, i.e.~for each ultrafilter $\mathfrak U$ 
it holds \\ $$\emptyset \lra \homm\Sets - {\mathfrak U}_\diag \rtt X_\mmU[+1]\lra X_\mmU$$

\item A metric space $X$ is pre-compact iff each ultrafilter is Cauchy, i.e.~for each ultrafilter $\mathfrak U$ it holds \\ 
$$ \homm\Sets - {\mathfrak U}_\diag \lra \homm\Sets -{\mathfrak U}_\cart \rtt X_\mmU\lra \bot$$

\item A metric space $X$ is complete iff each Cauchy ultrafilter converges, i.e.~for each ultrafilter $\mathfrak U$ 
it holds \\ $$\emptyset \lra \homm\Sets -{\mathfrak U}_\cart \rtt X_\mmU[+1]\lra X_\mmU$$

\item  A map $X\lra Y$ of topological spaces is proper iff for each ultrafilter $\mathfrak U$ it holds 
%X.[+1] \lra  X. /\_X.[+1] Y.[+1] (- (  {o<1}^cart |_| {o>1}^cart \lra  {o<->1}^cart  )^lr
$$\emptyset \lra \homm\Sets -{\mathfrak U}_\diag \rtt
X_\ttt[+1] \lra  X_\ttt \vee_{X_\ttt[+1]} Y_\ttt[+1]
$$
This lifting property is equivalent to \href{http://mishap.sdf.org/tmp/Bourbaki_General_Topology.djvu#page=110110110110110110110110110110110}{[Bourbaki, I\S10,2,Theorem 1d]}:
%THEOREM I. Let j': X 
% Y be a continuous mapping. Then the following 
%four statements are equivalent: 
%a) f is proper. 
%b) f is closed and fl (y) is quasi-compact for each yeY. 
%c) If 
% is a filter on X and if y e Y is a cluster point of f (
%) then there 
%is a cluster point x of 
% such that f (x) = y. 
%d) If U is an ultrafilter on X and if y e Y is a limit point of the ultrafilter 
%base f (U), then there is a limit point x of U such that f (x) = y. 
\newline\noindent\includegraphics[width=1\linewidth]{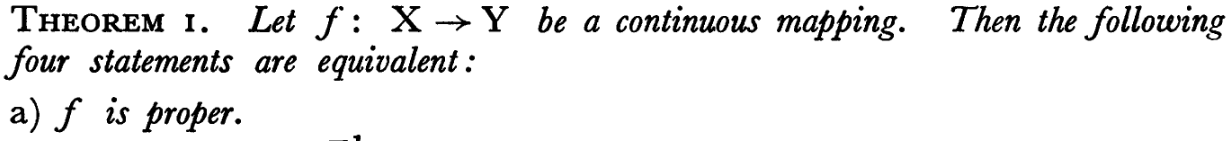}
\newline\noindent\includegraphics[width=1\linewidth]{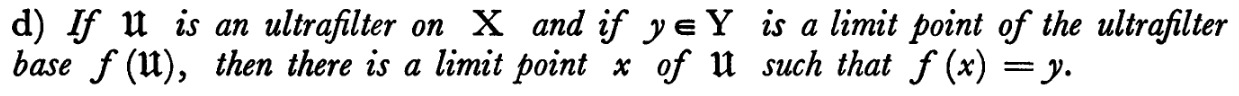}

\ei\end{exer}

\subsubsection{Concise reformulations in terms of iterated orthogonals/Quillen negations.}
The reformulations above lead to a concise expression in terms of iterated Quillen negations/orthogonals. 

Consider the two element set $\{o, 1 \}$ and let 
 $\{o\llrra 1 \}$ denote the filter with the unique neighbourhood on this set. 
Let $\{o<1\}$ and $ \{o>1\}$ denote the filters on this set generated by $\{o\}$, resp. $\{1\}$.

For a filter $F$, let $F_\cart$ denote the simplicial set $\homm{\Sets}- F_\cart$,
and let $\sqcup$ denote the disjoint union (equivalently, coproduct). 

Finally, let 
$\{o<1\}_\cart \sqcup \{o>1\}_\cart \lra  \{o\llrra 1 \}_\cart$ denote the obvious map. 

\begin{exer}%_Exercise_.._ 
\bi\item
   A topological space K is quasi-compact iff 
%  K.[+1]-->K. (- (  {o<1}^cart |_| {o>1}^cart \lra  {o<->1}^cart  )^lr 
$$ K_\ttt[+1]\lra K_\ttt \in  \left( \{o<1\}_\cart \sqcup \{o>1\}_\cart \lra  \{o\llrra 1 \}_\cart % {o<1}^cart |_| {o>1}^cart \lra  {o<->1}^cart  
\right)^\lr$$
\item

  A map $X\lra Y$ of topological spaces is proper iff
%X.[+1] \lra  X. /\_X.[+1] Y.[+1] (- (  {o<1}^cart |_| {o>1}^cart \lra  {o<->1}^cart  )^lr
$$
X_\ttt[+1] \lra  X_\ttt \vee_{X_\ttt[+1]} Y_\ttt[+1] \in 
 \left( \{o<1\}_\cart \sqcup \{o>1\}_\cart \lra  \{o\llrra 1 \}_\cart % {o<1}^cart |_| {o>1}^cart \lra  {o<->1}^cart  
\right)^\lr$$
\ei
%\begin{enonce}{Hint} 
(Hint: First check that a filter $F$ is an ultrafilter iff 
 %    {}-->F^diag /_  {o<1}^cart |_| {o>1}^cart \lra  {o<->1}^cart
$$\emptyset\lra F_\diag \rtt \{o<1\}_\cart \sqcup \{o>1\}_\cart \lra  \{o\llrra 1 \}_\cart $$
  (the lifting property means that the preimage of either $o$ or $1$ is a neighbourhood). 
  This is enough for the `if' implication.)%\end{enonce}
\end{exer}

\begin{exer}%_Exercise_.._(??) 
      Check whether the following holds for a reasonable class of metric spaces. 
      A metric space is complete iff $$ M_\mmU[+1]\lra M_\mmU \in  \{ \RR_\mmU[+1]\lra \RR_\mmU \}^{\lr}$$
      where $\RR$ denotes the real line with the usual metric. 
\end{exer}

%\begin{exer}%_Exercise_..(??) Find a map of "finite" (in what sense?) simplicial filters I1-->I2 such that 
%  a metric space is complete iff  $M_\mmU[+1]\lra M  \in  \{ I1-->I2 }^\lr $.
%\end{exer} 
%
%_Hint_: Try I_n := { the set of linear preorders on n elements } , 
%        with the finest filter such that I_0\lra I_n are continuous.
%
%        Then there is a map N^n\lra I_n , as each tuple (x1,..,xn) 
%        defines a linear preorder on n : i <= j iff x_i <= x_j. 
%        
%        The map appears to give rise to N^diag --> I. but not N^cart\lra I. .  
%        
%        not sure it works.
%
%

%Conjecture.  In the category of topological spaces,
%      (({0}-->{0->1})^r_{<5})^lr is the class of proper maps.

\begin{exer}[Taimanov theorem]%_Exercise_.._ 
Show that a map $f$ between $T_4$ topological spaces is proper iff in $Top$
       $f \in  \left((\{0\}\lra\{0\ra1\})^\rlr_{<5}\right)^\lr$. 
See  \href{http://mishap.sdf.org/mintsGE.pdf#page=11}{[Gavrilovich-Pimenov,\S2.2]} for explanations.
%%See  \href{http://mishap.sdf.org/mintsGE.pdf#page=67}{[ibid.,\S5,3]} 
%%e.g.~in $Top$  \bi\item 
%%$\{\bullet\}\lra X \in \{ \emptyset \lra \{\bullet\}\}^\text{rll}$ iff $X$ is connected and non-empty
%%\item $ \{ \emptyset \lra \{\bullet\}\}^\text{rr}$ is the class of subsets, i.e.~maps of form $X\subset Y$ (where the topology on $X$ is induced from $Y$)
%%\item  $ \{ \emptyset \lra \{\bullet\}\}^\text{rl}$ is the class of maps of form $X\lra X\sqcup D$ where $D$ is discrete
%%\item $  \{ \emptyset \lra \{\bullet\}\}^\text{r}$ is the class of surjections.
%%\ei
\end{exer}

\begin{exertodo} Check whether the following reformulations of topological properties hold in the category of topological spaces.
Can these expressions be interpreted in $\sFilth$ ? What would they mean for topological spaces ? 

For example, start by calculating the iterated orthogonals/negations below both in $Top$ and $\sFilth$ using the embedding of $Top$ into $\sFilth$.  
See  \href{http://mishap.sdf.org/mintsGE.pdf#page=67}{[Gavrilovich-Pimenov,\S5,3]} for the explnatation of the notation
and a list of other examples of  topological properties  expressed in terms of iterated Quillen orthogonals/negation.

\bi
%%
%%\item $  \{ \emptyset \lra \{\bullet\}\}^\text{r}$ is the class of surjections.
%%
%%
%%\item $ \{ \emptyset \lra \{\bullet\}\}^\text{rr}=$\verb| \{\{x\leftrightarrow  y\ra\}\lra\{x=y=c\}\}^\lrl|  is the class of subsets, i.e.~maps of form $X\subset Y$ (where the topology on $X$ is induced from $Y$)
%%\item  $ \{ \emptyset \lra \{\bullet\}\}^\text{rl}$ is the class of maps of form $X\lra X\sqcup D$ where $D$ is discrete
%%\item 
%%$\{\bullet\}\lra X \in \{ \emptyset \lra \{\bullet\}\}^\text{rll}$ iff $X$ is connected and non-empty
%%
%%\item $  X\lra\{\bullet\}\in \{ \emptyset \lra \{\bullet\}\}^\text{l}$ iff $X\neq\emptyset$
%%\item $ \{ \emptyset \lra \{\bullet\}\}^\text{ll}=(iso)$ is the class of isomorphisms
%%
%%\item $   \{ \emptyset \lra \{\bullet\}\}^\text{l} = \{ A\lra B : A\neq\emptyset\text{ or } A=B=\emptyset \,\}$ %iff $X\neq\emptyset$
%%\item $   X\lra\{\bullet\} \{ \emptyset \lra \{\bullet\}\}^\text{ll}$ iff $X=\emptyset$
%%
%%\item $  \{ \emptyset \lra \{\bullet\}\}^\text{lr}=  \{ \emptyset\lra Y :\ Y \text{ arbitrary} \}$
%%\item $  \{ \emptyset \lra \{\bullet\}\}^\text{lrr}$ is the class of maps which admit a section
%%\item $  \{ \emptyset \lra \{\bullet\}\}^\text{r}$
%%
\item \verb|{{}-->{o}}^r| is the class of surjections.
\item \verb|{{}-->{o}}^rr == {{x<->y->c}-->{x=y=c}}^l| is the class of subsets, i.e.~maps of form $X\subset Y$ (where the topology on $X$ is induced from $Y$)
\item \verb|{{}-->{o}}^rl|  is the class of maps of form $X\lra X\sqcup D$ where $D$ is discrete
\item $\{\bullet\}\lra X \in$ \verb|{{}-->{o}}^rll| iff $X$ is connected and non-empty
\item \verb|{{}-->{o}}^l| $= \{ A\lra B : A\neq\emptyset\text{ or } A=B=\emptyset \,\}$ 
\item \verb|{{}-->{o}}^lr| $ =\{ \emptyset\lra Y :\ Y \text{ arbitrary} \}$ 
\item \verb|{{}-->{o}}^lrr|  is the class of maps which admit a section
\item \verb|{{c}-->{o->c}}^l| is the class of maps with dense image
\item \verb|{{x,y}-->{x=y}}^r == {{x<->y}-->{x=y}^l| is the class of injective maps
\item \verb|{{x<->y->c}-->{x<->y=c}}^l == {{c}-->{o->c}}^lr| is the class of closed subsets
\item \verb| {{x<->y<-c}-->{x<->y=c}}^l| is the class of open subsets
\item $X\lra\{\bullet\}\in$ \verb|{{a<-b->c<-d->e}-->{b=c=d}}^l| iff $X$ is normal (T4)
\item (Tietze lemma) $\RR\lra\{\bullet\}\in$ \verb| {{a<-b->c<-d->e}-->{b=c=d},{a<-b->c}-->{a=b=c}}^lr|
\item in $Top$ the following expression means something similar to the Urysohn lemma without 
taking care of the necessary(!) conditions like being first countable:
 $$\verb| R-->{a<-b->c} (- {{a<-b->c<-d->e}-->{b=c=d}}^lr|$$ 
%%
%%
%%\begin{verbatim}
%%компактность:   {{o}-->{o->c}}^r_{<5}^lr   ' )<5() '\. 
%%плотный образ:  {{c}-->{o->c}}^l   .('\. 
%%сюрьекция:      {{}-->{o}}^r       .). 
%%инъкеция:       {{x,y}-->{x=y}}^r=={{x<->y}-->{x=y}^l   '~'(' == .-.). 
%%связность:      {{x,y}-->{x=y}}^l ..(.=.  {{}-->{o}}^rll    )((. 
%%дискретность:   {{}-->{o}}^rl             )(.  
%%подмножество: {{}-->{o}}^rr={{x<->y->c}-->{x=y=c}}^l  )). ==  ~\(.  
%%замкнутое подмножество:     {{x<->y->c}-->{x<->y=c}}^l=={{c}-->{o->c}}^lr 
%%открытое  подмножество:     {{x<->y<-c}-->{x<->y=c}}^l  '~'\ ( '~'=. 
%%нормальность (T4): {{a<-b->c<-d->e}-->{b=c=d}}^l    /V\(/\ 
%%лемма Титца: R-->{o} (- {{a<-b->c<-d->e}-->{b=c=d},{a<-b->c}-->{a=b=c}}^lr 
%%лемма Урысона: R-->{a<-b->c} (- {{a<-b->c<-d->e}-->{b=c=d}}^lr (не совсем!) 
%%ретракт: {{*-->{o}}^l    *(. 
%%окрестностой ретракт: {Y->-oo}-->{X->-oo} \in  {{*-->{o}}^l *(.  (почти!)
%%\end{verbatim}
%%
\ei\end{exertodo}

\begin{exertodo}\label{non-fork-extns}%_Exercise_.._ 
Being proper can also be defined as ``universally closed''. Formulate
an analogue of this definition in $\sFilt$. The following steps may be of use. 
\bi
\item Reformulate the condition that a map of topological spaces is closed in terms of neighbourhoods. 
Namely, a map $f:X\lra Y$ is closed iff %for every point $y_0 \in Y$, 
 every system $U_x\ni x$ of neighbourhoods there exist a system $V_y\ni y$ of neighbourhoods
such that $f^{-1}(V_y) \subset \bigcup\limits_{f(x)=y} U_x$, i.e.
\bi\item
$f(x') \in V_y$ implies that $x'\in U_x$ for some $x$ such that $f(x)=y$
\ei
\item Reformulate the tube lemma in a similar manner. 
\item Rewrite the above in terms of the simplicial neighbourhoods $f_\ttt:X_\ttt\lra Y_\ttt$.
\item Do the same for metric spaces. 
\item 
Ponder the syntactic similarity of the reformulations above to the characterisation of non-forking in terms of 
indiscernible sequences [Tent-Ziegler, 7.1.5] and to Definition~\ref{mapping_space} of neighbourhood structure of the mapping space.
\ei

\end{exertodo}

\begin{question} Find a compact proper definition of compact spaces and proper maps. Note that 
you probably want the following to be examples of compact spaces: 
(i) function spaces in \S\ref{arz-asc-spaces} coming from the Arzela-Ascoli theorems 
(ii) Stone spaces of indiscernible sequences in a  model in \S\ref{stone}.
\end{question}

\subsubsection{Convergence of sequences of functions}\label{ascoli}
Here we reformulate various notions of uniform convergence of a family of functions
as saying that a morphism in $\sFilth$ is well-defined.

Let $\{\NN\}$ denote the trivial filter on $\NN$ with a unique neighbourhood $\NN$ itself,
and $\NN_{cofinite}$ denote the filter of cofinite subsets of $\NN$.

A sequence  ${(f_i)_{i \in \NN}}$ of functions $f_i:X\lra M$ from a topological space $X$ to a metric space $M$
is 
 {\em equicontinuous} if either of the following equivalent conditions holds:
\bi
\item for every ${x \in X}$ and ${\varepsilon > 0}$,
there exists a neighbourhood ${U}$ of ${x}$ such that
${d_Y(f_i(x'), f_i(x)) \leq \varepsilon}$ for all ${i \in \NN}$ and ${x' \in U}$
\item the map $X_\ttt \times \{\NN\}_\const\lra M_\mU,\, (x,i)\longmapsto f_i(x)$ is well-defined
\item the map $X_\ttt \times (\NN_{cofinite})_\const\lra M_\mU,\, (x,i)\longmapsto f_i(x)$ is well-defined

\item the map $X_\ttt \times \{\NN\}_\diag  \lra M_\mU,\, (x,i)\longmapsto f_i(x)$ is well-defined

\item the map $X_\ttt \times (\NN_{cofinite})_\diag\lra M_\mU,\, (x,i)\longmapsto f_i(x)$ is well-defined

\ei

 If ${X = (X,d_X)}$ is also a metric space, replacing $X_\ttt$ by $X_\mmU$ above
gives us the notion of being {\em uniformly equicontinuous}. The family ${f_i}$
is {\em uniformly equicontinuous}
iff either of the following equivalent conditions holds:
\bi \item for every ${\varepsilon > 0}$ there exists a ${\delta > 0}$ such that
${d_Y(f_i(x'), f_i(x)) \leq \varepsilon}$ for all ${i \in \NN}$ and ${x', x \in X}$ with ${d_X(x,x') \leq \delta}$
\item the map $X_\mmU \times \{\NN\}_\const\lra M_\mU,\, (x,i)\longmapsto f_i(x)$ is well-defined
\item the map $X_\mmU \times (\NN_{cofinite})_\const\lra M_\mU,\, (x,i)\longmapsto f_i(x)$ is well-defined

\item the map $X_\mmU \times \{\NN\}_\diag  \lra M_\mU,\, (x,i)\longmapsto f_i(x)$ is well-defined

\item the map $X_\mmU \times (\NN_{cofinite})_\diag\lra M_\mU,\, (x,i)\longmapsto f_i(x)$ is well-defined
\ei

Replacing $\diag$ by $\cart$ gives us the notion of {\em uniformly Cauchy}. 
The family $f_i$ is {\em uniformly Cauchy} iff
either of the following equivalent conditions holds:
\bi \item for every ${\varepsilon > 0}$ there exists a ${\delta > 0}$ and $N>0$ such that
${d_Y(f_i(x'), f_j(x)) \leq \varepsilon}$ for all $i,j>N$ and ${x', x \in X}$ with ${d_X(x,x') \leq \delta}$.

\item the map $X_\mmU \times \{\NN\}_\cart  \lra M_\mU,\, (x,i)\longmapsto f_i(x)$ is well-defined

\item the map $X_\mmU \times (\NN_{cofinite})_\cart\lra M_\mU,\, (x,i)\longmapsto f_i(x)$ is well-defined
\ei

An equicontinuous family $f_i$  {\em converges} to a function $f_\infty$ iff
either of the following equivalent conditions holds:
\bi \item for every ${\varepsilon > 0}$ and $x\in X$ there exists a ${\delta > 0}$ and 
$N>0$ such that
${d_Y(f_\infty(x'), f_i(x')) \leq \varepsilon}$ for all $i>N$ and ${x \in X}$ with ${d_X(x,x') \leq \delta}$.

\item the following is a well-defined diagram:   
%$$X_\ttt[+1] \times (\NN_{cofinite})_\diag\lra M_\mU[+1] \,\,\,\,\, (x_0,i)\mapsto f(x_0)\,\,\,\,\,  (x,i)\longmapsto f_i(x)$$
%$$\text{ over  }X_\ttt \times (\NN_{cofinite})_\diag\lra M_\mU,\, (x,i)\longmapsto f_i(x)$$
%
$$ \xymatrix{  X_\ttt[+1] \times  (\NN_{cofinite} )_\diag   \ar@{..>}[r]|-------{} \ar@{->}[d]|-{[-1]\times \id} &   M_\mU[+1] \ar[d]|-{[-1]\times \id} \\  X_\ttt \times (\NN_{cofinite})_\diag  \ar[r]|---{(f_1,f_2,...)} &  M_\mU}$$
where the bottom morphism is  $(x,i)\longmapsto f_i(x)$, $i>0$,
and the top morphism is %necessarily of form 
$ (x_0,i)\longmapsto f_\infty(x_0)\,\,\,\,\,  (x,i)\longmapsto f_i(x)$, $i>0$.

Moreover, the top morphism is necessarily of this form for some function $f_\infty:X\lra M$. 
\ei

An uniformly equicontinuous family $f_i$  {\em uniformly converges} to a function $f_\infty$ iff
either of the following equivalent conditions holds:
\bi \item for every ${\varepsilon > 0}$ there exists %a ${\delta > 0}$ and 
$N>0$ such that
${d_Y(f(x), f_i(x)) \leq \varepsilon}$ for all $i>N$ and ${x \in X}$. %with ${d_X(x,x') \leq \delta}$.

\item the following is a well-defined diagram:%morphism  
%$$X_\mmU[+1] \times (\NN_{cofinite})_\diag\lra M_\mU[+1] \,\,\,\,\, (x_0,i)\mapsto f(x_0)\,\,\,\,\,  (x,i)\longmapsto f_i(x)$$
%$$\text{ over  }X_\mmU \times (\NN_{cofinite})_\diag\lra M_\mU,\, (x,i)\longmapsto f_i(x)$$
%
$$ \xymatrix{  X_\mU[+1] \times  (\NN_{cofinite} )_\diag   \ar@{..>}[r]|-------{} \ar@{->}[d]|-{[-1]\times \id} &   M_\mU[+1] \ar[d]|-{[-1]\times \id} \\  X_\mU \times (\NN_{cofinite})_\diag  \ar[r]|---{(f_1,f_2,...)} &  M_\mU}$$
where the bottom morphism is  $(x,i)\longmapsto f_i(x)$, $i>0$,
and the top morphism is %necessarily of form 
$ (x_0,i)\longmapsto f_\infty(x_0)\,\,\,\,\,  (x,i)\longmapsto f_i(x)$, $i>0$.

Moreover, the top morphism is necessarily of this form for some function $f_\infty:X\lra M$. 

\ei

A uniformly equicontinuous family $f_i$  has a subsequence which uniformly converges to a function $f$ iff
for some ultrafilter $\mathfrak U$ extending the filter $\NN_{cofinite}$ of cofinite subsets 
either of the following equivalent conditions holds:
\bi 
\item there is a sequence $(i_j)_{j\in \NN}$ such that 
for every ${\varepsilon > 0}$ there exists %a ${\delta > 0}$ and 
$N>0$ such that
${d_Y(f_\infty(x), f_i(x)) \leq \varepsilon}$ for all $i>N$ and ${x \in X}$. %with ${d_X(x,x') \leq \delta}$.
\item the following is a well-defined diagram for some  ultrafilter $\mathfrak U$ extending the filter $\NN_{cofinite}$ of cofinite subsets: 
%$$X_\mmU[+1] \times \mathfrak U_\diag\lra M_\mU[+1] \,\,\,\,\, (x_0,i)\mapsto f(x_0)\,\,\,\,\,  (x,i)\longmapsto f_i(x)$$
%$$\text{ over } X_\mU \times \mathfrak U_\diag\lra M_\mU,\, (x,i)\longmapsto f_i(x)$$
$$ \xymatrix{  X_\mU[+1] \times  \mathfrak U _\diag   \ar@{..>}[r]|-------{} \ar@{->}[d]|-{[-1]\times \id} &   M_\mU[+1] \ar[d]|-{[-1]\times \id} \\  X_\mU \times (\NN_{cofinite})_\diag  \ar[r]|---{(f_1,f_2,...)} &  M_\mU}$$
where the bottom morphism is  $(x,i)\longmapsto f_i(x)$, $i>0$,
and the top morphism is %necessarily of form 
$ (x_0,i)\longmapsto f_\infty(x_0)\,\,\,\,\,  (x,i)\longmapsto f_i(x)$, $i>0$.

Moreover, the top morphism is necessarily of this form for some function $f_\infty:X\lra M$.

\ei

\begin{exertodo} M.Dubashinsky suggested it might be possible to attempt to reformulate $\Gamma$-convergence in these terms.
``The natural setting of Γ$\Gamma$-convergence are lower semicontinuous functions'' [Braides, I,p.19],
and this suggests that a $\Gamma$-convergent sequence of functions $f_i:X\lra \RR$
something like a morphism $\NN_\leq \times X \lra \RR_\leq$ where $\NN_\leq:= \homm\leq-{\NN_\leq}$ and 
$\RR_\leq := \homm\leq-{\RR_\leq}$ are ssets of non-decreasing (or non-increasing?) sequences equipped 
with appropriate filters, or perhaps a simplex of an appropriate Skorokhod mapping space.  
\end{exertodo}

\subsubsection{Arzela-Ascoli theorems as diagram chasing} We now see that the Arzela-Ascoli theorem 
can be reformulated in terms of diagram chasing. The following exercise is a summary of the reformulations above.

\begin{exer} Check that the following diagrams represent the reformulation of the following Arzela-Ascoli theorem. 
\begin{quote} Theorem (Arzela-Ascoli).
Let $M$ be a complete metric space and $X$ be a compact metric space, and let
$f_i:X\lra Y$, $i\in \NN$, be a sequence of functions. 
%. Then every equicontinuous pointwise precompact 
%sequence  $f_i:X\lra Y$ of functions has a subsequence which converges uniformly, where pointwise precompact means
%that for each $x\in X$ the set $\{f_i(x)\}_i$ is precompact. 
Then the following are equivalent:
\bi\item[(i)] $(f_i)_{i \in \NN}$ has a convergent subsequence.

%is a precompact subset of {C(X \rightarrow Y)}.
\item[(ii)] ${(f_i)_{i \in \NN}}$ is pointwise precompact and equicontinuous.

\item[(iii)] ${(f_i)_{i \in \NN}}$ is pointwise precompact and uniformly equicontinuous.
\ei
\end{quote}

As diagram chasing:
 
\bi

\item ($X$ is compact) it holds $$ X_\mmU[+1]\lra X_\mmU \in 
 \left( \{o<1\}_\cart \sqcup \{o>1\}_\cart \lra  \{o\llrra 1 \}_\cart % {o<1}^cart |_| {o>1}^cart \lra  {o<->1}^cart  
\right)^\lr,$$
or, equivalently, for each ultrafilter $\mathfrak U$ it holds
  \\ $$\emptyset \lra \mathfrak U_\diag \rtt X_\mmU[+1]\lra X_\mmU$$

\item ($M$ is complete)  for each ultrafilter $\mathfrak U$ it holds \\ $$\emptyset \lra \mathfrak U_\cart \rtt M_\mmU[+1]\lra M_\mmU$$
or perhaps 
$$ M_\mmU[+1]\lra M_\mmU \in \left(  \RR_\mmU[+1]\lra \RR_\mmU \right)^\lr$$

\item (${(f_i)_{i \in \NN}}$ is pointwise precompact, i.e. for each point $x$ there is a subsequence such that $(f_{i_j}(x))_j$ converges)

 for each ultrafilter $\mathfrak U$ it holds
$$  \xymatrix{  X_\diag[+1] \times \mathfrak U_\diag   \ar@{..>}[r]|--{} \ar@{->}[d]|{} &   M_\mU[+1] \ar[d]|{} \\  X_\ttt \times (\NN_{cofinite})_\diag  \ar[r]|-{(f_i)} &  M_\mU}$$
 where  $X_\diag=(X_{\text{discrete}})_\mU$ denotes $X$ equipped with the filter of diagonals, i.e. a subset is a neighbourhood iff it contains the diagonal.

\item (${(f_i)_{i \in \NN}}$ being uniformly equicontinuous imply they converge uniformly)

 for each ultrafilter $\mathfrak U$ it holds
 $$\xymatrix{  X_\mmU[+1] \times \mathfrak U_\diag   \ar@{..>}[r]|--{} \ar@{->}[d] &   M_\mU[+1] \ar[d]|{} \\  X_\mU \times (\NN_{cofinite})_\diag  \ar[r]|-{(f_i)} &  M_\mU}$$

\item (${(f_i)_{i \in \NN}}$ being equicontinuous imply being uniformly equicontinuous)

 $$\xymatrix{  X_\ttt \times (\NN_{cofinite})_\diag   \ar[r]|--{} \ar@{->}[d] &   M_\mU \ar[d]|{} \\ 
  X_\mU \times (\NN_{cofinite})_\diag  \ar@{..>}[r]|-{} &  M_\mU}$$

\ei
%%
%%
%%\item[(ii)$\Longleftrightarrow$(iii)] 
%%$$  \xymatrix{  X_\diag \times \mathfrak U_\diag   \ar@{..>}[r]|--{} \ar@{->}[d]|{} &   M_\mU[+1] \ar[d]|{} \\  X_\ttt \times (\NN_{cofinite})_\diag  \ar[r]|-{(f_i)} &  M_\mU}
%%\text{ iff }
%% \xymatrix{  X_\mmU \times \mathfrak U_\diag   \ar@{..>}[r]|--{} \ar@{->}[d] &   M_\mU \ar[d]|{} \\  X_\ttt \times (\NN_{cofinite})_\diag  \ar[r]|---{(f_i)} &  M_\mU}$$
%%
%%\item[(iii)$\implies$(i)]
%%$$  \xymatrix{  \{x\}\times \mathfrak U_\diag   \ar@{..>}[r]|--{} \ar@{->}[d]|{\forall x\in X} &   M_\mU[+1] \ar[d]|{} \\  X_\mU \times (\NN_{cofinite})_\diag  \ar[r]|-{(f_i)} &  M_\mU}
%%\text{ implies }
%% \xymatrix{  X_\mmU[+1] \times \mathfrak U_\diag   \ar@{..>}[r]|--{} \ar@{->}[d] &   M_\mU[+1] \ar[d]|{} \\  X_\mU \times (\NN_{cofinite})_\diag  \ar[r]|-{(f_i)} &  M_\mU}$$
%%
%%\ei
\end{exer}

\subsubsection{Arzela-Ascoli theorems; compactness of function spaces}\label{arz-asc-spaces}

Above we saw that to give a convergent sequence of functions is the same as to give a certain morphism in $\sFilth$,
and that to take the limit of the sequence is to ``lift'' this morphism along shifts $X[+1]\xra{[-1]} X$ and  $M[+1]\xra{[-1]} M$.  
This makes various Arzela-Ascoli theorems on pre-compactness of function spaces reminiscent of base change. 

\begin{todo} Give a category theoretic approach to various Arzela-Ascoli theorems
including the Prokhorov theorem on tightness of measures and $\Gamma$-convergence.
\end{todo}

The following exercises may be of use. 

\begin{exertodo}
Calculate the Skorokhod mapping spaces related to the morphisms above. Are these spaces pre-compact? complete? 
Is the subspace of continuous functions relatively compact within the mapping space under reasonable assumptions, 
i.e.~is the subset $$\homm\sFilt X Y \subset \Homm\sFilt X Y$$ relatively compact, and in what precise meaning ? 

Calculate the following Skorokhod function spaces:
\bi\item 
 equicontinuous:
\bi
\item $\Homm\sFilth{X_\ttt \times \{\NN\}_\const}{M_\mU}$ 
\item $\Homm\sFilth{X_\ttt \times (\NN_{cofinite})_\const}{M_\mU}$%,\ (x,i)\longmapsto f_i(x)$ is well-defined

\item $\Homm\sFilth{X_\ttt \times \{\NN\}_\diag  }{M_\mU}$%,\ (x,i)\longmapsto f_i(x)$ is well-defined

\item  $\Homm\sFilth{X_\ttt \times (\NN_{cofinite})_\diag}{M_\mU}$% ,\ (x,i)\longmapsto f_i(x)$ is well-defined

\ei
\item 
% If ${X = (X,d_X)}$ is also a metric space, replacing $X_\ttt$ by $X_\mmU$ above
%gives us the notion of being {\em 
uniformly equicontinuous:
%}. The family ${f_i}$
%is {\em uniformly equicontinuous}
%iff either of the following equivalent conditions holds:
\bi 
\item $\Homm\sFilth{X_\mmU \times \{\NN\}_\const}{M_\mU}$ % (x,i)\longmapsto f_i(x)$ is well-defined
\item $\Homm\sFilth{X_\mmU \times (\NN_{cofinite})_\const}{M_\mU}$ % (x,i)\longmapsto f_i(x)$ is well-defined

\item $\Homm\sFilth{X_\mmU \times \{\NN\}_\diag  }{M_\mU}$ % (x,i)\longmapsto f_i(x)$ is well-defined

\item $\Homm\sFilth{X_\mmU \times (\NN_{cofinite})_\diag}{M_\mU}$ % (x,i)\longmapsto f_i(x)$ is well-defined
\ei
\item uniformly Cauchy:
\bi 
\item $\Homm\sFilth{X_\mmU \times \{\NN\}_\cart  }{M_\mU}$ % (x,i)\longmapsto f_i(x)$ is well-defined

\item $\Homm\sFilth{X_\mmU \times (\NN_{cofinite})_\cart}{M_\mU}$ % (x,i)\longmapsto f_i(x)$ is well-defined
\ei
\ei

\end{exertodo}

\begin{exertodo}
In measure theory Prokhorov's theorem relates tightness of measures to relative compactness (and hence weak convergence) in the space of probability measures. Reformulate the Prokhorov's theorem and give a uniform approach to both Prokhorov theorem and Arzela-Ascoli theorems.
\end{exertodo}

\begin{exertodo} Rewrite the Arzela-Ascoli theorems entirely in terms of iterated orthogonals/Quillen negations, 
taking M2-decompositions (i.e.~(co)fibrant replacement), and, more generally, rules for manipulating labels on morphisms.
%%.
%%For example, something along the following lines might work.
%%
%%Let $$(Comp) := \left(\{o<1\}_\cart \sqcup \{o>1\}_\cart \lra  \{o\llrra 1 \}_\cart\right)  $$
%%$$ (pre-Comp) := \left\{\mathfrak  U_\diag  \lra \mathfrak U_\cart : \ \mathfrak U:\text{is an ultrafilter } \right\}$$
%%
%%
%%
%%
%%
%%Then Arzela-Ascoli can reformulated as the following diagram chasing rule:
%%\begin{equation*}\begin{split}
%%  K_\ttt[+1]\xra{(Comp)^\lr} K_\ttt \text{ and } 
%%  M_\mU[+1]\xra{(Comp)^\lr} M_\mU  
%%\implies  \\ 
%% \forall\, \mathfrak U:\text{ultrafilter }  
%% \bot\lra K_\ttt\times  \mathfrak U_\diag \rtt & M_\mU[+1] \lra  M_\mU 
%%\end{split}\end{equation*}
%%
%%
%%
%%\begin{equation*}\begin{split}
%% K_\ttt[+1]\xra{(Comp)^\lr} K_\ttt
%%%\item[] %($M$ is pre-compact) 
%%\text{ and }  M_\mU\xra{(pre-Comp)^r}\top  
%%%\ for each ultrafilter $\mathfrak U$ 
%%\implies  \\ 
%% \forall\, \mathfrak U:\text{ultrafilter }  
%% K_\ttt\times  \mathfrak U_\diag \lra  & K_\ttt\times  \mathfrak U_\cart\, \rtt\,  M_\mU\lra \top
%%\end{split}\end{equation*}
%%
%%
%%Can the  latter be reformulated further as: 
%%$$ f \in (Comp)^\lr \text{ and } g \in {(pre-Comp)^\rl} \implies f\times g \in {(pre-Comp)^\rl} \ \ ?$$
%%
\end{exertodo}

\section{Model theory}\label{sec:4}

\subsection{Ramsey theory}\label{ramsey} 
Let $X:\sSet$ be a simplicial set, and $c:X^\text{nd}_n \lra C$ be a colouring of the set $X^\text{nd}_n$ 
of non-degenerate $n$-simplices, i.e.~an arbitrary function defined on the set of non-degenerate simplices of $X$ of dimension $n$. 
Call a simplex $s:X_N$ {$c$-homogeneous} iff all its non-degenerate faces of dimension $n$ have the same $c$-colour.  
Let $c(X)$ be the subsset of $X$ consisting of $c$-homogeneous simplices in $X$. 
\begin{exer}[Ramsey theorem] 
\bi
\item Verify that $c(X):\sSet$ is indeed a well-defined simplicial set; it is a disjoint union of subssets corresponding to different $c$-colours.  
\item (Ramsey theorem) Let the set of $c$-colours be finite. 
If $X:\sSet$ has non-degenerate simplices of arbitrarily high dimension, 
then so does $c(X)$. In more detail, if for unboundedly many $N$ there is a simplex in $X_N$ 
which is not a face of a simplex of lower dimension, then the same holds for $c(X)$. 
\item Take $X:\sSet$ to be $\homm\Sets-S$ where $S$ is an infinite set. 
Verify that the item above is the usual statement of Ramsey theorem: 
for each colouring of subsets of  $S$ of size $n$, there is an arbitrarily large subset of $S$ 
such that all its subsets of size $n$ have the same colour. 
\ei 
\end{exer} 

\begin{todo} Ponder the discussion of Ramsey theory-type theorems in [MLGrovov], [Gromov2014].
\end{todo}

\subsection{Indiscernible sequences in model theory}\label{model} Ramsey theory provides a basic tool in model theory known as the {\em indiscernible sequences}. 

\begin{defi}[$L$-indiscernability neighbourhood structure on a model $M^I$] Let $M$ be a model in a language $L$, and let $I_\leq$ be a linear order. 
For a $n$-ary formula $\varphi$ of $L$, we say that a sequence $(a_i)\in M^I$ is 
{\em $\varphi$-homogeneous} iff 
\bi
\item[]
 $ M \models \varphi(a_{i_1},...,a_{i_n}) \equiv  \varphi(a_{j_1},...,a_{j_n})$ whenever $i_1<i_2<...<i_n$, $j_1<j_2<...<j_n$, and 
all the $a_{i_k}$'s are distinct, and all the  $a_{j_k}$'s are distinct, 
i.e.~ $a_{i_k}\neq a_{i_l}$ and $a_{j_k}\neq a_{j_l}$ whenever $1\leq k<l\leq n$ %$a_{i_k}\neq a_{i_l}$ and $a_{j_k}\neq a_{j_l}$
\ei
For a type $\pi$, a sequence is {\em $\pi$-homogeneous} iff it is $\varphi$-homogeneous for each formula in $\pi$.
A {\em $\varphi$-neighbourhood of the diagonal in $M^I$} is the subset consisting of all the $\varphi$-homogeneous sequences. 
The {\em $L$-indiscernability filter on $M^I$} is generated by the $\varphi$-neighbourhoods of the diagonal, for $\varphi$ a formula of $L$. 
\end{defi}

The requirement that all the $a_{i_k}$'s are distinct, and all the  $a_{j_k}$'s are distinct, is needed to show
 that a non-decreasing map $f: J\lra I$ induces a continuous map  $f^*:M^I \lra M^J$ of indiscernability filters.

\def\Stone{{\ethi{\ethmath{gE}}}}
\def\Stone{{\ethi{\ethmath{rI}}}}

\begin{defi}[The Stone space of a model.]\label{stone} 
Let $M$ be a model in a language $L$, and let $A\subset M$ be a subset. Also assume that $M$ is $\card(A)^+$-saturated.  
Call the {\em Stone space $\Stone^M_1(A)=\Stone_1^M(A)$ of a model $M$ with parameters $A$} the sset $n_\leq \longmapsto \homm\Sets n M$ where 
$M^n= \homm\Sets n M $ is equipped with the $L(A)$-indiscernability neighbourhood structure. 
Similarly define  $\Stone_n^M(A)$ as the Stone space of $n$-tuples. 
\end{defi}

The following exercise is based on the syntactic similarity of the characterisation of non-forking in terms
of indiscernible sequences  [Tent-Ziegler, Lemma 7.1.5], the reformulation Exercise~\ref{non-fork-extns} of 'being a closed map' in terms of neighbourhoods,
and of  Definition~\ref{mapping_space} of neighbourhood structure of the mapping space.

\begin{exer}
Check that the usual Stone space $S_1(A)$ is the Hausdorff quotient of the topologisation of $\sFilt$-Stone space $\Stone^M_1(A)$
whenever $M$ is $\card(A)^+$-saturated.
\end{exer}

\begin{todo} 
\bi\item Is $\Stone^M_n(A)$ quasi-compact? complete? 
\item Is the map $\Stone^M_n(B)\lra \Stone^M_n(A)$ proper, cf.~Exercise~\ref{non-fork-extns} ? 
Is this related to the fact that every type has a non-forking extension? 
\item Interpret %the characterisation of non-forking in terms of indiscernible sequeces  
[Tent-Ziegler, Lemma 7.1.5] 
as a property of neighbourhoods of Stone spaces, e.g.~that a projection of a certain subset is closed, cf.~Exercise~\ref{non-fork-extns}.
\item Are there are non-trivial maps $\Stone^M_n(A)\lra \Stone^N_n(B)$ for different models $M$ and $N$; what does it mean model-theoretically?

\item More generally, ponder if $\sFilt$-Stone spaces allow to define interesting homological or homotopical invariants of models. 
\ei
\end{todo}

\begin{todo}
\bi\item Check whether a model $M$ is stable iff $\Stone^M_n(A)$ is symmetric
for every $A\subset M$, i.e.~$\Stone^M_n(A):\Dop\lra\Filt$ factors as
$\Stone^M_n(A): \Dop \lra FiniteSets^\text{op}\lra\Filt$
\item Are there similar characterisations of e.g.~simple or NIP theories in terms of their Stone spaces? 
\ei\end{todo}

\begin{todo}
Do these Stone spaces allow to express neatly the theory of forking? 
\bi \item For example, do pushforwards mentioned in [Simon, Exercise 9.12 (distality of $T^{eq}$)]
 are indeed pushforwards in $\sFilt$?
\item Reformulate the definition [Simon, Def.~9.28] of distality in terms of %indiscernible sequences %[Simon, Def.~9.28]
endomorphisms of $\sFilt$-Stone spaces or similar objects of $\sFilt$.
\newline\noindent\includegraphics[width=1\linewidth]{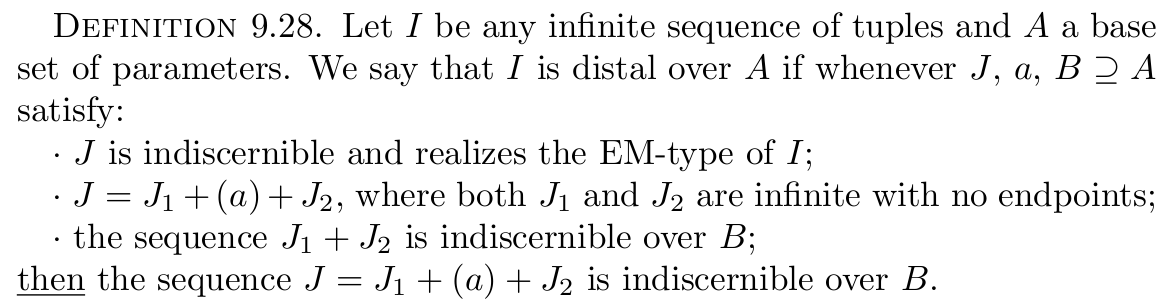} 
\ei

\end{todo}

\tiny
\subsubsection*{Acknowledgements} I thank M.Bays, K.Pimenov for many useful discussions, proofreading and generally 
taking interest in this work. Thanks are also due to  S.Ivanov, S.Kryzhevich, S.Sinchuk, V.Sosnilo, and A.Smirnov for a number of discussions,
and particularly to the participants of the seminar of A.Smirnov. I thank D.Rudskii for bringing the Levy-Prohorov metric to my attention. 
I thank D.Grayson for several insightful comments.
I express deep gratitude to
friends for creating an excellent social environment in St.Petersburg, invitations and enabling me to work. % there. 

This %work 
continues [Gavrilovich-Pimenov]; see also acknowledgements there.

%\theendnotes

\end{document}